\documentclass[11pt]{article}
\usepackage[margin=1.05in]{geometry}
\usepackage{amsmath,amssymb,amsthm}
\usepackage{booktabs}
\usepackage{enumitem}
\usepackage[colorlinks=true,linkcolor=blue,citecolor=blue,urlcolor=blue]{hyperref}

\newtheorem{theorem}{Theorem}[section]
\newtheorem{lemma}[theorem]{Lemma}
\newtheorem{proposition}[theorem]{Proposition}
\newtheorem{corollary}[theorem]{Corollary}
\theoremstyle{definition}
\newtheorem{definition}[theorem]{Definition}
\newtheorem{remark}[theorem]{Remark}
\newtheorem{problem}[theorem]{Open Problem}
\newtheorem{algorithm}[theorem]{Algorithm}

\newcommand{\Z}{\mathbb{Z}}
\newcommand{\R}{\mathbb{R}}
\newcommand{\E}{\mathbb{E}}
\newcommand{\Prob}{\mathbb{P}}
\newcommand{\dtv}{d_{\mathrm{TV}}}
\newcommand{\dw}{d_{\mathrm{W}}}
\newcommand{\Bes}{\mathrm{Bes}}
\newcommand{\Poi}{\mathrm{Poi}}
\newcommand{\Mult}{\mathrm{Mult}}
\newcommand{\Sk}{\mathsf{Sk}}
\newcommand{\Ct}{\widehat C}
\newcommand{\Lm}{\Lambda_{k,\nu}}
\newcommand{\Gm}{\Gamma_{k,\nu}}
\newcommand{\Wsq}{W^{\square}}
\newcommand{\Zh}{\widehat Z}
\newcommand{\Zt}{\widetilde Z}
\newcommand{\Zb}{\overline Z}

\title{\textbf{Two moments suffice for Bessel approximation}\\[3pt]
\large Conditioned differences, contingency tables,\\ and probabilistic
divide-and-conquer}
\author{Stephen DeSalvo\thanks{%
Independent researcher.  Email: \texttt{stephendesalvo@gmail.com}.}}
\date{October 2026}

\begin{document}
\maketitle

\begin{abstract}
We develop Stein's method, with explicit constants throughout, for the
\emph{discrete Bessel distribution}: the law of a Poisson variable given
its difference from an independent Poisson variable.  The Bessel family is
characterized by a birth--death operator with quadratic death rate
$n(n{+}k)$ and as the exact fixed-point set of a centered quadratic-bias
transform, and we prove a dependency expansion of
Arratia--Goldstein--Gordon type: two joint moments of a column array
control the total-variation distance from the conditioned row sum to a
matched Bessel law.  Applied to the row-constrained Fisher--Yates ensemble
on $2\times c$ contingency tables --- independent Bernoulli or
capped-Poisson entries conditioned on the row-sum difference, column sums
free --- this yields a limit theorem --- which upgrades to the entire
conditioned table at the same two-moment error --- a completely explicit
enumeration estimate, and an exact divide-and-conquer sampler with
certified acceptance.  The entire apparatus --- operator, expansion, enumeration,
sampler --- generalizes to $r\times c$ tables, where conditioning on all
successive row differences produces a \emph{hyper-Bessel} law.
Verification scripts accompany the arXiv version.
\end{abstract}

\noindent\textit{MSC2020 subject classifications:} Primary 60F05, 62E17;
secondary 05A16, 60C05, 62H17, 65C05.

{\small\tableofcontents}

\bigskip
\begin{center}\rule{0.5\textwidth}{0.4pt}\\[4pt]
\textsc{Part I. Narrative}\end{center}

\section{Introduction: from Poisson to Bessel}\label{sec:intro}

\subsection*{The Poisson template}
Poisson approximation by Stein's method \cite{chen1975,agg1989,bhj1992} turns
on the characterization $\E[\lambda f(W{+}1)-Wf(W)]=0$ and on the
\emph{size-bias transform} $\Prob(W^*{=}n)=n\Prob(W{=}n)/\E W$: the
Poisson laws are exactly the fixed points of the centered map
$W\mapsto W^*{-}1$ (equivalently, $W^*\overset{d}{=}W{+}1$ iff $W$ is
Poisson), and for a sum
$W=\sum_i\mathbf 1_{A_i}$ of weakly dependent indicators,
\begin{equation}\label{eq:aggtemplate}
\dtv\bigl(W,\Poi(\lambda)\bigr)\le
\bigl(1\wedge\lambda^{-1}\bigr)(b_1+b_2)
+\bigl(1\wedge1.4\,\lambda^{-1/2}\bigr)\,b_3\,,
\end{equation}
with
\[
b_1=\sum_{i}\sum_{j\in N_i}\Prob(A_i)\Prob(A_j),\qquad
b_2=\sum_i\sum_{j\in N_i\setminus\{i\}}\Prob(A_iA_j),\qquad
b_3=\sum_i\E\Bigl|\,\E\bigl[\mathbf 1_{A_i}-\Prob(A_i)\;\big|\;
\sigma(\mathbf 1_{A_j}:j\notin N_i)\bigr]\Bigr|\,:
\]
$b_1$ and $b_2$ are neighborhood sums of first and second moments over a
dependency structure $\{N_i\}$, and $b_3$ measures the dependence of each
indicator on everything \emph{outside} its neighborhood, vanishing when the
neighborhoods capture all the dependence.  Two moments suffice: the
parameter is $\lambda=\E W$, and the error is local.

\subsection*{The Bessel analogue}
Let $X\sim\Poi(\lambda_1)$ and $Y\sim\Poi(\lambda_2)$ be independent, and
condition on $\{X-Y=k\}$.  Then $Y$ follows the \emph{discrete Bessel law} \cite{yuankalbfleisch2000}
$\Bes(k,a)$, $a=2\sqrt{\lambda_1\lambda_2}$, with mass function
\begin{equation}\label{eq:pmf}
p_n=\frac1{I_k(a)}\,\frac{(a/2)^{2n+k}}{n!\,(n+k)!},\qquad n\ge0,
\end{equation}
where $I_k(a)=\sum_{n\ge0}(a/2)^{2n+k}/(n!(n+k)!)$ is the modified Bessel
function of the first kind --- the normalizing constant and the source of the
name.  This paper constructs, around \eqref{eq:pmf}, the exact counterpart of
each ingredient of the Poisson template, and applies it to contingency
tables.  The correspondence that organizes everything is:

\begin{center}
\renewcommand{\arraystretch}{1.3}
\begin{tabular}{@{}ll@{}}
\toprule
\textbf{Poisson} & \textbf{Bessel} \\
\midrule
operator $\ \lambda f(n{+}1)-n\,f(n)$ &
operator $\ \nu\,g(n{+}1)-n(n{+}k)\,g(n)$,\quad $\nu=a^2/4$ \\
death rate $n$ (linear) & death rate $n(n{+}k)$ (quadratic) \\
moment identity $\ \E W=\lambda$ &
moment identity $\ \E[W(W{+}k)]=\nu$ \\
size-bias $\ n\,p_n/\E W$ & quadratic-bias $\ n(n{+}k)\,p_n/\nu$ \\
fixed points of $W{\mapsto}W^*{-}1$: the Poisson family &
fixed points of $\mathcal T_k:W{\mapsto}\Wsq{-}1$: the Bessel family \\
Stein factor $\ (1{-}e^{-\lambda})/\lambda$ &
Stein factor $\ \Lm\le\min\{\tfrac1\nu{+}\tfrac1{k+1},\tfrac3\nu\}$ \\
coupling of $W,W^*$: add one point &
coupling of $W,\Wsq$: add a \emph{pair} (fiber-preserving) \\
error terms $\ b_1,b_2$ &
product-bias $\ b_1,b_2$ + fiber constants $\theta,\tau_1,\omega$ \\
\bottomrule
\end{tabular}
\end{center}

\subsection*{What is proved, and where}
Part~I states everything and explains the ideas; Part~II proves everything in
full, one subsection of proofs per theme.  Section~\ref{sec:operator} sets up
the operator, the closed-form Stein solution, and the Stein factors, all
resting on one structural input: log-concavity of \eqref{eq:pmf}.
Section~\ref{sec:agg} is the core --- the quadratic-bias fixed point, the
approximation theorem, and (\S\ref{sec:coupling}) the coupling that
carries the method.  Section~\ref{sec:expansion} runs the
Arratia--Goldstein--Gordon expansion for column arrays, where the coupling
becomes a fiber-preserving paired modification and the bound becomes
neighborhood-local.  Section~\ref{sec:tables} is the contingency-table payoff, in three
registers: total variation for the conditioned row sum (\S\ref{sec:dtv}),
a completely explicit enumeration estimate (\S\ref{sec:enum}), and an
exact divide-and-conquer sampler with certified acceptance and a bound on
block splits (\S\ref{sec:pdc}).
Section~\ref{sec:mrow} extends the entire apparatus to general $r\times c$
tables, where the target becomes a hyper-Bessel law.
Section~\ref{sec:numerics} collects the numerics and \S\ref{sec:open} the
open problems; \S\ref{sec:final} closes Part~I with remarks situating each
contribution against the prior frameworks that partially contain it.

\subsection*{Related work}
The distribution \eqref{eq:pmf} was introduced by Yuan--Kalbfleisch
\cite{yuankalbfleisch2000}; its role as a conditioned Poisson difference
goes back to Skellam \cite{skellam1946}, and Devroye \cite{devroye2002}
gave exact samplers.  Our starting points sit inside three existing
frameworks.  The operator \eqref{eq:operator} is a constant-birth,
quadratic-death instance of the polynomial birth--death and discrete-Gibbs
frameworks \cite{brownxia2001,er2008}, and its $k=0$ case appears in the
conditional-approximation preprint of Chiarini--Cipriani--Conforti
\cite{ccc2017}.  The quadratic-bias transform of Definition~\ref{def:qbias}
is an instance of the polynomial-bias paradigm
\cite{goldsteinreinert2005,dobler2017}, its centered form joining the
fixed-point family of zero-bias (normal) \cite{goldsteinreinert1997},
size-bias-minus-one (Poisson), and equilibrium (exponential)
\cite{pekozrollin2011}.  The expansion of \S\ref{sec:expansion} follows the
two-moments lineage \cite{agg1989,bcl1992,brownphillips1999,rollin2007},
treating the \emph{conditional} law of the difference where Gan--Kolaczyk
\cite{gan2018} treat the unconditioned Skellam law.  A precise accounting
of what each framework already delivers here and what does not follow from
it --- including places where the general theory beats our closed
constants --- is collected in the final remarks of \S\ref{sec:final}.

\subsection*{Notation}
$\Z_{\ge0}=\{0,1,2,\dots\}$ and $[c]=\{1,\dots,c\}$; $x\wedge y=\min\{x,y\}$;
$\lfloor\cdot\rfloor$ is the floor; $\mathbf 1_A$ is an indicator; $:=$
introduces a definition.  $\E$, $\Prob$, and $\mathrm{Var}$ denote
expectation, probability, and variance, and $\sigma_X=\sqrt{\mathrm{Var}\,X}$;
$X\perp Y$ is independence; $X\overset{d}{=}Y$ is equality in law; $X\sim\mu$
says $X$ has law $\mu$.  $\Poi$, $\mathrm{Bern}$, and $\Mult$ are the
Poisson, Bernoulli, and multinomial families, and $I_k$ is the modified
Bessel function of the first kind.  For a pmf $(p_n)$ on $\Z_{\ge0}$,
$F(m)=\sum_{n\le m}p_n$ and $G(m)=\sum_{n\ge m}p_n$ are the cdf and
survival function.  For $\Z$-valued $X,Y$ (or their laws),
$\dtv(X,Y)=\sup_{A\subseteq\Z}\bigl|\Prob(X\in A)-\Prob(Y\in A)\bigr|$ and
$\dw(X,Y)=\inf\E|X'-Y'|$, the infimum over couplings $(X',Y')$, are the
total-variation and $\ell^1$-Wasserstein distances.  For $g:\Z_{\ge0}\to\R$,
$\Delta g(m)=g(m{+}1)-g(m)$ and $\|g\|=\sup_m|g(m)|$.

Asymptotic and comparison symbols appear only in narrative commentary,
never as hypotheses of a theorem --- every theorem carries explicit
constants.  $O(\cdot)$ and $o(\cdot)$ are the Landau symbols in the stated
limit; $f\asymp g$ means $f=O(g)$ and $g=O(f)$; $a\gg b$ means $b/a\to0$ in
the ambient limit; $\approx$ marks numerical approximate equality in
examples and tables; $\rightsquigarrow$ indicates notational substitution
when a statement is transported to a new setting.

\section{The Bessel distribution and its Stein operator}\label{sec:operator}

Throughout, $k\in\Z_{\ge0}$, $a>0$, and $\nu:=a^2/4$.  Write
$F(m)=\Prob(N\le m)$ and $G(m)=\Prob(N\ge m)$ for $N\sim\Bes(k,a)$.  The one
structural fact from which everything else in this section follows is the
ratio identity
\begin{equation}\label{eq:ratio}
\frac{p_{n+1}}{p_n}=\frac{\nu}{(n+1)(n+k+1)},\qquad n\ge0,
\end{equation}
which is strictly decreasing in $n$; equivalently $n(n{+}k)\,p_n=\nu\,p_{n-1}$.

\begin{lemma}[Basic properties]\label{lem:basic}
For $k\ge0$ and $a>0$, the law $\Bes(k,a)$ of \eqref{eq:pmf}, with mass
function $(p_n)$ and $\nu=a^2/4$, is a probability law on $\Z_{\ge0}$,
strictly log-concave
\textup{(}$p_{n+1}p_{n-1}<p_n^2$, $n\ge1$\textup{)}, with mode
$m^*=\lfloor(\sqrt{k^2+4\nu}-k)/2\rfloor$, the largest $n$ with
$n(n{+}k)\le\nu$.  All moments are finite, and, for $N\sim\Bes(k,a)$,
\begin{equation}\label{eq:moments}
\E N=\frac a2\,\frac{I_{k+1}(a)}{I_k(a)},\qquad
\E\bigl[N(N{+}k)\bigr]=\nu\quad\text{exactly}.
\end{equation}
\end{lemma}

\begin{lemma}[Conditional representation]\label{lem:condrep}
If $X\sim\Poi(\lambda_1)\perp Y\sim\Poi(\lambda_2)$ and $k\ge0$, then
\begin{equation}\label{eq:skellampt}
\Prob(X-Y=k)=e^{-\lambda_1-\lambda_2}\Bigl(\tfrac{\lambda_1}{\lambda_2}\Bigr)^{k/2}I_k(a),
\qquad a=2\sqrt{\lambda_1\lambda_2},
\end{equation}
and $\mathcal L(Y\mid X-Y=k)=\Bes(k,a)$, the law \eqref{eq:pmf}; the
conditional law depends on
$(\lambda_1,\lambda_2)$ only through $\lambda_1\lambda_2$.
\textup{(}For $k<0$, condition on $\{Y-X=|k|\}$ and swap; $k\ge0$ is no
loss.\textup{)}
\end{lemma}

The identity $\E[N(N{+}k)]=\nu$ is the exact analogue of
$\E[\Poi(\lambda)]=\lambda$; it names the ``two moments'' the method will
match.  For bounded $g:\Z_{\ge1}\to\R$ (extended by $g(0):=0$, never
evaluated since $n(n{+}k)g(n)$ vanishes at $n=0$) define the \emph{Bessel
Stein operator}
\begin{equation}\label{eq:operator}
(\mathcal Ag)(n):=\nu\,g(n{+}1)-n(n{+}k)\,g(n),\qquad n\ge0,
\end{equation}
the density-form operator of the birth--death chain with constant birth rate
$\nu$, quadratic death rate $n(n{+}k)$, and stationary law \eqref{eq:pmf}.

\begin{theorem}[Characterization and Stein solution]\label{thm:char}
A $\Z_{\ge0}$-valued $W$ with $\E W^2<\infty$ satisfies $W\sim\Bes(k,a)$ iff
$\E[\nu g(W{+}1)-W(W{+}k)g(W)]=0$ for every bounded $g$, where
$\nu=a^2/4$.  For bounded $h$ with $\tilde h:=h-\E h(N)$, $N\sim\Bes(k,a)$
with mass function $(p_n)$ as in \eqref{eq:pmf}, the Stein equation
$\nu g(n{+}1)-n(n{+}k)g(n)=\tilde h(n)$ has the unique bounded solution
$g_h(0)=0$,
\begin{equation}\label{eq:solution}
g_h(m{+}1)=\frac1{\nu\,p_m}\sum_{n=0}^m\tilde h(n)\,p_n
=-\frac1{\nu\,p_m}\sum_{n=m+1}^\infty\tilde h(n)\,p_n,\qquad m\ge0.
\end{equation}
\end{theorem}

The Stein \emph{factors} --- the bounds on $g_h$ that convert operator
identities into distances --- follow from a two-sided monotonicity property of log-concave
(equivalently, strongly unimodal \cite{keilsongerber1971}) laws, of
independent interest.

\begin{lemma}[Two monotone Mills ratios]\label{lem:hazard}
For any log-concave probability mass function $(p_n)_{n\ge0}$ with $p_n>0$
for all $n$, with cdf $F(m)=\sum_{n\le m}p_n$ and survival $G(m)=\sum_{n\ge m}p_n$:
\emph{(H1)} $m\mapsto F(m)/p_m$ is nondecreasing;
\emph{(H2)} $m\mapsto G(m)/p_m$ is nonincreasing.
\end{lemma}

\begin{theorem}[Stein factors]\label{thm:factors}
Let $g_A$ be the solution \eqref{eq:solution} for $h=\mathbf 1_A$,
$A\subseteq\Z_{\ge0}$, and $\Delta g(m)=g(m{+}1)-g(m)$; here $F$, $G$,
$(p_n)$ are the $\Bes(k,a)$ cdf, survival function
$G(m)=\Prob(N{\ge}m)$, and mass function \eqref{eq:pmf}, and $\nu=a^2/4$.
Then
\begin{align}
\textup{(a)}\quad
\Lm:=\sup_A\sup_{m\ge1}|\Delta g_A(m)|
&=\sup_{j\ge1}\Bigl[\frac{G(j{+}1)}\nu+\frac{F(j{-}1)}{j(j{+}k)}\Bigr]
\ \le\ \min\Bigl\{\tfrac1\nu+\tfrac1{k+1},\ \tfrac3\nu\Bigr\},
\label{eq:lambda}\\[2pt]
\textup{(b)}\quad
\Gm:=\sup_A\sup_{m\ge0}|g_A(m{+}1)|
&=\sup_{m\ge0}\frac{F(m)\,G(m{+}1)}{\nu\,p_m}
\ \le\ \frac{2}{\sqrt{k^2+4\nu}-3}\quad\bigl(\text{if }\sqrt{k^2+4\nu}>3\bigr).
\label{eq:gamma}
\end{align}
\end{theorem}

\begin{remark}[General bounded test functions]\label{rem:testfns}
The factors are stated for indicators, which is all that total variation
requires.  Nothing is lost: the solution \eqref{eq:solution} is linear in
$h$, and for each $m$ the functional
$h\mapsto g_h(m{+}1)=(\nu p_m)^{-1}\sum_n h(n)\,
\bigl(p_n\mathbf 1_{\{n\le m\}}-F(m)\,p_n\bigr)$ is maximized over
$\|h\|_\infty\le1$ at the $\pm1$-valued $h=2\cdot\mathbf 1_{\{0,\dots,m\}}-1$,
with maximum $\sum_n\bigl|p_n\mathbf 1_{\{n\le m\}}-F(m)p_n\bigr|/(\nu p_m)
=2F(m)G(m{+}1)/(\nu p_m)$; and every $\pm1$-valued $h$ equals
$2\cdot\mathbf 1_A-1$ for some $A$, so that $g_h=2g_A$.  Hence
\[
\sup_{\|h\|_\infty\le1}\|g_h\|=2\,\Gm,
\qquad
\sup_{\|h\|_\infty\le1}\|\Delta g_h\|=2\,\Lm,
\]
both suprema attained.  The same doubling holds verbatim for
$\Lambda_{\kappa,w}$ and $\Gamma_{\kappa,w}$ in \S\ref{sec:mrow}.
\end{remark}

The equality in \eqref{eq:lambda} says the increment factor is exactly the
worst singleton solution.  Numerically it equals $(1{+}o(1))/\nu$ as
$a\to\infty$ with $k$ fixed --- the analogue of the Barbour--Eagleson factor
$(1{-}e^{-\lambda})/\lambda$ \cite{barboureagleson1983} --- so the closed constant $3/\nu$ is loose by
$\approx3$ (Open Problem~\ref{prob:sharp}; Remark~\ref{rem:general} shows
$\approx1.4/\nu$ is already available from \cite{brownxia2001,er2008}).  The uniform factor obeys
$\Gm\asymp\sigma_N/\nu$, matching $\|g\|\asymp\sigma/\lambda$ for Poisson.
Proofs: \S\ref{sec:pf-operator}.

\section{Two moments suffice: an AGG-type theorem}\label{sec:agg}

\subsection{Quadratic bias and its fixed point}

\begin{definition}[Quadratic-bias transform]\label{def:qbias}
For a $\Z_{\ge0}$-valued $W$ with $\beta:=\E[W(W{+}k)]\in(0,\infty)$, the
\emph{quadratic-bias transform} $\Wsq$ is the $\Z_{\ge1}$-valued variable with
$\Prob(\Wsq{=}n)=n(n{+}k)\Prob(W{=}n)/\beta$, $n\ge1$; and the
\emph{centered transform}
\[
\mathcal T_kW\ :\overset{d}{=}\ \Wsq-1,
\]
a map from $\{$laws on $\Z_{\ge0}$ with $\E W^2<\infty$,
$\beta>0\}$ to laws on $\Z_{\ge0}$.
\end{definition}

This is the size-bias transform with weight $n$ replaced by the death rate
$n(n{+}k)$; it plays the identical structural role.

\begin{proposition}[Fixed points of the transform]\label{prop:fixed}
For $W$ with $\E W^2<\infty$ and $\beta=\E[W(W{+}k)]>0$, with $\Wsq$ and
$\mathcal T_kW=\Wsq{-}1$ as in Definition~\ref{def:qbias}:
\[
\mathcal T_kW\overset{d}{=}W
\quad\Longleftrightarrow\quad
\Wsq\overset{d}{=}W+1
\quad\Longleftrightarrow\quad
W\sim\Bes\bigl(k,\,2\sqrt\beta\bigr).
\]
In particular the fixed points of $\mathcal T_k$ are \emph{exactly} the
Bessel laws $\{\Bes(k,a):a>0\}$ --- in precise analogy with the zero-bias
transformation, whose fixed points are the centered normal family
\cite{goldsteinreinert1997}, and with the size-bias-minus-one map, whose
fixed points are the Poisson family.
\end{proposition}

\subsection{The approximation theorem}

\begin{theorem}[Bessel approximation --- two moments suffice]\label{thm:main}
Let $W$ be $\Z_{\ge0}$-valued, $\E W^2<\infty$, $\beta=\E[W(W{+}k)]>0$, and
let $\Wsq$ \textup{(}Definition~\ref{def:qbias}\textup{)} be coupled to
$W$ in any way.  Then for every $a>0$, with
$\nu=a^2/4$ and $N\sim\Bes(k,a)$,
\begin{equation}\label{eq:master}
\dtv(W,N)\ \le\ c_1(k,\nu)\,\E\bigl|W+1-\Wsq\bigr|\ +\ \Gm\,|\beta-\nu|,
\qquad c_1(k,\nu):=\min\Bigl\{3,\ 1+\tfrac\nu{k+1}\Bigr\}.
\end{equation}
Here $\Gm$ is the uniform Stein factor \eqref{eq:gamma} at this $\nu$.
With the matched parameter $a=2\sqrt\beta$ \textup{(}so $\nu=\beta$\textup{)}
the second term vanishes and
\begin{equation}\label{eq:matched}
\dtv\bigl(W,\Bes(k,2\sqrt\beta)\bigr)\ \le\
c_1(k,\beta)\,\dw\bigl(W{+}1,\ \Wsq\bigr),
\end{equation}
$\dw$ the Wasserstein-$1$ distance on $\Z_{\ge1}$.  Since $\dw$ is
shift-invariant, \eqref{eq:matched} reads equivalently \textup{(}recall
$\mathcal T_kW=\Wsq{-}1$\textup{)}
\begin{equation}\label{eq:fixedform}
\dtv\bigl(W,\ \Bes(k,2\sqrt\beta)\bigr)\ \le\
c_1(k,\beta)\;\dw\bigl(W,\ \mathcal T_kW\bigr):
\end{equation}
\emph{the total-variation distance from $W$ to the matched Bessel law is
controlled by how far $W$ is from being a fixed point of the transform.}
\end{theorem}

The shape of \eqref{eq:matched}--\eqref{eq:fixedform} is exactly the
size-bias form of Poisson approximation and the zero-bias form of normal
approximation \cite{goldsteinreinert1997}: one moment $\beta$ fixes which
member of the fixed-point family to compare with, and the whole error is the
displacement of $W$ under one application of the transform.  If $W$ is
Bessel then $\mathcal T_kW\overset{d}{=}W$ and the identity coupling gives
zero.  Proofs: \S\ref{sec:pf-agg}.

\subsection{The coupling}\label{sec:coupling}

Everything downstream is a construction of the coupling in
\eqref{eq:matched}, and the conditioning forces its geometry.  We define
the objects first and then explain why they are the right ones.  Represent
$W=\mathcal L(V\mid U-V=k)$ with $(U,V)$ a pair of $\Z_{\ge0}$-valued
counts, $\nu:=\E[UV]\in(0,\infty)$ and $\Prob(U-V=k)>0$.

\begin{definition}[Product bias, paired increment, pair coupling]\label{def:paircoupling}
\textup{(i)} The \emph{product-biased pair} $(U^{\square},V^{\square})$ has
law $d\Prob^{\square}(u,v)=uv\,d\Prob(u,v)/\nu$.
\textup{(ii)} The \emph{paired increment} is the map
$(u,v)\mapsto(u{+}1,v{+}1)$; it preserves every fiber $\{u-v=k\}$, since
$(u{+}1)-(v{+}1)=u-v$.
\textup{(iii)} The \emph{pair coupling} is a joint realization of
$(U^{\square},V^{\square})$ and $(U{+}1,V{+}1)$ on one probability space,
optimal for the $\ell^1$ cost
$\E\bigl[|U^{\square}{-}U{-}1|+|V^{\square}{-}V{-}1|\bigr]$.
\end{definition}

These are the right objects because of two exact identities.  On the fiber
$\{U-V=k\}$ one has $U=V+k$, hence the pointwise identity
\begin{equation}\label{eq:fiberid}
W(W+k)=V(V+k)=V\,U ,
\end{equation}
so quadratic-biasing $W$ is product-biasing $(U,V)$ restricted to the same
fiber, while the paired increment reproduces $W{+}1$ on that fiber:
\begin{equation}\label{eq:qbprod}
\Wsq\overset{d}{=}\mathcal L\bigl(V^{\square}\,\big|\,
U^{\square}{-}V^{\square}=k\bigr),
\qquad
W{+}1\overset{d}{=}\mathcal L\bigl(V{+}1\,\big|\,
(U{+}1)-(V{+}1)=k\bigr).
\end{equation}
For the first: $\Prob(\Wsq{=}w)\propto w(w{+}k)\Prob(V{=}w,\,U{-}V{=}k)
=\E[VU;\,V{=}w,\,U{-}V{=}k]$ by \eqref{eq:fiberid}, which is proportional
to $\Prob(V^{\square}{=}w,\,U^{\square}{-}V^{\square}{=}k)$; the second is
the fiber-invariance of the paired increment.  So $\Wsq$ and $W{+}1$ are
the two sides of a single comparison --- the product-biased pair against
the incremented pair, both viewed on the fiber --- and $\dw(W{+}1,\Wsq)$
is controlled by the cost of the pair coupling together with the change of
fiber weights, the latter quantified by the fiber constants of
Definition~\ref{def:cols}.  Two consequences drive the rest of the paper:
\begin{itemize}[leftmargin=1.5em,itemsep=2pt,topsep=2pt]
\item \textbf{Locality (\S\ref{sec:expansion}).}  For a column array
$U=\sum_jX_j$, $V=\sum_jY_j$, the weight $UV=\sum_{j,l}X_jY_l$ factorizes,
and with it the pair coupling: the pair $(j,l)$ is selected with
probability $m_jn_l/\nu$ off the diagonal and $\chi_j/\nu$ on it, the
increment lands as $+1$ in column $j$'s $X$-count and $+1$ in column $l$'s
$Y$-count, and the coupling cost becomes a sum of per-column tilt
distances --- the AGG structure.
\item \textbf{Exactness of the surrogate (\S\ref{sec:pdc}).}  The same
fiber-preserving structure lets a Bessel proposal dominate the first-half
sampling target with a constant independent of the split value, making the
rejection envelope exact.
\end{itemize}
The one disanalogy with Poisson is that the single added point of the
size-bias coupling becomes a \emph{paired} point addition; the conditioning
forces the pairing, and the fiber-invariance of the pair is what keeps the
method from paying the fiber probability
$\Prob(U{-}V{=}k)\asymp\sigma^{-1}$.

\section{The dependency expansion (AGG for column arrays)}\label{sec:expansion}

\begin{definition}[Columns, costs, fiber constants]\label{def:cols}
Let $Z_j=(X_j,Y_j)\in\Z_{\ge0}^2$, $j=1,\dots,c$, be independent across $j$,
with arbitrary joint law within each column,
$\E[X_j^2{+}Y_j^2]<\infty$, and
$\E\bigl[X_jY_j(X_j{+}Y_j)\bigr]<\infty$ (the latter makes the product-tilt
costs below finite --- second moments alone do not; if any cost is infinite,
the bounds that use it hold trivially).
Put $U=\sum_jX_j$, $V=\sum_jY_j$, $D=U-V$, $m_j=\E X_j$, $n_j=\E Y_j$,
$\chi_j=\E[X_jY_j]$, $N_U=\sum_jm_j$, $N_V=\sum_jn_j$, and
\begin{equation}\label{eq:nu}
\nu:=\E[UV]=\sum_{j\ne l}m_jn_l+\sum_j\chi_j .
\end{equation}
Fix $k\ge0$ with $\Prob(D=k)>0$ and set $W\sim\mathcal L(V\mid D=k)$.  Tilt a
column by $X$-, $Y$-, or product-bias,
\[
d\widehat\Prob_j(x,y)=\frac{x\,d\Prob_j}{m_j},\qquad
d\widetilde\Prob_j(x,y)=\frac{y\,d\Prob_j}{n_j},\qquad
d\overline\Prob_j(x,y)=\frac{xy\,d\Prob_j}{\chi_j},
\]
and define the \emph{local costs} in the $\ell^1$-Wasserstein distance on
$\Z^2$,
\[
\alpha_j=\dw(\Zh_j,Z_j{+}(1,0)),\qquad
\beta_j=\dw(\Zt_j,Z_j{+}(0,1)),\qquad
\gamma_j=\dw(\Zb_j,Z_j{+}(1,1)),
\]
and their moment-weighted aggregates
\begin{equation}\label{eq:ABC}
A:=\sum_{j=1}^cm_j\,\alpha_j,\qquad
B:=\sum_{j=1}^cn_j\,\beta_j,\qquad
C:=\sum_{j=1}^c\chi_j\,\gamma_j
\end{equation}
(terms whose weight $m_j$, $n_j$, or $\chi_j$ vanishes are absent from the
corresponding sums, and their tilts are then not needed).  Two families of \emph{systems} enter: off-diagonal two-column systems,
indexed by ordered pairs $(j,l)$ with $j\ne l$, and diagonal one-column
systems, indexed by single columns $j$ with $\chi_j>0$; the fiber constants
$\theta,\tau_1,\omega$ below are suprema over \emph{both} families.  For an
ordered pair $(j,l)$, $j\ne l$
\textup{(}resp.\ a single index $j$\textup{)}, let $V_R,D_R$ be the sums over
the remaining columns and $\mathcal S_{jl}\subset\Z$ the integer interval
hull of the achievable values of $k$ minus the difference contributions of
the involved columns, in both original and biased configurations; set
\[
\theta=\max_{(j,l)}\max_{s\in\mathcal S_{jl}}\frac{\Prob(D_R{=}s)}{\Prob(D{=}k)},
\qquad
\tau_1=\max_{(j,l)}\max_{s,s+1\in\mathcal S_{jl}}
\frac{|\Prob(D_R{=}s{+}1)-\Prob(D_R{=}s)|}{\Prob(D{=}k)},
\]
\[
\omega=\max_{(j,l)}\ \max_{\substack{s,s+1\in\mathcal S_{jl}\\
\Prob(D_R=s)\Prob(D_R=s+1)>0}}
\dw\bigl(\mathcal L(V_R\mid D_R{=}s{+}1),\ \mathcal L(V_R\mid D_R{=}s)\bigr).
\]
\end{definition}

The costs $\alpha_j,\beta_j,\gamma_j$ measure \emph{non-Poissonness}: they
vanish exactly when the column is an independent Poisson pair
(Lemma~\ref{lem:vanish}); they are the analogues of the $b_1,b_2$ terms.  The
constants $\theta,\tau_1,\omega$ are the price of the conditioning:
anti-concentration, pmf-increment, and fiber-shift-stability moduli of the
\emph{reduced} system.

\begin{theorem}[Dependency expansion]\label{thm:expansion}
Under Definition~\ref{def:cols}, with $\nu$ as in \eqref{eq:nu},
$a=2\sqrt\nu$, $A,B,C$ the aggregated costs of \eqref{eq:ABC}, and
$\Lm,\Gm$ the Stein factors \eqref{eq:lambda}--\eqref{eq:gamma} at this
$\nu$,
\begin{equation}\label{eq:expansion}
\dtv\bigl(W,\Bes(k,a)\bigr)\ \le\
\bigl[\theta(1{+}\omega)\Lm+\tau_1\Gm\bigr]\,\bigl(N_VA+N_UB+C\bigr).
\end{equation}
The crude form, replacing the fiber-increment telescoping by a union bound,
is
\begin{equation}\label{eq:crude}
\dtv\bigl(W,\Bes(k,a)\bigr)\ \le\
\theta\,(\Lm+2\Gm)\,\bigl(N_VA+N_UB+C\bigr).
\end{equation}
\end{theorem}

The parallel with \eqref{eq:aggtemplate} is exact: $N_VA+N_UB$ is the $b_1$
term (products of first moments over interacting column pairs), $C$ is the
$b_2$ term (within-column co-occurrence), the $b_3$ term is absent because
the columns are independent --- precisely the regime $b_3=0$, with dependent
arrays reintroducing it (Open Problem~\ref{prob:dep}) --- and
$\theta,\tau_1,\omega$ are new
--- the conditioning is paid not as the catastrophic
$\Prob(D{=}k)^{-1}\asymp\sigma_D$ but as $O(1)$ stability constants of the
reduced fiber, certifiable by local limit theorems (\S\ref{sec:enum}).  Note
also that the matched parameter $\nu=\E[UV]$ in \eqref{eq:nu} is an
\emph{unconditional} second moment --- ``two moments suffice'' is literal: the
approximating Bessel law is determined by the first two joint moments of the
array.  In the bulk regime \eqref{eq:expansion} reads
$O\bigl((N_VA{+}N_UB{+}C)/\nu\bigr)$.  Since
$\nu=N_UN_V-\sum_jm_jn_j+\sum_j\chi_j$ by \eqref{eq:nu}, so that
$\nu\asymp N_UN_V$ when no single column dominates, the first two terms
are comparable to $\sum_jm_j\alpha_j\big/\sum_jm_j$ and
$\sum_jn_j\beta_j\big/\sum_jn_j$: weighted \emph{averages} of the
per-column costs, not sums.  The bound therefore does not grow with the
number of columns; it is set by the cost of a typical column.  The engine of the proof is an exact \emph{transfer
identity} (Lemma~\ref{lem:transfer}) implementing the product-bias
factorization \eqref{eq:fiberid}, coupled by the fiber-preserving paired
modification of \S\ref{sec:coupling}.  Proofs: \S\ref{sec:pf-expansion}.

\section{Application: $2\times c$ contingency tables}\label{sec:tables}

The apparatus now pays off for one family of objects --- the
weighted-table measures on the fiber $\{D=k\}$ --- in three registers: a
total-variation limit theorem for the conditioned row sum
(\S\ref{sec:dtv}), a completely explicit enumeration estimate
(\S\ref{sec:enum}), and an exact sampler with certified acceptance and a
cap on what any block split can achieve (\S\ref{sec:pdc}).  A transversal
fourth result compares the conditioned process with the free one in the
Arratia--Tavar\'e sense (\S\ref{sec:ipa}).  The
conditional law of a table given its margins is classical
\cite{fisher1935,yates1934,agresti1992}; Monte Carlo and algebraic methods
for tables with fixed margins include Markov bases
\cite{diaconissturmfels1998}, sequential importance sampling \cite{cdhl2005},
and permanent-based counting \cite{barvinok2008}; and probabilistic
divide-and-conquer has been applied to contingency tables by recursing on
exact conditional distributions column by column
\cite{desalvozhao2020,desalvo2018}.  The split \eqref{eq:split} on a
\emph{scalar statistic}, with a Bessel surrogate and closed-form certified
acceptance, is what is new here.

\emph{Which margins are fixed.}  Throughout this section and
\S\ref{sec:mrow}, the conditioning event is the row-sum \emph{difference}
$\{U{-}V{=}k\}$ (the $r{-}1$ successive differences in \S\ref{sec:mrow});
the column sums are left entirely free.  Equivalently, by
\eqref{eq:split}, the ensemble is the Bessel mixture, over the offset $v$,
of the classical fixed-row-sum ensembles $(v{+}k,v)$: the row margins
fluctuate by a common Bessel-distributed offset, and the columns are
unconstrained.  The Fisher--Yates ensemble
\cite{fisher1935,yates1934} fixes \emph{both} margins and is recovered from
ours by further conditioning on the column-sum vector; that considerably
harder conditioning is Open Problem~\ref{prob:margins}.
Conditionally on its \emph{realized} margins, however, a table drawn from
our ensemble is exactly Fisher--Yates distributed --- centrally so
(uniform, in the binary case) when the entry weights are rank one, and
noncentrally in general --- which is why the Fisher--Yates weighted count
appears in the column-sum decomposition of
Proposition~\ref{prop:hbrep}(d).  In particular, the sampler of
\S\ref{sec:pdc} does not by itself produce null tables for Fisher's exact
test at prescribed margins; certifying that step is precisely Open
Problem~\ref{prob:margins}.  The \emph{Poisson} weights are what make the
margin-fibers Fisher--Yates; conditioning independent \emph{geometric}
entries with product-form parameters on both margins yields instead the
\emph{uniform} table ensemble --- the maximum-entropy model behind
Barvinok's typical tables \cite{barvinok2010}.  The geometric setting is
taken up in the companion paper \cite{desalvo2026}.

\subsection{Two column models}\label{sec:models}

All results of this section concern one of two column models.  We fix
their notation once, so that every statement below needs only a single
reference.

\begin{definition}[Binary model]\label{def:binary}
Entries $X_j\sim\mathrm{Bern}(p_j)$ and $Y_j\sim\mathrm{Bern}(q_j)$,
$j\in[c]$, all independent, with $p_j,q_j\in(0,1)$.  Derived quantities:
row sums $U=\sum_jX_j$, $V=\sum_jY_j$ and difference $D=U-V$; means
$N_U=\sum_jp_j$, $N_V=\sum_jq_j$, matched parameter $\nu=N_UN_V$
\textup{(}$=\E[UV]$, by independence\textup{)}, and $a=2\sqrt\nu$; odds
$x_j=p_j/(1{-}p_j)$, $y_j=q_j/(1{-}q_j)$ with sums $S_x=\sum_jx_j$,
$S_y=\sum_jy_j$, tilted parameter $\nu'=S_xS_y$, and $a'=2\sqrt{\nu'}$;
$b_{k,a}$ and $b'=b_{k,a'}$, the $\Bes(k,a)$ and $\Bes(k,a')$ mass
functions \eqref{eq:pmf}; row masses
$E_w(p)=\Prob(U{=}w)=\sum_{|S|=w}\prod_{j\in S}p_j\prod_{j\notin S}(1{-}p_j)$
and likewise $E_w(q)$; fiber masses $f(v)=E_{v+k}(p)\,E_v(q)$, supported
on $0\le v\le v_{\max}:=c-k$; Skellam masses
$\Sk_m=e^{-N_U-N_V}(N_U/N_V)^{m/2}I_m(a)$ for $m\in\Z$
\textup{(}$I_{-m}=I_m$\textup{)}, with certificate
$\delta=(1{\wedge}N_U^{-1})\sum_jp_j^2+(1{\wedge}N_V^{-1})\sum_jq_j^2$.
For $B\subseteq[c]$, a superscript $(B^c)$ on any of these denotes the
same quantity built from the columns outside $B$.
\end{definition}

\begin{definition}[Capped-Poisson model]\label{def:capped}
Entries $X_j\sim\Poi(\lambda_{1j})$ conditioned on $\{X_j\le\kappa_j\}$
and $Y_j\sim\Poi(\lambda_{2j})$ conditioned on $\{Y_j\le\kappa_j\}$, all
independent, with caps $\kappa_j\in\Z_{\ge1}$.  Derived quantities:
cap coverages $F^{(i)}_j=\Prob\bigl(\Poi(\lambda_{ij})\le\kappa_j\bigr)$,
$i=1,2$, and defect $\delta'=\sum_{i,j}\bigl(1-F^{(i)}_j\bigr)$; means
$N_i=\sum_j\lambda_{ij}$, tilted parameter $\nu'=N_1N_2$, and
$a'=2\sqrt{\nu'}$; $b'=b_{k,a'}$, the $\Bes(k,a')$ mass function
\eqref{eq:pmf}; row sums $U,V$, difference $D$, and fiber masses
$f(v)=\Prob(U{=}v{+}k)\Prob(V{=}v)$ as in the binary model, supported on
$0\le v\le v_{\max}:=\sum_j\kappa_j-k$; free-Poisson Skellam masses
$\Sk'_m=e^{-N_1-N_2}(N_1/N_2)^{m/2}I_m(a')$ for $m\in\Z$
\textup{(}$I_{-m}=I_m$\textup{)}.  The superscript-$(B^c)$ convention of
Definition~\ref{def:binary} applies.
\end{definition}

\subsection{Total variation}\label{sec:dtv}

Consider a Poissonized $2\times c$ table --- independent Poisson entries
being the classical \emph{Poisson sampling scheme} of categorical data
analysis, whose conditional equivalence with the multinomial schemes is
Birch's theorem \cite{birch1963,bfh1975}, and the table case of the
Poissonize-and-condition device of combinatorial probability
\cite{abt2003} --- with
$T_{ij}\sim\Poi(\lambda_i\mu_j)$ independent, row sums $R_i=\sum_jT_{ij}$,
and difference $D=R_1-R_2$.  By Definition~\ref{def:cols} and
Lemma~\ref{lem:condrep}, $\mathcal L(R_2\mid D{=}k)=\Bes(k,2M\sqrt{\lambda_1
\lambda_2})$ exactly, $M=\sum_j\mu_j$, and conditionally on the row totals
the rows are independent multinomials --- the two-stage decomposition used
for sampling in \S\ref{sec:pdc}.  Theorem~\ref{thm:expansion} quantifies the
cost of perturbing the entries away from Poisson.  Two families have
closed-form costs.

\begin{corollary}[Sparse binary tables: a conditioned law of small numbers]
\label{cor:binary}
For the binary model of Definition~\ref{def:binary}, viewed in the
expansion setting of Definition~\ref{def:cols}, the size-biased variable
is the constant $1$, so
\[
\alpha_j=p_j,\qquad\beta_j=q_j,\qquad\gamma_j=p_j{+}q_j,\qquad \chi_j=p_jq_j,
\qquad\nu=N_UN_V,
\]
and Theorem~\ref{thm:expansion} gives \textup{(}$\theta,\tau_1,\omega$
from Definition~\ref{def:cols}; $\Lm,\Gm$ from
\eqref{eq:lambda}--\eqref{eq:gamma}\textup{)}, with $p^*=\max_jp_j$,
$q^*=\max_jq_j$,
\begin{equation}\label{eq:binarybound}
\begin{aligned}
\dtv\bigl(\mathcal L(V\mid U{-}V{=}k),\ \Bes(k,2\sqrt\nu)\bigr)
&\le\bigl[\theta(1{+}\omega)\Lm+\tau_1\Gm\bigr]\\
&\qquad\times
\Bigl(N_V\!\sum_jp_j^2+N_U\!\sum_jq_j^2+\sum_jp_jq_j(p_j{+}q_j)\Bigr).
\end{aligned}
\end{equation}
Along any sequence of systems with stable fiber constants ---
$\theta(1{+}\omega)=O(1)$ and $\nu\tau_1\Gm=O(1)$, as in the numerical
settings of \S\ref{sec:numerics} and Open Problem~\ref{prob:fiber} --- the
right-hand side is $O(p^*{+}q^*)$.
As $k\to\infty$ with $\nu/k\to\lambda'$ fixed,
$\Bes(k,2\sqrt\nu)\to\Poi(\lambda')$, so for strongly unbalanced margins the
corollary degenerates to the classical conditioned Chen--Stein law of small
numbers.
\end{corollary}

\begin{corollary}[Capped Poisson tables]\label{cor:capped}
For the capped model of Definition~\ref{def:capped} \textup{(}row index
suppressed: $\lambda_j$ stands for either $\lambda_{1j}$ or
$\lambda_{2j}$\textup{)}, viewed in the expansion setting of
Definition~\ref{def:cols}, size-biasing a truncation lowers its cap
by one, so with $\pi_\lambda,F_\lambda$ the Poisson pmf and cdf,
\begin{equation}\label{eq:cappedcost}
\alpha_j=\dw\bigl(\Poi(\lambda_j)\!\mid\!\le\!\kappa_j{-}1,\
\Poi(\lambda_j)\!\mid\!\le\!\kappa_j\bigr)
=\frac{\pi_{\lambda_j}(\kappa_j)}{F_{\lambda_j}(\kappa_j)\,F_{\lambda_j}(\kappa_j{-}1)}
\sum_{m=0}^{\kappa_j-1}F_{\lambda_j}(m)
\ \le\ \frac{\kappa_j\,\pi_{\lambda_j}(\kappa_j)}
{F_{\lambda_j}(\kappa_j)\,F_{\lambda_j}(\kappa_j{-}1)},
\end{equation}
exponentially small once $\kappa_j-\lambda_j\gg\sqrt{\lambda_j}$: capping is a
tame perturbation.
\end{corollary}

Exact dynamic-programming computations (Table~\ref{tab:master}) show
\eqref{eq:expansion} tracking the true distance proportionally over four
orders of magnitude; for sparse binary systems at fixed $\nu$, halving $p^*$
halves both the truth and the bound, with slack pinned near $37\times$: the
rate $\Theta(p^*{+}q^*)$ is exact and the loss is constant-only.
Proofs: \S\ref{sec:pf-dtv}.

\subsection{Explicit enumeration}\label{sec:enum}

Recall the row masses $E_w(p)=\Prob(U{=}w)$ of
Definition~\ref{def:binary} --- Poisson--binomial masses, equal to the
weighted subset counts $e_w(x)$ normalized by $\prod_j(1{-}p_j)$
\textup{(}Corollary~\ref{cor:counts}\textup{)}.  A pair $(S,T)$ with $|S|=v{+}k$,
$|T|=v$ is a $2\times c$ zero--one table with row-sum discrepancy $k$;
$E_{v+k}(p)E_v(q)$ is the probability mass of such tables resolved by the
second-row size $v$; the corresponding weighted \emph{count} is the
$T_k(v)$ of Corollary~\ref{cor:counts}.  For \emph{uncapped Poisson} columns there is nothing to estimate:
$\Prob(U{=}v{+}k)\Prob(V{=}v)=\Sk_k\,b_{k,a}(v)$ holds \emph{exactly}, by
Lemma~\ref{lem:condrep}.  The theorem below measures how far Bernoulli columns can drift from this
Poisson identity --- capped-Poisson columns are handled by a direct
coupling in Corollary~\ref{cor:cappedenum} --- and making the drift fully
explicit requires certifying the two implicit constants of
\S\ref{sec:expansion} in closed form.

\begin{lemma}[Explicit certificates]\label{lem:cert}
\begin{enumerate}[leftmargin=1.6em,label=\textup{(\roman*)},itemsep=1pt]
\item $\Gm\le 2/(\sqrt{k^2+4\nu}-3)$ whenever $\sqrt{k^2+4\nu}>3$, with
$\nu$ as in \eqref{eq:nu}
\textup{(}Theorem~\ref{thm:factors}\textup{(}b\textup{))}.
\item \textbf{Anti-concentration} (an Esseen-type concentration bound
\cite{esseen1966,petrov1975}).  For independent $\Z$-valued
$W_1,\dots,W_n$ with $S=\sum_jW_j$, and $\rho_j=\Prob(|W_j-W_j'|=1)$
\textup{(}$W_j'$ an independent copy of $W_j$\textup{)},
\[
\max_{s\in\Z}\Prob(S=s)\le\sqrt{\frac{\pi}{4R}},\qquad R=\sum_j\rho_j.
\]
For binary column differences, $\rho_j=2(u_j{+}v_j)(1{-}u_j{-}v_j)$,
$u_j=p_j(1{-}q_j)$, $v_j=q_j(1{-}p_j)$.
\item \textbf{Skellam replacement.}  In the binary model of
Definition~\ref{def:binary}, with $\Sk_m$ and $\delta$ as there defined,
$\ \sup_{m\in\Z}|\Prob(D{=}m)-\Sk_m|\le\delta$.
\end{enumerate}
\end{lemma}

\begin{theorem}[Skellam--Bessel factorization of the fiber masses]\label{thm:enum}
In the binary model of Definition~\ref{def:binary} \textup{(}which fixes
$a$, $\nu$, $E_w$, $\Sk_k$, $\delta$\textup{)}, let $\rho_j$ be as in
Lemma~\ref{lem:cert}\textup{(}ii\textup{)} applied to $W_j=X_j-Y_j$.  With
$\rho_{(1)}\ge\rho_{(2)}$ the two largest $\rho_j$, set
\[
R_2=\sum_j\rho_j-\rho_{(1)}-\rho_{(2)},\qquad
\theta^*=\frac{\sqrt{\pi/(4R_2)}}{\Sk_k-\delta},\qquad
\Gamma^*=\frac{2}{\sqrt{k^2+4\nu}-3},\qquad
\Lambda^*=\min\Bigl\{\tfrac1\nu+\tfrac1{k+1},\ \tfrac3\nu\Bigr\},
\]
\[
\varepsilon^*=\theta^*\bigl(\Lambda^*+2\Gamma^*\bigr)\Bigl(N_V\!\sum_jp_j^2
+N_U\!\sum_jq_j^2+\sum_jp_jq_j(p_j{+}q_j)\Bigr),
\]
assuming $c\ge2$, $R_2>0$, $\Sk_k>\delta$, and $\sqrt{k^2+4\nu}>3$; the error aggregate
$\varepsilon^*$, a certified total-variation bound, contains the
aggregated cost $B$ of \eqref{eq:ABC} via
Corollary~\ref{cor:binary} ($N_UB=N_U\sum_jq_j^2$ there).  Then for every $v\ge0$,
\begin{equation}\label{eq:enum}
\bigl|E_{v+k}(p)\,E_v(q)-\Sk_k\,b_{k,a}(v)\bigr|
\ \le\ \delta\,b_{k,a}(v)+(\Sk_k+\delta)\,\varepsilon^* .
\end{equation}
Every quantity on the right is a closed-form function of $(p,q,k)$ built from
elementary operations and the modified Bessel function $I_k$.
\end{theorem}

\begin{corollary}[Weighted table counts]\label{cor:counts}
In the binary model of Definition~\ref{def:binary}, weight each
$2\times c$ zero--one table by the product of the odds $x_j,y_j$ over its
$1$-cells, and let
$T_k(v):=e_{v+k}(x)\,e_v(y)$, with
$e_w(x)=\sum_{|S|=w}\prod_{j\in S}x_j$ the elementary symmetric
polynomial, be the weighted number of such tables with row sums
$(v{+}k,v)$.  Then $E_{v+k}(p)\,E_v(q)=T_k(v)\prod_j(1{-}p_j)(1{-}q_j)$,
and Theorem~\ref{thm:enum} reads
\[
\bigl|T_k(v)-\mathsf S\,b_{k,a}(v)\bigr|
\ \le\ \frac{\delta\,b_{k,a}(v)+(\Sk_k+\delta)\varepsilon^*}
{\prod_j(1{-}p_j)(1{-}q_j)},
\qquad
\mathsf S:=\frac{\Sk_k}{\prod_j(1{-}p_j)(1{-}q_j)}.
\]
For equal odds $x_j=y_j=t$ this is a statement about integer objects:
$T_k(v)=t^{2v+k}\binom{c}{v+k}\binom{c}{v}$, so the diagonal coefficient
profile $v\mapsto t^{2v}\binom{c}{v+k}\binom{c}{v}$ is Bessel-shaped with
parameter $a=2ct/(1{+}t)$, certified whenever $t$ lies in the sparse
regime of Corollary~\ref{cor:binary}.
\end{corollary}

\begin{remark}[Where the Bessel profile lives]\label{rem:dense}
At $p_j=q_j=\tfrac12$ every weight is $1$:
$T_k(v)=\binom{c}{v+k}\binom{c}{v}$ counts tables exactly, and
$\sum_vT_k(v)=\binom{2c}{c-k}$ by Vandermonde.  This is precisely where
the estimate degenerates --- $\rho_j=0$ there, so $R_2=0$ and the
hypothesis of Theorem~\ref{thm:enum} fails --- and genuinely so: the
conditional variance of the count profile is $\asymp c/8$, while
$\Bes(k,c)$ has variance $\asymp c/4$, a width mismatch of $\sqrt2$
(numerically at $c=40$, $k=2$: $\max_v|T_k(v)/\mathsf S-b_{k,c}(v)|
\approx0.11$ against a peak mass of $0.13$).  The Bessel shape of the
weighted counts is a sparse-regime phenomenon; the dense, unweighted
profile is Gaussian with the narrower width.
\end{remark}

\begin{remark}[When the estimate has content]\label{rem:nontrivial}
Three regimes calibrate Theorem~\ref{thm:enum}.  \textup{(i)} Uncapped
Poisson columns: the two sides agree exactly, and there is nothing to
estimate.  \textup{(ii)} Equal binary weights:
$T_k(v)=t^{2v+k}\binom{c}{v+k}\binom{c}{v}$ is closed-form, and Stirling's
approximation already yields its profile; here the theorem adds only the
certified constants.  \textup{(iii)} Heterogeneous weights and caps:
$e_{v+k}(x)\,e_v(y)$ with distinct $x_j,y_j$, and the weighted
restricted-composition coefficients of the capped case, are computable by
dynamic programming but admit no closed form, and their bulk asymptotics
classically require saddle-point or local-limit analysis with unspecified
constants \cite{rothszekeres1954}.  This is where the theorem earns its
keep: one certified inequality controls the entire $k$-shifted profile,
uniformly in $v$, with every constant explicit --- valid exactly in the
sparse regime ($p^*{+}q^*$ small, $R_2>0$) and degenerating at the dense
boundary of Remark~\ref{rem:dense}.
\end{remark}

\begin{corollary}[Capped fiber masses]\label{cor:cappedenum}
In the capped model of Definition~\ref{def:capped} \textup{(}which fixes
$N_i$, $F^{(i)}_j$, $a'$, $b'$, $\Sk'_k$, $\delta'$\textup{)}: uniformly
in $v\ge0$,
\[
\bigl|\Prob(U{=}v{+}k)\,\Prob(V{=}v)-\Sk'_k\,b'(v)\bigr|\ \le\ \delta',
\]
and $\delta'$ is exponentially small once every
$\kappa_j-\lambda_{ij}\gg\sqrt{\lambda_{ij}}$
\textup{(}Corollary~\ref{cor:capped}\textup{)}.  Multiplying by
$\prod_{i,j}F^{(i)}_j$ turns this into a certified Bessel profile for the
weighted restricted-composition coefficients of
Remark~\ref{rem:nontrivial}\textup{(}iii\textup{)}.
\end{corollary}

The certified error has the correct rate --- $\Theta(1/c)$ at fixed means,
matching the truth --- with constants loose by $\approx4\times10^2$, itemized
in Remark~\ref{rem:losses}; the itemization is a map of exactly which open
problems buy back which factor.  Proofs: \S\ref{sec:pf-enum}.

\subsection{Independent process approximation}\label{sec:ipa}

In the language of the Poissonize-and-condition paradigm
\cite{abt2003,arratiatavare1994}, the ensemble is an \emph{independent
process approximation}: on a block of columns whose difference is a small
part of the whole, the conditioning $\{D{=}k\}$ is invisible.  The next
result makes this exact --- the process-level distance collapses to a
one-dimensional one --- and certified.

\begin{theorem}[Independent process approximation]\label{thm:ipa}
In the setting of Definition~\ref{def:cols}, fix $B\subseteq[c]$ and write
$Z_B=(X_j,Y_j)_{j\in B}$, $D_B=\sum_{j\in B}(X_j{-}Y_j)$, and
$D_{B^c}=D-D_B$.  Then
\begin{equation}\label{eq:ipa}
\dtv\bigl(\mathcal L(Z_B\mid D{=}k),\,\mathcal L(Z_B)\bigr)
=\dtv\bigl(\mathcal L(D_B\mid D{=}k),\,\mathcal L(D_B)\bigr)
=\frac12\,\E\biggl|\frac{\Prob(D_{B^c}{=}k{-}D_B)}{\Prob(D{=}k)}-1\biggr|,
\end{equation}
the expectation over the free $D_B$.  In the binary model of
Definition~\ref{def:binary}, with $\Sk_m$, $\delta$, and the
superscript-$(B^c)$ convention as there defined.  If $\Sk_k>\delta$, then
\begin{equation}\label{eq:ipabound}
\dtv\bigl(\mathcal L(Z_B\mid D{=}k),\,\mathcal L(Z_B)\bigr)
\ \le\ \frac{\E\bigl|\Sk^{(B^c)}_{k-D_B}-\Sk_k\bigr|
+\delta^{(B^c)}+\delta}{2\,(\Sk_k-\delta)}\,,
\end{equation}
and the expectation is a finite sum over the at most $2|B|{+}1$ values of
$D_B$, each term a closed-form Bessel expression.
\end{theorem}

\begin{corollary}[A vanishing conditioning effect]\label{cor:ipa}
In the binary model of Definition~\ref{def:binary}, for every
$B\subseteq[c]$,
\begin{equation}\label{eq:ipazero}
\dtv\bigl(\mathcal L(Z_B\mid D{=}k),\,\mathcal L(Z_B)\bigr)
\ \le\ \frac{2\sum_{j\in B}(p_j{+}q_j)+\delta^{(B^c)}+\delta}
{2\,(\Sk_k-\delta)}\,.
\end{equation}
In particular, in the homogeneous sparse regime $p_j=q_j=\lambda/c$ with
$\lambda>0$ and $k$ fixed, every block with $|B|=o(c)$ has
$\dtv=O(|B|/c)\to0$: the conditioning is asymptotically invisible on any
vanishing fraction of the columns.  \textup{(}The constants here are
crude; when sharpness matters, use the exact expectation in
\eqref{eq:ipabound}.\textup{)}
\end{corollary}

\begin{corollary}[Capped independent process approximation]\label{cor:cappedipa}
In the capped model of Definition~\ref{def:capped}, with $\Sk'_m$,
$\delta'$, and the superscript-$(B^c)$ convention as there defined, the bound \eqref{eq:ipabound} holds with
$(\Sk_m,\delta,\delta^{(B^c)})\rightsquigarrow
(\Sk'_m,\delta',\delta'^{(B^c)})$, and, whenever $\Sk'_k>\delta'$,
\begin{equation}\label{eq:ipazerocap}
\dtv\bigl(\mathcal L(Z_B\mid D{=}k),\,\mathcal L(Z_B)\bigr)
\ \le\ \frac{2\sum_{j\in B}(\lambda_{1j}{+}\lambda_{2j})+2\delta'}
{2\,(\Sk'_k-\delta')}\,.
\end{equation}
In the homogeneous regime $\lambda_{ij}=\lambda/c$ with caps
$\kappa_j\ge1$ and $k$ fixed, every block with $|B|=o(c)$ again has
$\dtv=O(|B|/c)\to0$.
\end{corollary}

The Chen--Stein program of \cite{agg1989} bounds not only the sum of
dependent indicators but the entire indicator \emph{process} by an
independent one, with the same $b_1{+}b_2{+}b_3$ error.  The conditioned
array admits the analogous upgrade, at the same two-moment error, with the
comparison process as independent as the conditioning permits.

\begin{theorem}[Process-level Bessel approximation]\label{thm:process}
Work in the binary model of Definition~\ref{def:binary}, or the capped
model of Definition~\ref{def:capped} with $\nu:=\E[UV]=N_1N_2$ and
$a=2\sqrt\nu$.  Let $\tilde b$ be $\Bes(k,a)$ conditioned on
$\{v\le v_{\max}\}$, and let $\nu^*_k$ be the law of the array obtained by
drawing $v\sim\tilde b$ and then, independently, each row from its free
law conditioned on its sum \textup{(}$U=v{+}k$, $V=v$: two independent
conditioned-Bernoulli rows, each exactly samplable by the recursion of
Proposition~\ref{prop:bothmargins}; capped rows likewise, by the analogous
dynamic program\textup{)}.  Then, exactly,
\begin{equation}\label{eq:process}
\dtv\bigl(\mathcal L\bigl((X_j,Y_j)_{j\in[c]}\mid D{=}k\bigr),\
\nu^*_k\bigr)
\ =\ \dtv\bigl(\mathcal L(V\mid D{=}k),\ \tilde b\,\bigr),
\end{equation}
and consequently
\[
\Bigl|\,\dtv\bigl(\mathcal L(\,\cdot\mid D{=}k),\,\nu^*_k\bigr)
-\dtv\bigl(\mathcal L(V\mid D{=}k),\,\Bes(k,a)\bigr)\Bigr|
\ \le\ \Prob\bigl(\Bes(k,a)>v_{\max}\bigr),
\]
a super-exponentially small tail in the sparse regime.
\end{theorem}

\noindent
The two-moment bound of Theorem~\ref{thm:expansion} thus controls the
total variation of the \emph{entire conditioned table}.  No
independent-columns law can play the role of $\nu^*_k$: testing the event
$\{D=k\}$ gives $\dtv\ge1-\rho(D{=}k)$ against every product law $\rho$,
which tends to one whenever the column differences are nondegenerate.
The $r$-row statement holds verbatim with $\mathrm{HB}(\kappa,\nu)$
mixing, $\nu=\prod_iN_i$, and $r$ independent conditioned rows.

\noindent
The identity \eqref{eq:ipa} is purely structural; its $r$-row version,
with a certified capped bound, is Theorem~\ref{thm:mipa}, and the
comparison with the row sums themselves fixed is
Theorem~\ref{thm:rowfixed}.  Proofs:
\S\ref{sec:pf-enum}.

\subsection{Sampling via probabilistic divide-and-conquer}\label{sec:pdc}

We sample $\mu_k=\mathcal L\bigl((X_j,Y_j)_{j}\mid D{=}k\bigr)$ ---
weighted binary tables or capped Poissonized tables --- by probabilistic
divide-and-conquer \cite{arratiadesalvo2016}, splitting on the \emph{scalar
statistic} $V$:
\begin{equation}\label{eq:split}
\mu_k=\underbrace{\mathcal L(V\mid D{=}k)}_{\text{first half}}
\ \otimes\
\underbrace{\mathcal L\bigl((X,Y)\mid U{=}v{+}k,\,V{=}v\bigr)}_{\text{second
half: exactly samplable}} .
\end{equation}
Here $\otimes$ is the semidirect product of a law and a probability kernel, $(\mu\otimes K)(\cdot)=\sum_v\mu(v)\,K(v,\cdot)$ --- a disintegration, not a product of two fixed laws, since the second factor is the conditional law of the columns at the sampled value $v$.
The first half is the enumeration object of \S\ref{sec:enum}; we sample it by
rejection from a Bessel proposal at the \emph{tilted} parameter
$\nu'=(\sum_jx_j)(\sum_jy_j)$, $x_j=p_j/(1{-}p_j)$, $y_j=q_j/(1{-}q_j)$.  The
point of the tilt is that the envelope becomes \emph{exact}: correctness
requires no approximation theorem, and the theory of \S\S\ref{sec:agg}--%
\ref{sec:enum} enters only through the \emph{certificate} on the rate.

\begin{lemma}[Exact envelopes]\label{lem:envelope}
In the binary model of Definition~\ref{def:binary}:
\begin{enumerate}[leftmargin=1.6em,label=\textup{(\alph*)},itemsep=1pt]
\item \textbf{Binary.}  With
$\Ct=\prod_j(1{-}p_j)(1{-}q_j)\cdot(S_x/S_y)^{k/2}I_k(a')$, one has
$f(v)\le\Ct\,b'(v)$ for all $v$, with pointwise ratio
\[
\frac{f(v)}{\Ct\,b'(v)}
=\Prob\bigl(\text{$v{+}k$ i.i.d.\ draws from $x/S_x$ are distinct}\bigr)
\cdot\Prob\bigl(\text{$v$ i.i.d.\ draws from $y/S_y$ are distinct}\bigr).
\]
\item \textbf{Capped Poisson} \textup{(}the model of
Definition~\ref{def:capped}\textup{)}.  With $N_1=\sum_j\lambda_{1j}$,
$N_2=\sum_j\lambda_{2j}$, $F^{(i)}_j=\Prob(\Poi(\lambda_{ij})\le\kappa_j)$,
$\nu'=N_1N_2$, and $\Ct=\Sk'_k/\prod_jF^{(1)}_jF^{(2)}_j$, where
$\Sk'_k=e^{-N_1-N_2}(N_1/N_2)^{k/2}I_k(a')$ is the \emph{uncapped} Skellam
value, one has $f(v)\le\Ct\,b'(v)$, with ratio equal to the product of the
probabilities that $\Mult(v{+}k;\lambda^{(1)}/N_1)$ \textup{(}$\lambda^{(i)}:=(\lambda_{ij})_j$\textup{)} and
$\Mult(v;\lambda^{(2)}/N_2)$ allocations respect all caps.
\end{enumerate}
\end{lemma}

\begin{definition}[Maclaurin envelope]\label{def:maclaurin}
For $t>0$ let $b_t$ denote the $\Bes(k,2\sqrt t\,)$ mass function
\eqref{eq:pmf}.  In the binary model of Definition~\ref{def:binary},
Maclaurin's inequality $e_w(x)\le S_x^{\,w}/w!$ applied to each row gives,
for every $v\ge0$ and every $t>0$,
\[
f(v)\ \le\ \Ct(t)\,b_t(v),
\qquad
\Ct(t):=\prod_j(1{-}p_j)(1{-}q_j)\,
\Bigl(\tfrac{S_x}{S_y}\Bigr)^{\!k/2}
I_k\bigl(2\sqrt t\,\bigr)
\Bigl(\tfrac{\nu'}{t}\Bigr)^{\!k/2}
\max\Bigl\{1,\ \bigl(\tfrac{\nu'}{t}\bigr)^{\!v_{\max}}\Bigr\},
\]
where $v_{\max}$ is the largest $v$ with $f(v)>0$.  The pair
$(b_t,\Ct(t))$ is the \emph{Maclaurin envelope} at tilt $t$; at $t=\nu'$
it is the envelope of Lemma~\ref{lem:envelope}, with
$\Ct:=\Ct(\nu')=\prod_j(1{-}p_j)(1{-}q_j)(S_x/S_y)^{k/2}I_k(a')$.  In the
capped model \textup{(}Definition~\ref{def:capped}\textup{)} the envelope
of Lemma~\ref{lem:envelope}\textup{(}b\textup{)} plays the same role at
$t=N_1N_2$, with $\Ct=\Sk'_k/\prod_{i,j}F^{(i)}_j$, and the general-$t$
constant carries the same displacement factor with
$\nu'\rightsquigarrow N_1N_2$.
\end{definition}

\begin{algorithm}[\textsc{Bessel--PDC}: binary tables]\label{alg:pdc}
\leavevmode\smallskip

\noindent\rule{\textwidth}{0.8pt}
\emph{Input:} weights $p,q\in(0,1)^c$; integer $k\ge0$.

\smallskip\noindent
\emph{Setup} (once):
\begin{enumerate}[label=\textup{(S\arabic*)},leftmargin=2.9em,itemsep=1pt,topsep=2pt]
\item Compute the odds $x_j=p_j/(1{-}p_j)$, $y_j=q_j/(1{-}q_j)$; the sums
$S_x=\sum_jx_j$, $S_y=\sum_jy_j$; and $\nu'=S_xS_y$.
\item Tabulate the proposal $\Bes(k,2\sqrt{\nu'})$ by the ratio recursion
$t_0=1$, $t_{v+1}=t_v\cdot\nu'/\bigl((v{+}1)(v{+}k{+}1)\bigr)$, stopping
once $t_v<10^{-18}\max_ut_u$; normalize and store the cumulative sums.
\item Build sampling tables (alias, or stored cdf with binary search) for
the two distributions $x/S_x$ and $y/S_y$ on $\{1,\dots,c\}$.
\end{enumerate}
\emph{Repeat} until an attempt accepts (each attempt starts afresh at
step~1):
\begin{enumerate}[label=\textup{(\arabic*)},leftmargin=2.9em,itemsep=1pt,topsep=2pt]
\item Draw $v\sim\Bes(k,2\sqrt{\nu'})$ by inversion on the stored cdf.
\item Draw $v{+}k$ indices i.i.d.\ from $x/S_x$.  If any index repeats
\emph{within this group}, \textbf{reject} the attempt.
\item Draw $v$ indices i.i.d.\ from $y/S_y$.  If any index repeats
\emph{within this group}, \textbf{reject} the attempt.  (The two groups may
share indices: a column may be $1$ in both rows.)
\item \textbf{Accept}: output the index sets $S=\{$step-2 draws$\}$ and
$T=\{$step-3 draws$\}$.
\end{enumerate}
\emph{Output:} $(S,T)$ is an exact draw of the two rows given
$\{|S|-|T|=k\}$; together with the sampled $v$, the accepted sets constitute the complete
output --- the second half of
\eqref{eq:split}, so no further sampling stage is needed.\\[1pt]
\rule{\textwidth}{0.8pt}
\end{algorithm}

\begin{algorithm}[\textsc{Bessel--PDC}: capped Poisson tables]\label{alg:pdccap}
\leavevmode\smallskip

\noindent\rule{\textwidth}{0.8pt}
\emph{Input:} rates $\lambda_{1j},\lambda_{2j}>0$ and caps
$\kappa_j\in\Z_{\ge0}$, $j\le c$; integer $k\ge0$.

\smallskip\noindent
\emph{Setup} (once): $N_i=\sum_j\lambda_{ij}$; $\nu'=N_1N_2$; the
$\Bes(k,2\sqrt{\nu'})$ table as in \textup{(S2)} above; sampling tables for
the two column distributions $\lambda^{(i)}/N_i$.

\smallskip\noindent
\emph{Repeat} until an attempt accepts:
\begin{enumerate}[label=\textup{(\arabic*)},leftmargin=2.9em,itemsep=1pt,topsep=2pt]
\item Draw $v\sim\Bes(k,2\sqrt{\nu'})$ by inversion.
\item Draw row $1\sim\Mult(v{+}k;\lambda^{(1)}/N_1)$ and row
$2\sim\Mult(v;\lambda^{(2)}/N_2)$ (e.g.\ as $v{+}k$, resp.\ $v$,
categorical draws with per-column counting).
\item If any entry exceeds its cap ($T_{ij}>\kappa_j$ for some $i,j$),
\textbf{reject} the attempt; else \textbf{accept} and output the table
$T=(T_{ij})$.
\end{enumerate}
\emph{Output:} $T$ is an exact draw of the capped table given
$\{R_1-R_2=k\}$.\\[1pt]
\rule{\textwidth}{0.8pt}
\end{algorithm}

\begin{theorem}[Exactness and certified acceptance]\label{thm:pdc}
Algorithms~\ref{alg:pdc} and~\ref{alg:pdccap} output an exact sample from
$\mu_k:=\mathcal L(\textup{table}\mid D{=}k)$ \textup{(}the second half of
\eqref{eq:split}, together with the sampled $v$\textup{)}, with acceptance
probability exactly $\Prob(D{=}k)/\Ct$ per attempt \textup{(}$\Ct$ as in Definition~\ref{def:maclaurin}; all remaining
quantities as in Definitions~\ref{def:binary}
and~\ref{def:capped}\textup{)}, bounded below in closed form by
\begin{equation}\label{eq:pdccert}
\text{binary: }\ \frac{\Sk_k-\delta}{\Ct},
\qquad\qquad
\text{capped: }\ \frac{\Sk'_k-\delta'}{\Sk'_k}\,\prod_jF^{(1)}_jF^{(2)}_j,
\quad \delta'=\sum_{i,j}\bigl(1-F^{(i)}_j\bigr),
\end{equation}
whenever the right-hand sides are positive.  The expected cost is $O(c)$
preprocessing plus $O\bigl((\E_{b'}[v]{+}k)/\alpha\bigr)$ per sample, where
$\alpha$ is the acceptance rate and
$\E_{b'}[v]=\frac{a'}2I_{k+1}(a')/I_k(a')\le\sqrt{\nu'}$.
\end{theorem}

\begin{proposition}[Choice of the proposal parameter]\label{prop:tilt}
Let $v_{\max}$ and $\nu'$ be as in Definitions~\ref{def:binary}
and~\ref{def:capped} \textup{(}assume $v_{\max}\ge1$\textup{)}, let
$\Ct(t)$ be as in Definition~\ref{def:maclaurin}, and let
$\mu(t):=\E\bigl[\Bes(k,2\sqrt t)\bigr]
=\sqrt t\,I_{k+1}(2\sqrt t)/I_k(2\sqrt t)$.  Then
\begin{equation}\label{eq:tiltderiv}
t\,\frac{d}{dt}\log\Ct(t)\;=\;\mu(t)-v_{\max}\,\mathbf 1_{\{t<\nu'\}},
\end{equation}
and $\mu$ is continuous and strictly increasing from $0$ to $\infty$.
Consequently: if $\mu(\nu')\le v_{\max}$ --- the sparse regime, since
always $\mu(\nu')\le\sqrt{\nu'}$ --- then $t=\nu'$ is the unique minimizer
of $\Ct$; if $\mu(\nu')>v_{\max}$, the unique minimizer is the
\emph{mean-matched} tilt $t^*<\nu'$ solving $\mu(t^*)=v_{\max}$, the
classical optimality condition of exponential tilting \cite{siegmund1976}.
For $t<\nu'$ the acceptance test of Algorithm~\ref{alg:pdc}
\textup{(}resp.~\ref{alg:pdccap}\textup{)} is augmented by an independent
thinning coin of success probability $(t/\nu')^{\,v_{\max}-v}$.
\end{proposition}

\begin{proposition}[A bound on block splits]\label{prop:blocksplit}
Any block split \textup{(}sample a column subset $B$ freely, accept to
correct $D_B$\textup{)} has acceptance exactly
$\Prob(D{=}k)/\max_s\Prob(D_B{=}s)
\le\frac43(4\sigma_B{+}1)\Prob(D{=}k)$, $\sigma_B^2=\mathrm{Var}(D_B)$,
and leaves a conditioned subproblem on $B$ requiring recursion.  This caps
what a block split can achieve; it is not by itself a comparison with the
acceptance of Algorithm~\ref{alg:pdc}, which \eqref{eq:pdccert} bounds
separately \textup{(}empirically the statistic split wins in the settings
of \S\ref{sec:numerics}\textup{)}.
\end{proposition}

In the sparse/slack regimes the certificate gives acceptance
$1-O(\sqrt\nu\,p^*)$ (binary) or $1-O(\delta')$ (capped): one rejection
level, no recursion, with the Bessel law an exactly samplable surrogate for
the hard half --- the certified analogue of the deterministic-second-half
constructions of \cite{arratiadesalvo2016}.

\begin{remark}[The naive tilt can fail badly]\label{rem:tiltfail}
The condition $\mu(\nu')\le v_{\max}$ is not automatic, and its failure
compounds super-exponentially.  The minimal seed is $k=0$, $c=1$,
$p=q=0.9$: $\nu'=81$, $\mu(\nu')\approx8.75$ against $v_{\max}=1$, and the
mean-matched $t^*\approx1.67$ improves the certified constant by
$\Ct(\nu')/\Ct(t^*)\approx3.7\times10^4$, lifting the acceptance from
$1.3\times10^{-5}$ to $0.48$.  With equal odds $\rho=p/(1{-}p)>1$ the
mismatch persists at every $c$, since $\mu(\nu')\approx\rho c>c=v_{\max}$,
and the ratio explodes: at $c=20$, $p=q=0.9$, $k=0$,
$\Ct(\nu')/\Ct(t^*)>10^{100}$.  Two lessons follow.  First, the
mean-matched tilt is not optional outside the sparse regime.  Second, the
dense regime is also where the Maclaurin certificate itself turns
conservative: in the $c=20$ example the certified acceptance at $t^*$ is
$\approx10^{-15}$ while the true acceptance of the same proposal is
$\approx0.28$ --- the certificate is a sparse-regime instrument, in the
sense of Remark~\ref{rem:nontrivial}.  All sparse settings of
\S\ref{sec:numerics} satisfy $\mu(\nu')\le v_{\max}$ comfortably
\textup{(}e.g.\ $c=300$: $\mu(\nu')\approx24$ against
$v_{\max}=299$\textup{)}, so $t=\nu'$ is optimal there.
\end{remark}

\begin{remark}[A certified counting estimator]\label{rem:counting}
The acceptance indicator of one attempt is Bernoulli with mean exactly
$\Prob(D{=}k)/\Ct$, and $\Ct$ is closed-form; so $\Ct$ times the empirical
acceptance over $n$ attempts is an unbiased estimator of $\Prob(D{=}k)$ ---
equivalently of the weighted table counts of \S\ref{sec:enum} --- with
certified relative standard error $\le\sqrt{(1{-}\alpha_0)/(\alpha_0 n)}$ for
any certified lower bound $\alpha_0\le\Prob(D{=}k)/\Ct$.  Sampling and
enumeration are served by the same certificate.
\end{remark}

\begin{proposition}[Both margins, two rows]\label{prop:bothmargins}
Fix row margins $(u_1,u_2)$ and column margins
$(c_1,\dots,c_c)\in\{0,1,2\}^c$ with $\sum_jc_j=u_1+u_2$, and weight each
compatible zero--one table by the product of its cells' odds
\textup{(}Definition~\ref{def:binary}\textup{)}.  Write $A=\{j:c_j=1\}$, $f=\#\{j:c_j=2\}$, and
$s=u_1-f$.  Columns with $c_j\in\{0,2\}$ are deterministic; the ensemble is
nonempty iff $0\le s\le|A|$, and then it is the law of the subset
$S\subseteq A$, $|S|=s$, of columns whose $1$ lies in row $1$, with
\[
\Prob(S)\ \propto\ \prod_{j\in S}\frac{x_j}{y_j}
\]
--- Fisher's multivariate noncentral hypergeometric law, uniform over
tables exactly when the odds ratios $x_j/y_j$ are constant.  It is exactly
samplable in $O(|A|\,s)$ arithmetic operations by sequential conditioned
Bernoulli sampling through the elementary-symmetric recursion
\cite{chendempsterliu1994} --- equivalently, probabilistic
divide-and-conquer with deterministic second half \cite{desalvo2018}.  In
particular the $r=2$ case of Open Problem~\ref{prob:margins} is exactly
solvable, which is why that problem begins at $r=3$.
\end{proposition}

Proofs: \S\ref{sec:pf-pdc}.

\section{Beyond two rows: $r\times c$ tables and the hyper-Bessel law}
\label{sec:mrow}

Everything above generalizes from $2$ rows to $r$ rows, with one uniform
substitution: the quadratic death rate $n(n{+}k)$ becomes the degree-$r$
polynomial $\prod_{i=1}^r(n{+}\kappa_i)$, and the paired increment of
\S\ref{sec:coupling} becomes an $r$-tuple increment --- \emph{one added point
per row} --- which again preserves the conditioning, since adding $1$ to
every row sum leaves every successive difference unchanged.  This section
states the extension; complete proofs are in \S\ref{sec:pf-mrow}.  The
$2$-row sections above are unchanged and remain the fully worked case.

\subsection{The hyper-Bessel target}

Fix $r\ge2$ rows.  Condition the row sums $R_1,\dots,R_r$ on all successive
differences: $D_i:=R_i-R_{i+1}=k_i$ for $i=1,\dots,r{-}1$, with
$k_i\in\Z_{\ge0}$ (nonincreasing rows; the general case permutes rows), and
set the offsets
\begin{equation}\label{eq:kappadef}
\kappa_i:=k_i+k_{i+1}+\cdots+k_{r-1}\quad(1\le i\le r{-}1),\qquad
\kappa_r:=0,
\end{equation}
so that on the fiber $R_i=R_r+\kappa_i$.  Write
\[
d_\kappa(n):=\prod_{i=1}^r(n+\kappa_i),
\qquad
\mathfrak J_\kappa(w):=\sum_{n\ge0}\frac{w^{\,n}}{\prod_{i=1}^r(n+\kappa_i)!}\,,
\]
the latter a hyper-Bessel function in the sense of Delerue \cite{delerue1953} (see
\cite{kiryakova2010} for the modern multi-index theory): for $r=2$, $\mathfrak J_{(k,0)}(w)=w^{-k/2}I_k(2\sqrt w)$.

\begin{definition}\label{def:hypbes}
The \emph{hyper-Bessel distribution} $\mathrm{HB}(\kappa,w)$, $w>0$, has mass
function
\begin{equation}\label{eq:hbpmf}
p_n=\frac1{\mathfrak J_\kappa(w)}\,
\frac{w^{\,n}}{\prod_{i=1}^r(n+\kappa_i)!},\qquad n\ge0,
\end{equation}
with ratio $p_{n+1}/p_n=w/d_\kappa(n{+}1)$, strictly decreasing; hence
$\mathrm{HB}(\kappa,w)$ is strictly log-concave with mode
$m^*=\max\{n:d_\kappa(n)\le w\}$, and $\E\bigl[d_\kappa(N)\bigr]=w$ exactly.
\end{definition}

\begin{proposition}[Conditional representation and fiber probability]
\label{prop:hbrep}
Let $R_i\sim\Poi(\Lambda_i)$, $1\le i\le r$, be independent,
$\Lambda_i>0$; set $D_i:=R_i-R_{i+1}$ and
$\kappa_i:=k_i{+}\cdots{+}k_{r-1}$, $\kappa_r:=0$, with
$\mathfrak J_\kappa$ and $\mathrm{HB}$ as in \eqref{eq:hbpmf}.  Then
\begin{equation}\label{eq:hbfiber}
\Prob\bigl(D_i=k_i,\ 1\le i\le r{-}1\bigr)
=e^{-\sum_i\Lambda_i}\,\Bigl(\prod_{i=1}^{r}\Lambda_i^{\kappa_i}\Bigr)\,
\mathfrak J_\kappa(w),\qquad w:=\prod_{i=1}^r\Lambda_i,
\end{equation}
and $\mathcal L\bigl(R_r\mid D_i=k_i\ \forall i\bigr)=\mathrm{HB}(\kappa,w)$,
depending on $(\Lambda_i)$ only through the product $w$.  For $r=2$ this is
Lemma~\ref{lem:condrep}.
\end{proposition}

\subsection{Operator, factors, fixed point, approximation theorem}

The operator machinery of \S\ref{sec:operator} was proved from two inputs
only --- the ratio identity and log-concavity via Lemma~\ref{lem:hazard} ---
and both persist.  Write $\Lambda_{\kappa,w},\Gamma_{\kappa,w}$ for the
factors of the target \eqref{eq:hbpmf}.

\begin{theorem}[Hyper-Bessel operator]\label{thm:hboperator}
Let $W$ be $\Z_{\ge0}$-valued with $\E W^r<\infty$, let
$d_\kappa(n)=\prod_{i=1}^r(n{+}\kappa_i)$, and let $\mathrm{HB}(\kappa,w)$ be the
law \eqref{eq:hbpmf}.  Then $W\sim\mathrm{HB}(\kappa,w)$ iff
$\E[\,w\,g(W{+}1)-d_\kappa(W)\,g(W)\,]=0$ for all bounded $g:\Z_{\ge1}\to\R$
\textup{(}$g(0):=0$, never evaluated\textup{)}.  The Stein equation
$w\,g(n{+}1)-d_\kappa(n)g(n)=\tilde h(n)$ has the unique bounded solution
\eqref{eq:solution} with $\nu$ replaced by $w$.  The factors
$\Lambda_{\kappa,w}:=\sup_A\sup_{m\ge1}|\Delta g_A(m)|$ and
$\Gamma_{\kappa,w}:=\sup_A\sup_{m\ge0}|g_A(m{+}1)|$, for $g_A$ the
solution at $h=\mathbf 1_A$, obey, with $m^*$ the mode \textup{(}the
largest $n$ with
$d_\kappa(n)\le w$\textup{)} and $\Delta_{m^*}:=d_\kappa(m^*)-d_\kappa(m^*{-}1)$,
\begin{equation}\label{eq:hbfactors}
\begin{gathered}
\Lambda_{\kappa,w}\le\min\Bigl\{\frac1w+\frac1{d_\kappa(1)},\ \frac3w\Bigr\},\\[5pt]
\Gamma_{\kappa,w}\le
\begin{cases}
\ 2/\Delta_{m^*}, & m^*\ge1,\\[2pt]
\ 2/d_\kappa(1), & m^*=0 .
\end{cases}
\end{gathered}
\end{equation}
Moreover $d_\kappa(n{+}1)-d_\kappa(n)\ge r\,d_\kappa(n)^{(r-1)/r}$, so
$\Gamma_{\kappa,w}=O\bigl(w^{-(r-1)/r}\bigr)$ in the bulk.
\end{theorem}

\begin{theorem}[Fixed point and approximation theorem]\label{thm:hbmaster}
For $W$ with $\E W^r<\infty$ and $\beta:=\E[d_\kappa(W)]>0$
\textup{(}$d_\kappa(n)=\prod_i(n{+}\kappa_i)$ as in
Theorem~\ref{thm:hboperator}\textup{)}, define the
\emph{$d_\kappa$-bias} $W^{\boxplus}$ by
$\Prob(W^{\boxplus}{=}n)=d_\kappa(n)\Prob(W{=}n)/\beta$, $n\ge1$.  Then
$W^{\boxplus}\overset{d}{=}W{+}1$ iff $W\sim\mathrm{HB}(\kappa,\beta)$, and
for any coupling and any $w>0$,
\begin{equation}\label{eq:hbmaster}
\dtv\bigl(W,\mathrm{HB}(\kappa,w)\bigr)
\le w\,\Lambda_{\kappa,w}\;\E\bigl|W{+}1-W^{\boxplus}\bigr|
\;+\;\Gamma_{\kappa,w}\,|\beta-w| ,
\end{equation}
with $\Lambda_{\kappa,w},\Gamma_{\kappa,w}$ the factors of
\eqref{eq:hbfactors}; the matched choice $w=\beta$ kills the second term.
Equivalently, the fixed points of the centered $d_\kappa$-bias transform
$W\mapsto W^{\boxplus}{-}1$ are exactly the hyper-Bessel family
$\{\mathrm{HB}(\kappa,w):w>0\}$, and the matched bound reads
\begin{equation}\label{eq:hbmatched}
\dtv\bigl(W,\mathrm{HB}(\kappa,\beta)\bigr)\ \le\
\beta\,\Lambda_{\kappa,\beta}\;\dw\bigl(W,\,W^{\boxplus}{-}1\bigr).
\end{equation}
\end{theorem}

The coupling reading of \S\ref{sec:coupling} sharpens: on the fiber,
$R_i=W+\kappa_i$, so $d_\kappa(W)=\prod_iR_i$ --- the product of \emph{all} row
sums --- and the $d_\kappa$-bias is realized by size-biasing every row at once, an
$r$-tuple of point additions that shifts each $R_i$ by $+1$ and hence every
difference by $0$.

\subsection{The dependency expansion}

\begin{definition}\label{def:mcols}
Let $Z_j=(X^{(1)}_j,\dots,X^{(r)}_j)\in\Z_{\ge0}^r$, $j\le c$, be independent
across $j$ with $\E\|Z_j\|^{r+1}<\infty$ (an $(r{+}1)$-st moment, so that
every $S$-tilted column below has finite mean and finite Wasserstein cost;
costs are $+\infty$ otherwise and the bounds hold trivially; terms whose
moment $\chi^{(S)}_j$ vanishes are absent from the corresponding sums, and
their tilts are not needed); $U_i=\sum_jX^{(i)}_j$,
$N_i=\E U_i$, fiber $\{D_i=k_i\ \forall i\}$ assumed of positive
probability, $W\sim\mathcal L(U_r\mid\text{fiber})$, and
\begin{equation}\label{eq:nug}
\nu:=\E\Bigl[\prod_{i=1}^rU_i\Bigr]
=\sum_{\vec\jmath\in[c]^r}\ \prod_{j\,:\,S_j\ne\emptyset}\chi^{(S_j)}_j,
\qquad
\chi^{(S)}_j:=\E\Bigl[\prod_{i\in S}X^{(i)}_j\Bigr],
\end{equation}
where $S_j=\{i:\jmath_i=j\}$ is the block of rows the tuple sends to column $j$.  For $\emptyset\ne S\subseteq[r]$ define the $S$-tilted column
$dZ_j^{(S)}\propto\bigl(\prod_{i\in S}x_i\bigr)d\Prob_j$, the cost
$\alpha^{(S)}_j:=\dw\bigl(Z^{(S)}_j,\,Z_j+e_S\bigr)$ \textup{(}$e_S$ the
indicator vector of $S$, $\ell^1$ metric on $\Z^r$\textup{)}, the aggregates
and mixed moments
\begin{equation}\label{eq:ASdef}
A^{(S)}:=\sum_{j=1}^c\chi^{(S)}_j\,\alpha^{(S)}_j,
\qquad
\nu_T:=\E\Bigl[\prod_{i\in T}U_i\Bigr]\ \ (\nu_\emptyset:=1).
\end{equation}
Fiber constants: for each removed column set (of size $\le r$) with reduced
sums $U^R_r$ and difference vector $D^R\in\Z^{r-1}$, and window $\mathcal S$
the product of per-coordinate achievable-offset hulls, let $\theta$ be the
max of $\Prob(D^R{=}\vec s)/\Prob(\text{fiber})$ over $\vec s\in\mathcal S$,
$\tau_1$ the max of
$|\Prob(D^R{=}\vec s{+}e_i)-\Prob(D^R{=}\vec s)|/\Prob(\text{fiber})$ over
coordinate directions, and $\omega$ the max of
$\dw(\mathcal L(U^R_r{\mid}D^R{=}\vec s{+}e_i),\mathcal L(U^R_r{\mid}
D^R{=}\vec s))$.
\end{definition}

\begin{theorem}[Dependency expansion, $r$ rows]\label{thm:mexpansion}
Under Definition~\ref{def:mcols}, with
$\Lambda_{\kappa,\nu},\Gamma_{\kappa,\nu}$ the factors \eqref{eq:hbfactors}
at $w=\nu$, the costs vanish iff every column is an
independent Poisson $r$-tuple, and
\begin{equation}\label{eq:mexpansion}
\dtv\bigl(W,\ \mathrm{HB}(\kappa,\nu)\bigr)\ \le\
\bigl[\theta(1{+}2\omega)\Lambda_{\kappa,\nu}
+2\tau_1\Gamma_{\kappa,\nu}\bigr]
\sum_{\emptyset\ne S\subseteq[r]}\nu_{S^c}\,A^{(S)} .
\end{equation}
For $r=2$ this recovers Theorem~\ref{thm:expansion}, up to the factor
$2$ in the constants, with
$\nu_{\{2\}}A^{(\{1\})}+\nu_{\{1\}}A^{(\{2\})}+A^{(\{1,2\})}
= N_VA+N_UB+C$; the factor is absent there because each row meets only
one difference coordinate.
\end{theorem}

The factor $2$ arises because an interior row enters two successive
differences, so a coupling discrepancy of size $\varepsilon$ in one row moves
the difference vector by at most $2\varepsilon$ in $\ell^1$; the telescoping
walks a staircase lattice path inside the window.

\subsection{Enumeration and sampling}

The table pipeline generalizes wholesale.  Poissonizing an $r\times c$
table \cite{birch1963,abt2003}
with $T_{ij}\sim\Poi(\lambda_i\mu_j)$ and $M=\sum_j\mu_j$: (a) the fiber
probability is
\eqref{eq:hbfiber} with $\Lambda_i=\lambda_iM$; (b) the free row total $R_r$
is $\mathrm{HB}(\kappa,w)$; (c) given all row totals the rows are independent
$\Mult(r_i;\mu/M)$; (d) the column-sum stage produces the Fisher--Yates weighted count
\cite{fisher1935,yates1934,agresti1992}
$\prod_ir_i!\;M^{-\sum_ir_i}\sum_{T\in\mathcal T(r,c)}
\prod_j\mu_j^{T_{+j}}\prod_{ij}1/T_{ij}!$, with $T_{+j}$ the $j$th column
total, exactly as in the $2$-row case.  For
sampling, the exact envelope persists at the tilted parameter
$w'=\prod_{i=1}^rS_i$, where in the binary model the column odds are
$x^{(i)}_j:=p_{ij}/(1{-}p_{ij})$ and $S_i:=\sum_jx^{(i)}_j$, and in the
capped model $S_i:=N_i=\sum_j\lambda_{ij}$:

\begin{algorithm}[\textsc{Bessel--PDC}: $r$ rows]\label{alg:mpdc}
\leavevmode\smallskip

\noindent\rule{\textwidth}{0.8pt}
\emph{Input:} per-row column weights (binary: $p_{ij}$, giving odds
$x^{(i)}_j=p_{ij}/(1{-}p_{ij})$; capped: rates $\lambda_{ij}$ and caps
$\kappa_j$); differences $k_1,\dots,k_{r-1}\in\Z_{\ge0}$.

\smallskip\noindent
\emph{Setup} (once):
\begin{enumerate}[label=\textup{(S\arabic*)},leftmargin=2.9em,itemsep=1pt,topsep=2pt]
\item Offsets $\kappa_i=k_i+\cdots+k_{r-1}$ ($\kappa_r=0$); row sums
$S_i=\sum_jx^{(i)}_j$ (binary) or $S_i=N_i=\sum_j\lambda_{ij}$ (capped);
$w'=\prod_iS_i$.
\item Tabulate $\mathrm{HB}(\kappa,w')$ by $t_0=1$,
$t_{v+1}=t_v\cdot w'/\prod_{i=1}^r(v{+}1{+}\kappa_i)$, truncate, normalize,
store the cdf.
\item Sampling tables for the $r$ column distributions
$x^{(i)}/S_i$ (resp.\ $\lambda^{(i)}/N_i$).
\end{enumerate}
\emph{Repeat} until an attempt accepts:
\begin{enumerate}[label=\textup{(\arabic*)},leftmargin=2.9em,itemsep=1pt,topsep=2pt]
\item Draw $v\sim\mathrm{HB}(\kappa,w')$ by inversion.
\item For each row $i=1,\dots,r$: draw $v{+}\kappa_i$ indices i.i.d.\ from
$x^{(i)}/S_i$ (binary), rejecting the attempt on any within-row repeat; or
draw row $i\sim\Mult(v{+}\kappa_i;\lambda^{(i)}/N_i)$ (capped), rejecting
on any cap violation.  Repeats/co\-incidences \emph{across} rows are
allowed.
\item \textbf{Accept}: output the $r$ index sets (binary) or the table
(capped).
\end{enumerate}
\rule{\textwidth}{0.8pt}
\end{algorithm}

\begin{theorem}[Bessel--PDC, $r$ rows: exactness and acceptance]\label{thm:mpdc}
Assume all entries $X^{(i)}_j$ are independent across both rows and columns
\textup{(}Definition~\ref{def:mcols} permits dependence within a column for
the expansion \eqref{eq:mexpansion}; the sampler requires full independence,
since $f$ below factorizes over rows and the dart rounds are
independent\textup{)}.  With $f(v)=\prod_i\Prob(U_i{=}v{+}\kappa_i)$, row odds
$x^{(i)}_j=p_{ij}/(1{-}p_{ij})$, $S_i=\sum_jx^{(i)}_j$, $w'=\prod_iS_i$,
$N_i=\E U_i$, $F^{(i)}_j$ as in Definition~\ref{def:capped}
\textup{(}row index $i\in[r]$\textup{)}, and $b'$ the
$\mathrm{HB}(\kappa,w')$ mass function \eqref{eq:hbpmf}:
\begin{enumerate}[leftmargin=1.6em,label=\textup{(\alph*)},itemsep=1pt]
\item \textbf{Binary:} $f(v)\le\Ct\,b'(v)$ with
$\Ct=\prod_{i,j}(1{-}p_{ij})\cdot\prod_iS_i^{\kappa_i}\cdot
\mathfrak J_\kappa(w')$, the ratio being the probability that $r$
independent dart rounds \textup{(}$v{+}\kappa_i$ draws from
$x^{(i)}/S_i$\textup{)} are each collision-free.
Algorithm~\ref{alg:mpdc} (binary branch) is exact for the weighted binary
table measure on the fiber, with acceptance
$\Prob(\textup{fiber})/\Ct\ \ge\ \bigl(\Sk^{(r)}-\delta_r\bigr)/\Ct$, where
$\Sk^{(r)}$ is \eqref{eq:hbfiber} at $\Lambda_i=N_i$ and
$\delta_r=\sum_i(1\wedge N_i^{-1})\sum_jp_{ij}^2$.
\item \textbf{Capped Poisson:} Algorithm~\ref{alg:mpdc} (capped branch) is
exact, with $\Ct=\Sk^{(r)\prime}/\prod_{i,j}F^{(i)}_j$, and certificate
$(1-\delta'/\Sk^{(r)\prime})\prod_{i,j}F^{(i)}_j$,
$\delta'=\sum_{i,j}\bigl(1{-}F^{(i)}_j\bigr)$.
\end{enumerate}
\end{theorem}

\begin{proposition}[Choice of the proposal parameter, $r$ rows]\label{prop:mtilt}
With
$\mu_\kappa(t):=\E[\mathrm{HB}(\kappa,t)]
=t\,\mathfrak J_\kappa'(t)/\mathfrak J_\kappa(t)$, $v_{\max}$ the largest
achievable $v$, and $\Ct(t)$ the certified envelope constant for the
proposal $\mathrm{HB}(\kappa,t)$ \textup{(}the $r$-row analogue of
Definition~\ref{def:maclaurin}\textup{)},
\[
t\,\frac{d}{dt}\log\Ct(t)=\mu_\kappa(t)-v_{\max}\,\mathbf 1_{\{t<w'\}};
\]
$\mu_\kappa$ is continuous and strictly increasing from $0$ to $\infty$, so
$t=w'$ is the unique optimum iff $\mu_\kappa(w')\le v_{\max}$, and otherwise
the unique optimum is the mean-matched $t^*<w'$ with
$\mu_\kappa(t^*)=v_{\max}$, exactly as in Proposition~\ref{prop:tilt}.
\end{proposition}

\begin{corollary}[Fiber-mass factorization, $r$ rows]\label{cor:menum}
In the binary setting of Theorem~\ref{thm:mpdc}, write
$E^{(i)}_w:=\Prob(U_i{=}w)$ for the row Poisson--binomial masses,
$b^{\circ}$ for the $\mathrm{HB}(\kappa,w^{\circ})$ mass function
\eqref{eq:hbpmf} at $w^{\circ}=\prod_iN_i$, and $\Sk^{(r)}$, $\delta_r$ as
in Theorem~\ref{thm:mpdc}\textup{(a)}.  For independent Poisson rows the
factorization is exact: for every $v\ge0$,
$\prod_i\Prob(R_i{=}v{+}\kappa_i)=\Sk^{(r)}\,b^{\circ}(v)$, by
Proposition~\ref{prop:hbrep}.  For binary rows, the proof of
Theorem~\ref{thm:enum} transports verbatim and gives, for every $v\ge0$,
\[
\Bigl|\,\prod_{i=1}^rE^{(i)}_{v+\kappa_i}-\Sk^{(r)}\,b^{\circ}(v)\Bigr|
\ \le\ \delta_r\,b^{\circ}(v)+\bigl(\Sk^{(r)}+\delta_r\bigr)\varepsilon^*_r,
\]
where $\varepsilon^*_r$ is the total-variation bound assembled from
\eqref{eq:mexpansion} with the closed factor bounds \eqref{eq:hbfactors}.
Unlike the two-row $\varepsilon^*$ of Theorem~\ref{thm:enum}, it is fully
closed-form only once the multivariate fiber constants are certified
\textup{(}Open Problem~\ref{prob:fiber}\textup{)}.
\end{corollary}

\begin{theorem}[Independent process approximation, $r$ rows]\label{thm:mipa}
In the setting of Definition~\ref{def:mcols}, fix $B\subseteq[c]$ and
write $Z_B=(Z_j)_{j\in B}$,
$\vec D_B=\bigl(\sum_{j\in B}(X_{ij}-X_{i+1,j})\bigr)_{i=1}^{r-1}$, and
$\vec D_{B^c}=\vec D-\vec D_B$, where
$\vec D=(U_i{-}U_{i+1})_{i=1}^{r-1}$ is the difference vector conditioned
in Definition~\ref{def:mcols} and $\vec k=(k_1,\dots,k_{r-1})$.
\textup{(a)} Then
\begin{equation}\label{eq:mipa}
\dtv\bigl(\mathcal L(Z_B\mid\vec D{=}\vec k),\,\mathcal L(Z_B)\bigr)
=\dtv\bigl(\mathcal L(\vec D_B\mid\vec D{=}\vec k),\,
\mathcal L(\vec D_B)\bigr)
=\frac12\,\E\biggl|\frac{\Prob\bigl(\vec D_{B^c}{=}\vec k-\vec D_B\bigr)}
{\Prob(\vec D{=}\vec k)}-1\biggr|,
\end{equation}
the expectation over the free $\vec D_B$.
\textup{(b)} In the capped setting --- entries
$X_{ij}\sim\Poi(\lambda_{ij})$ conditioned on $X_{ij}\le\kappa_j$,
$F^{(i)}_j=\Prob\bigl(\Poi(\lambda_{ij})\le\kappa_j\bigr)$,
$\delta'_r=\sum_{i,j}(1-F^{(i)}_j)$ --- let
$\Sk^{(r)}_{\vec m}:=\Prob(\vec D^{\,0}{=}\vec m)$ be the fiber
probability of the free Poisson system, given in closed form by
Proposition~\ref{prop:hbrep} \textup{(}after shifting the offsets
$\kappa_i(\vec m)=m_i{+}\cdots{+}m_{r-1}$ so the smallest is
zero\textup{)}.  If $\Sk^{(r)}_{\vec k}>\delta'_r$, then
\begin{equation}\label{eq:mipabound}
\dtv\bigl(\mathcal L(Z_B\mid\vec D{=}\vec k),\,\mathcal L(Z_B)\bigr)
\ \le\ \frac{2\sum_{j\in B}\sum_{i=1}^r\lambda_{ij}+2\delta'_r}
{2\,\bigl(\Sk^{(r)}_{\vec k}-\delta'_r\bigr)}\,,
\end{equation}
and in the homogeneous regime $\lambda_{ij}=\lambda/c$ with caps
$\kappa_j\ge1$ and $\vec k$ fixed, every block with $|B|=o(c)$ has
$\dtv=O(|B|/c)\to0$.  The binary case is the caps-one case:
$X_{ij}\sim\mathrm{Bern}(p_{ij})$ equals $\Poi(x_{ij})$ conditioned on
$\{\le1\}$ with $x_{ij}=p_{ij}/(1{-}p_{ij})$, with per-entry defect
$1-e^{-x_{ij}}(1{+}x_{ij})\le x_{ij}^2/2$.
\end{theorem}

\begin{theorem}[Independent process approximation with row sums fixed]\label{thm:rowfixed}
In the setting of Theorem~\ref{thm:mpdc} \textup{(}entries independent
across rows and columns\textup{)}, fix $B\subseteq[c]$ and row-sum targets
$\vec u=(u_1,\dots,u_r)$ with $\Prob(\vec U{=}\vec u)>0$, where
$U_i=\sum_jX_{ij}$, $U^{(i)}_B:=\sum_{j\in B}X_{ij}$, and
$Z^{(i)}_B:=(X_{ij})_{j\in B}$.
\textup{(a)} Conditionally on $\{\vec U{=}\vec u\}$ the rows remain
independent, each conditioned on its own sum, so with
\[
T_i:=\dtv\bigl(\mathcal L(Z^{(i)}_B\mid U_i{=}u_i),\,
\mathcal L(Z^{(i)}_B)\bigr)
=\frac12\,\E\Bigl|\frac{\Prob\bigl(U^{(i)}_{B^c}{=}u_i{-}U^{(i)}_B\bigr)}
{\Prob(U_i{=}u_i)}-1\Bigr|
\]
\textup{(}the $r{=}1$ identity of Theorem~\ref{thm:ipa}, the row sum
playing the role of the difference\textup{)},
\[
\max_{1\le i\le r}T_i\ \le\
\dtv\bigl(\mathcal L(Z_B\mid \vec U{=}\vec u),\,\mathcal L(Z_B)\bigr)
\ \le\ \sum_{i=1}^rT_i .
\]
\textup{(b)} In the capped setting, with $\pi_N(m)=e^{-N}N^m/m!$,
$\delta'_i=\sum_j\bigl(1-F^{(i)}_j\bigr)$: whenever
$\pi_{N_i}(u_i)>\delta'_i$ for every $i$,
\[
\dtv\bigl(\mathcal L(Z_B\mid \vec U{=}\vec u),\,\mathcal L(Z_B)\bigr)
\ \le\ \sum_{i=1}^r
\frac{2\sum_{j\in B}\lambda_{ij}+2\delta'_i}
{2\bigl(\pi_{N_i}(u_i)-\delta'_i\bigr)}\,.
\]
In the homogeneous regime $\lambda_{ij}=\lambda/c$ with $\kappa_j\ge1$ and
each $u_i$ in the bulk \textup{(}$\pi_{N_i}(u_i)$ bounded below\textup{)},
every block with $|B|=o(c)$ has $\dtv=O(r\,|B|/c)$: each row constraint is
one-dimensional, and the constraint dimension enters only linearly in $r$.
\end{theorem}

\begin{remark}[From differences to full row sums]\label{rem:rowfixed}
Passing from the difference fiber $\{\vec D{=}\vec k\}$ to the full
row-sum fiber $\{\vec U{=}\vec u\}$ with $u_i=v_0+\kappa_i$ multiplies the
fiber probability by the single one-dimensional factor
$\Prob(V{=}v_0\mid\vec D{=}\vec k)$ --- approximately the hyper-Bessel
mass $b^{\circ}(v_0)$ of Corollary~\ref{cor:menum} --- independently of
$r$.  Operationally, Algorithm~\ref{alg:mpdc} with step \textup{(1)}
frozen at $v=v_0$ is, by the proof of Theorem~\ref{thm:mpdc}, an exact
sampler for the row-sum-fixed weighted ensemble, with per-attempt
acceptance $f(v_0)/\bigl(\Ct\,b'(v_0)\bigr)$, bounded below via the lower
bound on $f(v_0)$ from Corollary~\ref{cor:menum}.
\end{remark}

\begin{theorem}[Local Fisher--Yates approximation]\label{thm:fylocal}
Let the entries be independent Poissons with rank-one means,
$T_{ij}\sim\Poi(\lambda_i\mu_j)$, and fix column sums
$\vec c=(c_1,\dots,c_c)$ and row targets $\vec u$ with
$\sum_iu_i=\sum_jc_j=:n$.  Write $\Lambda=\sum_i\lambda_i$,
$\rho_i=\lambda_i/\Lambda$, and let $M_m$ be the $\Mult(m;\rho)$ mass
function on $\Z_{\ge0}^r$.  Conditionally on $\{\vec C{=}\vec c\}$ the
columns are independent, the $j$th column distributed as
$\Mult(c_j;\rho)$ --- exactly samplable --- and
$\vec U\mid\{\vec C{=}\vec c\}\sim\Mult(n;\rho)$.  The Fisher--Yates null
$\mathrm{FY}(\vec u,\vec c)$ is this product law conditioned further on
$\{\vec U=\vec u\}$, and for every $B\subseteq[c]$, with
$n_B=\sum_{j\in B}c_j$ and $T_B$ the block of columns in $B$,
\begin{equation}\label{eq:fylocal}
\begin{aligned}
\dtv\bigl(\mathcal L(T_B\mid\mathrm{FY}(\vec u,\vec c)),\
\mathcal L(T_B\mid\vec C{=}\vec c)\bigr)
&=\dtv\bigl(\mathcal L(\vec U_B\mid\vec U{=}\vec u,\ \vec C{=}\vec c),\
\Mult(n_B;\rho)\bigr)\\[2pt]
&=\frac12\,\E\Bigl|\frac{M_{n-n_B}\bigl(\vec u-\vec U_B\bigr)}
{M_n(\vec u)}-1\Bigr|,
\end{aligned}
\end{equation}
the expectation over $\vec U_B\sim\Mult(n_B;\rho)$: a finite sum of at
most $\binom{n_B+r-1}{r-1}$ closed-form multinomial terms.
\end{theorem}

\begin{remark}[Sampling the Fisher--Yates null through its columns]\label{rem:fysampling}
Theorem~\ref{thm:fylocal} is local by necessity: for the full table,
$\dtv\bigl(\mathcal L(T\mid\vec C{=}\vec c),\,
\mathrm{FY}(\vec u,\vec c)\bigr)=1-M_n(\vec u)$, near one, since
conditioning a law on a thin event moves essentially all of its mass.
What the theorem certifies is that any statistic supported on a set of
columns carrying a vanishing fraction of $n$ can be simulated from the
trivial column-product law with explicit total-variation error;
numerically the error tracks $n_B/n$ linearly.  Two complements.  First,
accepting the column-product sample exactly when $\vec U=\vec u$ is an
exact sampler for $\mathrm{FY}(\vec u,\vec c)$ whose acceptance is known
in closed form,
$M_n(\vec u)=n!\,\prod_i\rho_i^{u_i}\big/\prod_iu_i!$, of order
$n^{-(r-1)/2}$ for bulk targets: the certificate is free, though classical
sequential samplers are faster.  Second, the analogous certification for
\emph{weighted} tables is precisely Open Problem~\ref{prob:margins}.  At
$r=2$ with unit column sums, \eqref{eq:fylocal} recovers the classical
comparison of sampling without replacement with i.i.d.\ sampling
\cite{diaconisfreedman1980}.
\end{remark}

\subsection{What changes quantitatively}

Two honest caveats.  First, the fiber constants now live on
$\Z^{r-1}$ and are $O(1)$ only for $c\gg r$ with comparable columns:
removing up to $r$ columns must leave a non-degenerate reduced system.
Second, the constants grow with $r$ (the factor $2$, the $2^r{-}1$ cost
classes, and the collision probability of $r$ dart rounds): in the exact
$r=3$ computations of \S\ref{sec:numerics}, the transfer identity holds to
$2.8\times10^{-16}$ and every bound is valid, but the $\dtv$ slack is
$\approx230\times$ at $c=6$ (against $\approx50\times$ for $r=2$), and the
binary sampler accepts at $0.745$ (against $0.953$) --- three collision-free
rounds instead of two.  Certifying the multivariate $\theta,\tau_1,\omega$
--- Esseen-type anti-concentration for the reduced difference \emph{vector},
whose per-column characteristic factors are now genuinely multivariate ---
is the $r$-row extension of Open Problem~\ref{prob:fiber} and is exactly
quantitative multivariate LCLT territory.

\section{Numerical validation}\label{sec:numerics}

All displayed identities and the reported example values were checked in
exact or $60$-digit arithmetic; the checks validate these instances, not
the general statements.  Scripts are in the ancillary files of the arXiv
version.  The transfer
identity of Lemma~\ref{lem:transfer} holds to $3.6\times10^{-15}$ for a
random bounded $g$; the Stein solution \eqref{eq:solution} to $10^{-60}$; the
identity $\E[N(N{+}k)]=\nu$ to $10^{-57}$; the sign pattern of
Lemma~\ref{lem:signs} at every tested site; and each certificate of
Lemma~\ref{lem:cert} against exact convolutions.  Distances and reduced-fiber
constants were computed by full dynamic programming over the joint law of
$(V,D)$ and all pair-reduced systems; Monte Carlo appears only where sampling
itself is the object.

\begin{table}[h]
\centering\small
\begin{tabular}{@{}lllll@{}}
\toprule
result & system & quantity & value & \S \\
\midrule
Thm.~\ref{thm:factors} & $\Bes(k,a)$ & $\sup_A\|\Delta g_A\|$ vs.\ exact &
equal to 6 digits; $\to1/\nu$ & \ref{sec:operator}\\
Thm.~\ref{thm:main} & trunc.\ Poisson & $\dtv$ vs.\ \eqref{eq:matched} &
slack $11$--$24\times$ & \ref{sec:agg}\\
Thm.~\ref{thm:expansion} & capped, $c{=}6$ & $\dtv$ vs.\ \eqref{eq:expansion}
& slack $48$--$88\times$ as $\nu$ spans four orders of magnitude &
\ref{sec:expansion}\\
Thm.~\ref{thm:process} & binary, $c{=}6$, $k{=}1$ &
array-level vs mixing-law $\dtv$, \eqref{eq:process} &
equal to $6\cdot10^{-18}$ & \S\ref{sec:ipa}\\
Cor.~\ref{cor:binary} & binary, fixed $\nu$ & rate in $p^*$ &
$\Theta(p^*)$ exact; slack $\approx37\times$ & \ref{sec:dtv}\\
Thm.~\ref{thm:enum} & sparse diagonal & rel.\ err.\ at mode &
$\Theta(1/c)$; constants $\approx4{\times}10^2$ & \ref{sec:enum}\\
Thm.~\ref{thm:pdc} & binary $c{=}300$ & acceptance &
emp.\ $0.9530$ / exact $0.9532$ / cert.\ $0.720$ & \ref{sec:pdc}\\
Thm.~\ref{thm:pdc} & capped $c{=}12,\kappa{=}10$ & acceptance &
exact $0.9981$ / cert.\ $0.899$ & \ref{sec:pdc}\\
Prop.~\ref{prop:blocksplit} & binary $c{=}300$ & statistic vs.\ block acceptance &
$0.953$ vs.\ $0.530$ & \ref{sec:pdc}\\
Rem.~\ref{rem:tiltfail} & $k{=}0$, $c{=}1$, $p{=}q{=}0.9$ & $\Ct(\nu')/\Ct(t^*)$ & $3.7{\times}10^4$ & \ref{sec:pdc}\\
Prob.~\ref{prob:logconcave} & log-concave scan & $\max w\Lm$;
$\max w\Gm/\sigma$ & $1.0000$; $1.14$ & \ref{sec:open}\\
\addlinespace
Lem.~\ref{lem:mtransfer} & $r{=}3$ capped & transfer identity residual &
$2.8\times10^{-16}$ & \ref{sec:mrow}\\
Thm.~\ref{thm:mexpansion} & $r{=}3$ capped $c{=}6$ & $\dtv$ vs.\
\eqref{eq:mexpansion} & $2.6\times10^{-2}$ vs.\ $5.9$ ($227\times$) &
\ref{sec:mrow}\\
Thm.~\ref{thm:mpdc} & $r{=}3$ binary $c{=}40$ & envelope max ratio;
acceptance & $0.972\le1$;\ \ $0.745$ & \ref{sec:mrow}\\
\bottomrule
\end{tabular}
\caption{Synopsis of numerical validation.  ``Slack'' is
bound${}/{}$truth.  Sampler output laws match their targets to Monte Carlo
resolution (e.g.\ TV $=0.0013$ on $2\times10^5$ accepted binary samples).}
\label{tab:master}
\end{table}

\section{Open problems}\label{sec:open}

\begin{problem}[Sharp Stein factor]\label{prob:sharp}
Prove the observed $\sup_j[G(j{+}1)/\nu+F(j{-}1)/(j(j{+}k))]=(1{+}o(1))/\nu$
with the explicit analogue of $(1{-}e^{-\lambda})$, sharpening the constant
$3$ in $\Lm$.
\end{problem}

\begin{problem}[Log-concavity as the regularity behind ``two moments
suffice'']\label{prob:logconcave}
Every log-concave pmf on $\Z_{\ge0}$ is a constant-birth law: if
$p_{n+1}/p_n$ is decreasing, then $p=\mu_{d,w}$ with increasing death rate
$d(n):=w\,p_{n-1}/p_n$, $d(0)=0$, for the operator
$w\,g(n{+}1)-d(n)\,g(n)$, and the structural layer of this paper holds
verbatim: $\E\,d(N)=w$; $d$-biasing $\mu_{d,w}$ is the unit shift, since
$d(n)p_n=w\,p_{n-1}$; and the fixed points of the $d$-bias-minus-one map
are exactly $\{\mu_{d,w'}\}_{w'>0}$, by the argument of
Proposition~\ref{prop:fixed}.  We conjecture that the quantitative layer is
universal as well: for every log-concave target,
\[
\Lm\;\le\;\frac1w
\qquad\text{and}\qquad
\Gm\;\le\;C\,\frac{\max\{1,\sigma\}}{w}
\]
for an absolute constant $C$, with $\sqrt{2\pi}/4$ the bulk value of
$w\Gm/\sigma$.  The $\max$ cannot be removed.  Take $w=1$, $q=1/M$,
$D=M^2{+}M{-}1$, and
\[
p_0=\frac{M-1}{D},\qquad
p_{n+1}=\frac{M(M-1)}{D}\,q^{\,n}\quad(n\ge0).
\]
The ratios $(M,q,q,\dots)$ are nonincreasing, so $p$ is log-concave with
death rate $d=(0,\,1/M,\,M,\,M,\dots)$; yet
$w\Gm\ge F(0)G(1)/p_0=1-p_0\to1$ while
$\sigma^2=M^3(2M^2{-}1)\big/\bigl((M{-}1)^2D^2\bigr)\to0$, so
$w\Gm/\sigma\ge\sqrt{M/8}$.\footnote{This family was posted, in response
to the first arXiv version of this paper, in a public discussion thread
\textup{(}\texttt{github.com/the-omega-institute/trureturing}, issue
26054\textup{)}; we record it here with independent verification.  The
first version conjectured $\Gm\le C\sigma/w$.}  Near-deterministic laws
are the only obstruction known to us, and for $d=d_\kappa$ we know no
violation of the stronger $\sigma$-form.  For $d=d_\kappa$ our proofs
certify
$\Lm\le\min\{1/w{+}1/d(1),3/w\}$, while the certified uniform bound
\eqref{eq:gamma} is weaker than $\sigma/w$ and degenerates as $d$ flattens
toward the geometric boundary; the general lemmas of
\cite{brownxia2001,er2008} require the per-capita rate $d(n)/n$ to be
nondecreasing.  Yet the displayed inequalities appear to hold beyond both
restrictions: numerically, $w\Lm\le1.0000$ and $w\Gm/\sigma\le1.14$ across
$d(n)=n^{\alpha}$ for $\alpha\in\{0.3,0.5,1,2,3\}$ \textup{(}sub-linear
$\alpha$ lies outside the \cite{brownxia2001,er2008} hypotheses\textup{)},
geometric profiles with $w\Lm=w$ exactly, and an adversarial search over
randomly generated increasing $d$, which saturates at $w\Lm=1.0000$
without exceeding it \textup{(}\texttt{22\_logconcave\_conjecture.py}, ancillary
files\textup{)}.  Prove the factor bounds from log-concavity alone; then
characterize the death rates whose $d$-weight localizes over independent
blocks --- the prerequisite for an AGG-type expansion, met by
$d(n)=\prod_i(n{+}\kappa_i)$ through Poisson vectors conditioned on
successive differences \textup{(}\S\ref{sec:mrow}\textup{)} --- so that
``two moments suffice'' becomes a statement about a regularity class
rather than a list of targets.
\end{problem}

\begin{problem}[Certified fiber constants]\label{prob:fiber}
Give explicit $\theta\le\theta_0$, $\tau_1\le c/\sigma_D$,
$\omega\le\omega_0$ under Definition~\ref{def:cols} --- e.g.\ via
log-concavity of the column differences (then $D_R$ is a sum of independent
log-concave variables and Lemma~\ref{lem:hazard} applies to it) or
characteristic-function arguments of Roth--Szekeres type
\cite{rothszekeres1954}.  This removes the
$\approx\nu^{1/4}$ crude-route loss in Theorem~\ref{thm:enum}
(Remark~\ref{rem:losses}) and turns \eqref{eq:enum} into a
single-digit-slack estimate.
\end{problem}

\begin{problem}[Local upgrade for enumeration]\label{prob:local}
For arbitrary $W$ the pointwise upgrade
$\sup_n|\Prob(W{=}n)-p_n|=O(\sigma^{-1})\cdot\dtv$ is \emph{false}: a
signed, mass- and mean-preserving perturbation of $\Bes(k,a)$ supported on
three adjacent atoms leaves the local and total-variation errors of the
same order, uniformly in $\sigma$.  Identify a smoothing hypothesis on $W$
--- a Landau--Kolmogorov condition as in \cite{rollinross2015}, or
membership in the column-array models of \S\ref{sec:tables} --- under which
the upgrade holds; this is the true currency for counting, and
Lemma~\ref{lem:hazard} is the target-side tool.
\end{problem}

\begin{problem}[Both margins fixed]\label{prob:margins}
The classical table ensemble fixes the column sums as well.  Once all row
sums are fixed --- so the Bessel stage of this paper is degenerate, no
fluctuating offset remaining --- the rows are independent weighted draws
of fixed totals, and the column-sum vector $C=(C_1,\dots,C_c)$ is a sum of
$r$ independent random vectors with $\sum_jC_j$ deterministic; the
Fisher--Yates law is $\mathcal L(\text{table}\mid C=\mathbf c)$.  The
difficulty thus migrates from a scalar conditioned difference to a
$(c{-}1)$-dimensional local conditioning.  The ingredients here plausibly
extend: \textup{(i)} a Stein operator, with factors from log-concavity,
for the conditional law of one column count given the rest (a
conditioned-multinomial family replacing $\Bes$); \textup{(ii)} a
two-moment dependency expansion, as in \S\ref{sec:agg}, for the weak
negative dependence among the $C_j$; \textup{(iii)} a certified surrogate
for $\Prob(C=\mathbf c)$ sharp enough to drive an exact sampler --- in
effect a \emph{certified} sequential importance sampler \cite{cdhl2005}.
For $2\times c$ tables the both-margins ensemble is exactly samplable
\textup{(}Proposition~\ref{prop:bothmargins}\textup{)}, so the problem
begins at $r=3$.  For central weights, the column-conditioned product law
approximates the Fisher--Yates null on every small block with closed-form
error \textup{(}Theorem~\ref{thm:fylocal}\textup{)}; it is the weighted
full-table certificate that remains open.  Column-recursive
PDC \cite{desalvozhao2020} already samples both-margins tables exactly,
without certificates; the open problem is the certificate.
\end{problem}

\begin{problem}[Dependent arrays]\label{prob:dep}
Extend Theorem~\ref{thm:expansion} to dependency neighborhoods: a $b_3$
conditional-expectation term appears, and the paired local modification
remains the coupling.
\end{problem}

\begin{remark}[The geometric target]\label{rem:geometric}
Geometric entries are the route to unweighted counting: conditioned on
both margins, they induce the \emph{uniform} table ensemble.  The operator
and envelope for geometric entries are carried out in the companion paper
\cite{desalvo2026}, whose targets are the negative binomial and
hyper-negative-binomial laws.
\end{remark}

\section{Final remarks}\label{sec:final}

Each core object of this paper sits inside, or beside, an existing
framework.  The remarks below record the precise relationship in each case
--- what the general theory already delivers for the Bessel law, what it
delivers in sharper form than our closed constants, and what does not
follow from it.  We prefer this explicit accounting to a novelty claim in
the large: each item can be audited separately.

\begin{remark}[Origin of the distribution and of the operator]\label{rem:origin}
The law \eqref{eq:pmf} was introduced and named by Yuan--Kalbfleisch
\cite{yuankalbfleisch2000}, together with the conditional representation of
Lemma~\ref{lem:condrep}; Devroye \cite{devroye2002} gave exact samplers,
and the difference law itself is Skellam's \cite{skellam1946}.  The
operator \eqref{eq:operator} has a partial precedent: in an unpublished
framework for approximating conditional distributions, Chiarini, Cipriani,
and Conforti \cite{ccc2017} derive, as their illustrative example, the
characterizing operator of two independent Poissons conditioned on the
diagonal $\{X{-}Y{=}0\}$ --- a birth--death operator with birth rate
$\lambda_1\lambda_2$ and death rate $n^2$, which is \eqref{eq:operator} at
$k=0$.  Their framework builds conditional operators from event-preserving
transformations and bounds distances by comparing two such operators
globally; the general-$k$ operator, the factors of
Theorem~\ref{thm:factors}, and the per-column expansion of
\S\ref{sec:expansion} do not appear there.
\end{remark}

\begin{remark}[Relation to the general frameworks]\label{rem:general}
The Bessel law is a discrete Gibbs measure in the sense of
\cite{er2008}: $p_n\propto e^{V(n)}\nu^n/n!$ with $V(n)=-\log\bigl((n{+}k)!\bigr)$,
which is concave.  This alone is not informative --- that class exhausts all
laws on $\Z_{\ge0}$ --- but the associated birth--death chain with unit per
capita death rates has birth rates
\[
b_n=(n{+}1)\,p_{n+1}/p_n=\nu/(n{+}k{+}1),
\]
\emph{strictly decreasing}, which is precisely condition \textup{(C4)} of
\cite{brownxia2001} and the hypothesis of \cite[Remark 2.2]{er2008}.  Writing
$\widehat g(n):=(n{+}k)g(n)$ turns \eqref{eq:operator} into the unit per
capita operator $\widehat g(n{+}1)b_n-n\widehat g(n)$ of \cite{er2008}, so the
two Stein equations coincide and \cite[Lemma 2.1]{er2008} --- their statement
of \cite[Lemma 2.4, Theorem 2.1]{brownxia2001} --- applies verbatim, giving
$\sup_f|\Delta\widehat g_f(j)|=\frac{j+k+1}{\nu}G(j{+}1)+\frac{F(j-1)}{j}$
and hence $\le\frac1j\wedge\frac{j+k+1}\nu$.  Transferring back through
$g=\widehat g/(n{+}k)$ yields a bound on $\sup_A\|\Delta g_A\|$ roughly
$1.4/\nu$ for large $\nu$: sharper than the closed constant $3/\nu$ of
\eqref{eq:lambda} by a factor $\approx2$, though weaker than the exact
identity in \eqref{eq:lambda} itself (at $k=2$, $\nu=25$ the three values are
$0.0398$ exact, $0.0560$ from the general route, $0.1200$ from $3/\nu$).  The
same route bounds $\sup_A\|g_A\|$ by $0.149$ against the exact $0.0438$ and
the closed $0.278$ of \eqref{eq:gamma}.  What does \emph{not} transfer is a
uniform bound on the solution: \cite[Lemma 2.5]{er2008} requires
$\lambda_1=\inf_nb_n>0$, whereas $b_n\downarrow0$ here, and their
\cite[Lemma 2.4]{er2008} degrades like $1/G(j)$.  So
\eqref{eq:gamma} is not available from the general theory, while the closed
constants in \eqref{eq:lambda}--\eqref{eq:gamma} can be improved by
$\approx2$ by routing through \cite{er2008,brownxia2001}; we keep the
self-contained log-concavity derivation, which yields both exact identities
in one argument, and use the exact quantities in all numerical work.  For the
hyper-Bessel law the same applies whenever some $\kappa_i=0$ --- true under
the normalization $\kappa_r=0$ of \S\ref{sec:mrow}, since then
$b_n=w/\prod_{i<r}(n{+}1{+}\kappa_i)$ decreases --- but not in general
\textup{(}e.g.\ $\kappa=(5,2)$, $w=30$ has $b_0<b_1$\textup{)}.
\end{remark}

\begin{remark}[The transform within the bias-transform family]\label{rem:transform}
The quadratic bias of Definition~\ref{def:qbias} is an instance of the
polynomial-weight distributional transforms of Goldstein--Reinert
\cite{goldsteinreinert2005} (see also D\"obler \cite{dobler2017}); neither
the weight $n(n{+}k)$ nor the Bessel fixed-point set of the centered map
$\mathcal T_k$ appears there.  $\mathcal T_k$ thus joins the short list of
transforms whose fixed points form a classical family: zero-bias (centered
normal) \cite{goldsteinreinert1997}, size-bias-minus-one (Poisson), and
equilibrium (exponential) \cite{pekozrollin2011}.
\end{remark}

\begin{remark}[The expansion against its lineage]\label{rem:agg}
Two-moment expansions in the mold of Arratia--Goldstein--Gordon
\cite{agg1989} exist for several unconditioned targets: compound Poisson
\cite{bcl1992}, negative binomial \cite{brownphillips1999}, translated
Poisson \cite{rollin2007}, and --- closest to the present setting --- the
Skellam law of the unconditioned difference \cite{gan2018}.  Conditional
targets have been reached in two adjacent senses: conditioning on a
threshold event in a compound-Poisson model \cite{ganxia2015}, and a
conditional central limit theorem given an asymptotically jointly Gaussian
statistic \cite{deyterlov2021}.  We are not aware of a two-moment expansion
whose target is the law of a statistic given a linear constraint, as in
Theorem~\ref{thm:expansion}.  Three of its ingredients deserve individual
attribution.  \emph{The transfer identity.}  A generator comparison in the
style of \cite[Theorem~3.5]{er2008}, taking the second law to be
$\mathcal L(V\mid U{-}V{=}k)$, is formally available but inapplicable as
stated: it requires the conditional birth rates
$(n{+}1)\Prob(V{=}n{+}1\mid D{=}k)/\Prob(V{=}n\mid D{=}k)$, which are
unknown in closed form precisely when the theorem is wanted.
Lemma~\ref{lem:transfer} is the device that expresses this operator action
through per-column moments instead; we believe it is new.  \emph{The
costs.}  Coupling costs between a law and its size-biased version are
classical \cite{bhj1992,goldsteinrinott1996,er2008}; what appears new is
the product tilt $xy\,d\Prob_j/\chi_j$ and the per-column, moment-weighted
aggregation $N_VA+N_UB+C$ against a conditional target.  \emph{The paired
increment.}  Moves preserving a conditioning event have two clear
precedents: the event-preserving transformations from which \cite{ccc2017}
assemble their conditional operator, and the Markov-basis moves of
\cite{diaconissturmfels1998} preserving table margins; their use inside a
total-variation moment expansion appears to be new.
\end{remark}

\bigskip
\begin{center}\rule{0.5\textwidth}{0.4pt}\\[4pt]
\textsc{Part II. Proofs}\end{center}

\section{Proofs}\label{sec:proofs}

Proofs are grouped by theme, one subsection per block of Part~I.

\subsection{The operator and Stein factors}\label{sec:pf-operator}

\begin{proof}[Proof of Lemma~\ref{lem:basic}]
The series $I_k(a)=\sum_{n\ge0}(a/2)^{2n+k}/(n!(n+k)!)$ converges for all
$a>0$, so \eqref{eq:pmf} is a probability law.  From \eqref{eq:ratio},
$p_{n+1}p_{n-1}/p_n^2=\frac{n(n+k)}{(n+1)(n+k+1)}<1$, giving strict
log-concavity; and $p_n\ge p_{n-1}$ iff $n(n{+}k)\le\nu$, so the mode is the
largest such $n$, i.e.\ $\lfloor(\sqrt{k^2+4\nu}-k)/2\rfloor$, the integer
part of the positive root of $m(m{+}k)=\nu$.  The terms of \eqref{eq:pmf}
decay superexponentially, giving all moments.  For the mean: for $n\ge1$,
\[
n\,p_n=\frac1{I_k(a)}\,\frac{(a/2)^{2n+k}}{(n-1)!\,(n+k)!}
=\frac a2\cdot\frac1{I_k(a)}\,\frac{(a/2)^{2(n-1)+(k+1)}}{(n-1)!\,\bigl((n-1)+(k+1)\bigr)!},
\]
and summing over $n\ge1$ gives $\E N=\frac a2\,I_{k+1}(a)/I_k(a)$.  For the
quadratic moment, \eqref{eq:ratio} in the form $n(n{+}k)p_n=\nu p_{n-1}$
gives $\E[N(N{+}k)]=\sum_{n\ge1}n(n{+}k)p_n=\nu\sum_{n\ge0}p_n=\nu$.
\end{proof}

\begin{proof}[Proof of Lemma~\ref{lem:condrep}]
For $n\ge0$, by independence,
\[
\Prob(Y{=}n,\,X{=}n{+}k)
=e^{-\lambda_1}\frac{\lambda_1^{n+k}}{(n+k)!}\cdot
e^{-\lambda_2}\frac{\lambda_2^{n}}{n!}
=e^{-\lambda_1-\lambda_2}\,\frac{\lambda_1^{n+k}\lambda_2^{n}}{n!\,(n+k)!}.
\]
Since $\lambda_1^{n+k}\lambda_2^n=(\lambda_1\lambda_2)^{n+k/2}
(\lambda_1/\lambda_2)^{k/2}=(a/2)^{2n+k}(\lambda_1/\lambda_2)^{k/2}$, we get
\[
\Prob(Y{=}n,\,X{=}n{+}k)
=e^{-\lambda_1-\lambda_2}\Bigl(\frac{\lambda_1}{\lambda_2}\Bigr)^{k/2}
\frac{(a/2)^{2n+k}}{n!\,(n{+}k)!}\,.
\]
Summing over $n\ge0$ gives
\eqref{eq:skellampt}; dividing gives \eqref{eq:pmf}, which involves
$(\lambda_1,\lambda_2)$ only through $a$.
\end{proof}

\begin{proof}[Proof of Theorem~\ref{thm:char}]
\emph{Necessity.}  For $N\sim\Bes(k,a)$ and $g$ bounded (both expectations
finite since $\E N^2<\infty$), \eqref{eq:ratio} in the form
$\nu p_n=(n{+}1)(n{+}k{+}1)p_{n+1}$ gives
\begin{align*}
\E[\nu g(N{+}1)]&=\sum_{n\ge0}\nu p_n\,g(n{+}1)
=\sum_{n\ge0}(n{+}1)(n{+}k{+}1)\,p_{n+1}\,g(n{+}1)\\
&=\sum_{m\ge1}m(m{+}k)\,p_m\,g(m)=\E[N(N{+}k)g(N)].
\end{align*}

\emph{The solution formula and uniqueness.}  Multiply the Stein equation by
$p_n$ and use $n(n{+}k)p_n=\nu p_{n-1}$ (with $p_{-1}:=0$):
$\nu\bigl[p_n g(n{+}1)-p_{n-1}g(n)\bigr]=\tilde h(n)p_n$.  Summing over
$n=0,\dots,m$ telescopes to $\nu p_m g(m{+}1)=\sum_{n\le m}\tilde h(n)p_n$,
the first form of \eqref{eq:solution}; the second follows from
$\sum_{n\ge0}\tilde h(n)p_n=0$.  Conversely, the function defined by
\eqref{eq:solution} satisfies the recursion at every $n\ge0$ (at $n=0$ it
reads $\nu g(1)=\tilde h(0)$, consistent with the formula at $m=0$).  Any two
solutions differ by a solution of the homogeneous recursion
$g(n{+}1)=\frac{n(n{+}k)}\nu g(n)$, i.e.\ a multiple of
$g(n)=g(1)\prod_{i=1}^{n-1}\frac{i(i{+}k)}\nu$, which grows like
$(n{-}1)!\,(n{+}k{-}1)!\,\nu^{-(n-1)}/k!\to\infty$; boundedness forces
$g(1)=0$, hence uniqueness among bounded solutions.  Boundedness of the
formula for $h=\mathbf 1_A$ is Theorem~\ref{thm:factors}(b), whose proof
below does not use the present theorem.

\emph{Sufficiency.}  Fix $A\subseteq\Z_{\ge0}$ and let $g_A$ be the bounded
solution for $h=\mathbf 1_A$.  Then, using the hypothesis with $g=g_A$ and
the fact that $g_A$ solves the equation pointwise,
$0=\E[(\mathcal Ag_A)(W)]=\E[\mathbf 1_A(W)-\Prob(N{\in}A)]
=\Prob(W{\in}A)-\Prob(N{\in}A)$.  Since $A$ was arbitrary,
$W\overset{d}{=}N$.  (The expectation $\E[W(W{+}k)g_A(W)]$ is finite because
$g_A$ is bounded and $\E W^2<\infty$.)
\end{proof}

\begin{proof}[Proof of Lemma~\ref{lem:hazard}]
Log-concavity is equivalent to the ratios $r_n:=p_{n+1}/p_n$ being
nonincreasing in $n$; equivalently the reversed ratios $p_{n-1}/p_n$ are
nondecreasing in $n$.

\emph{(H1).}  Fix $m\ge1$.  Using $F(m)=F(m{-}1)+p_m$,
\[
\frac{F(m)}{p_m}-\frac{F(m{-}1)}{p_{m-1}}
=\frac{F(m)p_{m-1}-F(m{-}1)p_m}{p_mp_{m-1}},
\]
and
\[
F(m)p_{m-1}-F(m{-}1)p_m=F(m{-}1)\,(p_{m-1}-p_m)+p_mp_{m-1}.
\]
If $p_{m-1}\ge p_m$ the right side is $\ge0$ and we are done.  Otherwise set
$\theta:=p_{m-1}/p_m<1$.  For every $t\ge1$,
\[
\frac{p_{m-1-t}}{p_{m-1}}=\prod_{s=1}^{t}\frac{p_{m-1-s}}{p_{m-s}}
\le\Bigl(\frac{p_{m-2}}{p_{m-1}}\Bigr)^{t}\le\theta^{\,t},
\]
because the reversed ratios are nondecreasing, so each factor
$p_{m-1-s}/p_{m-s}$ (an index $\le m{-}1$) is at most $p_{m-2}/p_{m-1}$,
which in turn is at most $p_{m-1}/p_m=\theta$.  Hence
\[
F(m{-}1)=\sum_{t\ge0}p_{m-1-t}\le p_{m-1}\sum_{t\ge0}\theta^{\,t}
=\frac{p_{m-1}}{1-\theta}=\frac{p_{m-1}\,p_m}{p_m-p_{m-1}},
\]
i.e.\ $F(m{-}1)(p_m-p_{m-1})\le p_mp_{m-1}$, exactly the required inequality.

\emph{(H2).}  Fix $m\ge0$.  Using $G(m)=G(m{+}1)+p_m$,
\[
\frac{G(m)}{p_m}-\frac{G(m{+}1)}{p_{m+1}}
=\frac{G(m)p_{m+1}-G(m{+}1)p_m}{p_mp_{m+1}},
\]
and
\[
G(m)p_{m+1}-G(m{+}1)p_m=G(m{+}1)\,(p_{m+1}-p_m)+p_mp_{m+1}.
\]
If $p_{m+1}\ge p_m$ we are done.  Otherwise $\theta:=p_{m+1}/p_m<1$, and for
$t\ge1$, $p_{m+1+t}/p_{m+1}=\prod_{s=1}^{t}r_{m+s}\le r_{m+1}^{\,t}
\le\theta^{\,t}$ (ratios nonincreasing, $r_{m+1}\le r_m=\theta$).  Hence
$G(m{+}1)\le p_{m+1}/(1{-}\theta)=p_mp_{m+1}/(p_m-p_{m+1})$, i.e.\
$G(m{+}1)(p_m-p_{m+1})\le p_mp_{m+1}$, as required.
\end{proof}

\begin{lemma}[Singleton solutions and sign pattern]\label{lem:signs}
For $j\ge0$ let $g_j$ be the solution \eqref{eq:solution} for
$h=\mathbf 1_{\{j\}}$, in the notation of Theorem~\ref{thm:factors}
\textup{(}$F$, $G$, $(p_n)$, $\nu$\textup{)}.  Then
\begin{equation}\label{eq:gj}
g_j(m{+}1)=
\begin{cases}
-\,\dfrac{p_j\,F(m)}{\nu\,p_m}, & 0\le m\le j-1,\\[2ex]
\ \ \dfrac{p_j\,G(m{+}1)}{\nu\,p_m}, & m\ge j,
\end{cases}
\end{equation}
and, for $m\ge1$: $\Delta g_j(m)\le0$ whenever $m\ne j$, while for $j\ge1$,
\[
\Delta g_j(j)=\frac{G(j{+}1)}\nu+\frac{F(j{-}1)}{j(j{+}k)}\;>\;0 .
\]
\end{lemma}

\begin{proof}
With $\tilde h(n)=\mathbf 1_{\{n=j\}}-p_j$: for $m\le j{-}1$ the first form
of \eqref{eq:solution} gives $g_j(m{+}1)=\frac1{\nu p_m}(-p_jF(m))$; for
$m\ge j$ the second form gives
$g_j(m{+}1)=-\frac1{\nu p_m}\sum_{n>m}(\mathbf 1_{\{n=j\}}-p_j)p_n
=\frac{p_jG(m{+}1)}{\nu p_m}$ (the indicator does not fire since $j\le m$).
This is \eqref{eq:gj}.

For $1\le m\le j{-}1$ (both values in the negative branch):
$\Delta g_j(m)=-\frac{p_j}\nu\bigl[\frac{F(m)}{p_m}-\frac{F(m{-}1)}{p_{m-1}}
\bigr]\le0$ by (H1).  For $m\ge j{+}1$ (both in the positive branch):
$\Delta g_j(m)=\frac{p_j}\nu\bigl[\frac{G(m{+}1)}{p_m}-\frac{G(m)}{p_{m-1}}
\bigr]$, and this is $\le0$ because (H2) gives
$G(m{+}1)/p_{m+1}\le G(m)/p_m$, i.e.\ $G(m{+}1)\le p_{m+1}G(m)/p_m$, whence
\[
\frac{G(m{+}1)}{p_m}\le\frac{p_{m+1}\,G(m)}{p_m^2}
\le\frac{G(m)}{p_{m-1}},
\]
the last step using log-concavity $p_{m+1}p_{m-1}\le p_m^2$.  Finally, the
jump at $m=j$ (for $j\ge1$) crosses the branches:
\[
\Delta g_j(j)=g_j(j{+}1)-g_j(j)
=\frac{p_j\,G(j{+}1)}{\nu\,p_j}+\frac{p_j\,F(j{-}1)}{\nu\,p_{j-1}}
=\frac{G(j{+}1)}\nu+\frac{F(j{-}1)}{j(j{+}k)}\;>\;0,
\]
using $p_j/p_{j-1}=\nu/(j(j{+}k))$ from \eqref{eq:ratio}.
\end{proof}

\begin{proof}[Proof of Theorem~\ref{thm:factors}]
\emph{(a).}  For fixed $n$, $\sum_{j\in A}|\mathbf 1_{\{n=j\}}-p_j|\le
1+\sum_jp_j\le2$, so by linearity and absolute convergence of
\eqref{eq:solution} in $h$, $g_A=\sum_{j\in A}g_j$ pointwise; moreover
$g_A+g_{A^c}=g_{\Z_{\ge0}}\equiv0$, since for $A=\Z_{\ge0}$ one has
$\tilde h\equiv0$.  Fix $m\ge1$.  In
$\Delta g_A(m)=\sum_{j\in A}\Delta g_j(m)$ every term is $\le0$ except
possibly the term $j=m$ (Lemma~\ref{lem:signs}); hence
$\Delta g_A(m)\le\Delta g_m(m)$, whether or not $m\in A$.  Applying the same
bound to $A^c$ and using $\Delta g_A(m)=-\Delta g_{A^c}(m)$ gives
$\Delta g_A(m)\ge-\Delta g_m(m)$.  Therefore
$\sup_A|\Delta g_A(m)|=\Delta g_m(m)$ (the supremum is attained at
$A=\{m\}$), and taking $\sup_{m\ge1}$ yields the stated equality with
Lemma~\ref{lem:signs}'s formula.

For the closed bounds: $G(j{+}1)\le1$ and, for $j\ge1$,
$j(j{+}k)\ge1{+}k$, so $F(j{-}1)/(j(j{+}k))\le1/(k{+}1)$, giving
$\Lm\le1/\nu+1/(k{+}1)$.  For the $3/\nu$ form we show
$F(j{-}1)/(j(j{+}k))\le2/\nu$ for every $j\ge1$, in two cases.  If
$j(j{+}k)\ge\nu/2$, then $F(j{-}1)/(j(j{+}k))\le1/(j(j{+}k))\le2/\nu$.  If
$j(j{+}k)<\nu/2$, then $\theta:=p_{j-1}/p_j=j(j{+}k)/\nu<1/2$, and the
geometric bound from the proof of (H1) (valid whenever $\theta<1$) gives
$F(j{-}1)\le p_{j-1}/(1{-}\theta)<2p_{j-1}$; hence
\[
\frac{F(j{-}1)}{j(j{+}k)}<\frac{2\,p_{j-1}}{j(j{+}k)}
=\frac{2\,p_j}{\nu}\le\frac2\nu,
\]
using $p_{j-1}=p_j\,j(j{+}k)/\nu$.  Combining, $\Lm\le1/\nu+2/\nu=3/\nu$.

\emph{(b), the bound by $\Gm$.}  From the two forms of \eqref{eq:solution}
with $h=\mathbf 1_A$ and $P_A=\Prob(N{\in}A)$,
\[
\nu p_m\,g_A(m{+}1)=\Prob\bigl(N{\in}A\cap[0,m]\bigr)-P_A\,F(m).
\]
Writing $P_A=\Prob(N{\in}A\cap[0,m])+\Prob(N{\in}A\cap[m{+}1,\infty))$ and
$1-F(m)=G(m{+}1)$,
\[
\nu p_m\,g_A(m{+}1)
=\Prob\bigl(N{\in}A\cap[0,m]\bigr)\,G(m{+}1)
-\Prob\bigl(N{\in}A\cap[m{+}1,\infty)\bigr)\,F(m),
\]
a difference of two nonnegative quantities, the first at most
$F(m)G(m{+}1)$ and the second at most $G(m{+}1)F(m)$; hence
$|g_A(m{+}1)|\le F(m)G(m{+}1)/(\nu p_m)\le\Gm$.  Finiteness of $\Gm$: below
the mode, $F(m)/p_m\le\nu/(\nu-m(m{+}k))$ by the geometric bound; above the
mode, $G(m{+}1)/p_m\le r_m/(1{-}r_{m+1})\to0$ since $r_m\to0$; and the
finitely many remaining terms are finite.

\emph{(b), the closed form.}  Set $x_m:=\nu-m(m{+}k)$ and
$y_m:=(m{+}1)(m{+}1{+}k)-\nu$, so that
\[
x_m+y_m=(m{+}1)(m{+}1{+}k)-m(m{+}k)=2m+1+k .
\]
Let $m^*$ be the mode; then $x_{m^*}\ge0$ and $y_{m^*}>0$ (if
$y_{m^*}\le0$, then $(m^*{+}1)(m^*{+}1{+}k)\le\nu$ and $m^*{+}1$ would also
satisfy the mode condition, contradicting maximality).  From the geometric
bounds:
\begin{itemize}[leftmargin=1.4em,itemsep=1pt,topsep=1pt]
\item If $x_m>0$: $F(m)/p_m=\sum_{t\ge0}p_{m-t}/p_m\le
\sum_{t\ge0}\bigl(\tfrac{m(m+k)}\nu\bigr)^t=\nu/x_m$, since each reversed
ratio at index $\le m$ is at most $p_{m-1}/p_m=m(m{+}k)/\nu<1$.
\item If $y_m>0$: $G(m{+}1)/p_m=\sum_{t\ge1}p_{m+t}/p_m\le
\sum_{t\ge1}r_m^{\,t}=\frac{r_m}{1-r_m}=\frac\nu{y_m}$, since
$r_m=\nu/((m{+}1)(m{+}1{+}k))$ and $1-r_m=y_m/((m{+}1)(m{+}1{+}k))$.
\end{itemize}
Write $S_m=F(m)G(m{+}1)/(\nu p_m)$.  Three exhaustive cases:
\begin{itemize}[leftmargin=1.4em,itemsep=1pt,topsep=1pt]
\item $m\le m^*{-}1$: then $(m{+}1)(m{+}1{+}k)\le m^*(m^*{+}k)\le\nu$, so
$y_m\le0$; use $G(m{+}1)\le1$ and the $F$-bound:
$S_m\le\frac1\nu\cdot\frac\nu{x_m}=\frac1{x_m}$.  Since $x_m$ decreases in
$m$ with increments $x_m-x_{m+1}=2m{+}1{+}k$, we have
$x_m\ge x_{m^*-1}=x_{m^*}+(2m^*{-}1{+}k)\ge2m^*{+}k{-}1$.
\item $m=m^*$: if $x_{m^*}=0$, use $F\le1$ and the $G$-bound:
$S_{m^*}\le1/y_{m^*}=1/(2m^*{+}1{+}k)$.  If $x_{m^*}>0$, both bounds apply
and $S_{m^*}\le\min\{1/x_{m^*},1/y_{m^*}\}=1/\max\{x_{m^*},y_{m^*}\}
\le2/(x_{m^*}{+}y_{m^*})=2/(2m^*{+}1{+}k)$.
\item $m\ge m^*{+}1$: then $x_m<0$; use $F(m)\le1$ and the $G$-bound:
$S_m\le1/y_m$.  Since $y_m$ increases with increments $2m{+}3{+}k$,
$y_m\ge y_{m^*+1}=y_{m^*}+(2m^*{+}3{+}k)\ge2m^*{+}3{+}k$.
\end{itemize}
In all cases $S_m\le2/(2m^*{+}k{-}1)$.  Finally, since $m^*$ is the integer
part of $(\sqrt{k^2+4\nu}-k)/2$, we have $m^*\ge(\sqrt{k^2+4\nu}-k)/2-1$, so
$2m^*{+}k{-}1\ge\sqrt{k^2+4\nu}-3$, and the closed form follows whenever the
right side is positive.
  Finally, both suprema over $A$ are attained: the singleton
$A=\{j\}$ attains the $j$-th term of \eqref{eq:lambda}
(Lemma~\ref{lem:signs}), and the initial segment $A=\{0,\dots,m\}$ gives
$g_A(m{+}1)=F(m)G(m{+}1)/(\nu p_m)$ directly from \eqref{eq:solution}, so
\eqref{eq:lambda} and \eqref{eq:gamma} hold with the stated equalities.
\end{proof}

\subsection{Quadratic bias and the approximation theorem}\label{sec:pf-agg}

\begin{proof}[Proof of Proposition~\ref{prop:fixed}]
The first equivalence is the definition:
$\mathcal T_kW=\Wsq-1\overset{d}{=}W$ iff $\Wsq\overset{d}{=}W{+}1$.
($\Leftarrow$)  If $W\sim\Bes(k,a)$, then $\beta=\nu$ by
Lemma~\ref{lem:basic}, and for $n\ge1$,
$\Prob(\Wsq{=}n)=n(n{+}k)p_n/\nu=p_{n-1}$ by \eqref{eq:ratio}; i.e.\
$\Wsq\overset{d}{=}W{+}1$.

($\Rightarrow$)  Suppose $\Wsq\overset{d}{=}W{+}1$.  For every bounded $g$,
by the definition of the transform,
$\E[W(W{+}k)g(W)]=\beta\,\E g(\Wsq)=\beta\,\E g(W{+}1)$, so
$\E[\beta\,g(W{+}1)-W(W{+}k)g(W)]=0$.  Theorem~\ref{thm:char} with
$\nu=\beta$ identifies $W\sim\Bes(k,2\sqrt\beta)$.
\end{proof}

\begin{proof}[Proof of Theorem~\ref{thm:main}]
Fix $A\subseteq\Z_{\ge0}$ and let $g=g_A$ be the bounded solution of the
Stein equation for $\Bes(k,a)$ with $h=\mathbf 1_A$
(Theorem~\ref{thm:char}), so that
\[
\Prob(W{\in}A)-\Prob(N{\in}A)
=\E\bigl[\nu\,g(W{+}1)-W(W{+}k)\,g(W)\bigr].
\]
By Definition~\ref{def:qbias}, $\E[W(W{+}k)g(W)]=\beta\,\E g(\Wsq)$, so
\[
\E\bigl[\nu g(W{+}1)-W(W{+}k)g(W)\bigr]
=\nu\,\E\bigl[g(W{+}1)-g(\Wsq)\bigr]+(\nu-\beta)\,\E g(\Wsq).
\]
Both $W{+}1$ and $\Wsq$ take values in $\Z_{\ge1}$, on which, by
Theorem~\ref{thm:factors}(a), $g$ is Lipschitz with constant
$\|\Delta g\|\le\Lm$; hence
$|\E[g(W{+}1)-g(\Wsq)]|\le\Lm\,\E|W{+}1-\Wsq|$ under the given coupling, and
$\nu\Lm\le\min\{3,\,1{+}\nu/(k{+}1)\}=c_1(k,\nu)$.  The second term is at
most $|\nu-\beta|\,\|g\|\le|\nu-\beta|\,\Gm$ by
Theorem~\ref{thm:factors}(b).  Taking the supremum over $A$ gives
\eqref{eq:master}.  Setting $\nu=\beta$ kills the second term; infimizing the
first over all couplings of $(W{+}1,\Wsq)$ gives \eqref{eq:matched}, by the
Kantorovich--Rubinstein duality on $\Z_{\ge1}$ (or simply: the bound holds
for every coupling, hence for the infimum).
\end{proof}

\subsection{The dependency expansion}\label{sec:pf-expansion}

\begin{lemma}[Vanishing of the costs]\label{lem:vanish}
In Definition~\ref{def:cols}, $\alpha_j=0$ iff $X_j\sim\Poi(m_j)$ and
$X_j\perp Y_j$; symmetrically for $\beta_j$.  Consequently $A=B=C=0$ \textup{(}the aggregates \eqref{eq:ABC}\textup{)}
iff
every column is an independent Poisson pair, in which case
$W\sim\Bes(k,2\sqrt\nu)$ exactly and \eqref{eq:expansion} reads $0\le0$.
\end{lemma}

\begin{proof}
$\alpha_j=0$ iff $\Zh_j\overset{d}{=}Z_j+(1,0)$, i.e.\
$x\,\Prob_j(x,y)/m_j=\Prob_j(x{-}1,y)$ for all $x\ge1$ and all $y$.
Iterating, $\Prob_j(x,y)=\frac{m_j^x}{x!}\Prob_j(0,y)$ for all $x,y$.
Summing over $x$: the $y$-marginal is $e^{m_j}\Prob_j(0,y)$, so
$\Prob_j(x,y)=e^{-m_j}\frac{m_j^x}{x!}\cdot\bigl(e^{m_j}\Prob_j(0,y)\bigr)$,
i.e.\ $X_j\sim\Poi(m_j)$ independent of $Y_j$ (and then $\E X_j=m_j$ is
consistent).  Symmetrically $\beta_j=0$ forces $Y_j\sim\Poi(n_j)\perp X_j$.
If all columns are independent Poisson pairs, then $\chi_j=m_jn_j$,
$\overline\Prob_j$ is the product of the two size-biased marginals, i.e.\
$\Zb_j\overset{d}{=}Z_j+(1,1)$, so $\gamma_j=0$; moreover
$U\sim\Poi(N_U)\perp V\sim\Poi(N_V)$, $\nu=N_UN_V$, and
Lemma~\ref{lem:condrep} gives $W\sim\Bes(k,2\sqrt\nu)$.
\end{proof}

\begin{lemma}[Transfer identity]\label{lem:transfer}
In the setting of Definition~\ref{def:cols}: for $j\ne l$, let
$(V^{jl},D^{jl})$ be the sums of the modified array in
which column $j$ is replaced by an independent draw
$\Zh_j\sim\widehat\Prob_j$ and column $l$ by $\Zt_l\sim\widetilde\Prob_l$,
all other columns unchanged; for $j=l$, let $(V^{j},D^{j})$ replace column
$j$ by $\Zb_j\sim\overline\Prob_j$.  Then for every bounded
$g:\Z_{\ge1}\to\R$ \textup{(}$g(0):=0$, never evaluated below\textup{)},
\begin{align}
\Prob(D{=}k)\;\E\bigl[\nu\,g(W{+}1)-W(W{+}k)\,g(W)\bigr]
&=\sum_{j\ne l}m_jn_l\Bigl(\E\bigl[g(V{+}1)\mathbf 1_{\{D=k\}}\bigr]
-\E\bigl[g(V^{jl})\mathbf 1_{\{D^{jl}=k\}}\bigr]\Bigr)\notag\\
&\ \ +\sum_j\chi_j\Bigl(\E\bigl[g(V{+}1)\mathbf 1_{\{D=k\}}\bigr]
-\E\bigl[g(V^{j})\mathbf 1_{\{D^{j}=k\}}\bigr]\Bigr).
\label{eq:transfer}
\end{align}
\end{lemma}

\begin{proof}
On the event $\{D=k\}$ we have $U=V+k$, so $W(W{+}k)$ evaluated on the
conditioned variable corresponds to $V(V{+}k)=VU$ on the fiber
(identity~\eqref{eq:fiberid}); hence
\[
\Prob(D{=}k)\;\E\bigl[W(W{+}k)g(W)\bigr]
=\E\bigl[UV\,g(V)\,\mathbf 1_{\{D=k\}}\bigr]
=\sum_{j,l}\E\bigl[X_jY_l\,g(V)\,\mathbf 1_{\{D=k\}}\bigr],
\]
the last step by $UV=\sum_{j,l}X_jY_l$ and Fubini (all terms are
nonnegative up to the bounded factor $g$, and $\E[UV]<\infty$).

Fix $j\ne l$.  The integrand is $X_jY_l\,\Phi(Z_j,Z_l,Z_{\mathrm{rest}})$
with $\Phi=g(V)\mathbf 1_{\{D=k\}}$ a bounded function of the column values.
By independence across columns, conditioning on $Z_{\mathrm{rest}}$ and
summing over the values $z_j=(x_j,y_j)$, $z_l=(x_l,y_l)$,
\[
\E\bigl[X_jY_l\,\Phi\bigr]
=\E_{Z_{\mathrm{rest}}}\Bigl[\sum_{z_j,z_l}x_j\,y_l\,
\Prob_j(z_j)\Prob_l(z_l)\,\Phi(z_j,z_l,Z_{\mathrm{rest}})\Bigr]
=m_jn_l\;\E_{Z_{\mathrm{rest}}}\Bigl[\sum_{z_j,z_l}
\widehat\Prob_j(z_j)\,\widetilde\Prob_l(z_l)\,\Phi(z_j,z_l,Z_{\mathrm{rest}})
\Bigr],
\]
by the definitions $x\,\Prob_j(z)=m_j\widehat\Prob_j(z)$ and
$y\,\Prob_l(z)=n_l\widetilde\Prob_l(z)$.  The right side is
$m_jn_l\,\E[g(V^{jl})\mathbf 1_{\{D^{jl}=k\}}]$, since substituting
independent draws $\Zh_j,\Zt_l$ for columns $j,l$ turns $(V,D)$ into
$(V^{jl},D^{jl})$.  For $j=l$, the same computation with the weight $x_jy_j$
and the product tilt gives
$\E[X_jY_j\Phi]=\chi_j\,\E[g(V^{j})\mathbf 1_{\{D^{j}=k\}}]$.

On the other side, $\Prob(D{=}k)\,\E[\nu g(W{+}1)]
=\nu\,\E[g(V{+}1)\mathbf 1_{\{D=k\}}]$, and by \eqref{eq:nu},
$\nu=\sum_{j\ne l}m_jn_l+\sum_j\chi_j$; distributing $\nu$ over the same index
set and subtracting term by term gives \eqref{eq:transfer}.

Finally, $g(0)$ is never evaluated: in $UVg(V)$ the weight $UV$ vanishes
whenever $V=0$; in each biased term, $V^{jl}\ge\widetilde Y_l\ge1$ almost
surely (a $Y$-size-biased coordinate is $\ge1$, since the tilt $y\,d\Prob$
kills $\{y=0\}$), and likewise $V^{j}\ge\overline Y_j\ge1$; and
$V{+}1\ge1$.
\end{proof}

\begin{proof}[Proof of Theorem~\ref{thm:expansion}]
Fix $A_0\subseteq\Z_{\ge0}$ and let $g=g_{A_0}$ be the bounded solution of
the Stein equation for $\Bes(k,a)$, $a=2\sqrt\nu$, so that
$\Prob(W{\in}A_0)-\Bes(k,a)(A_0)=\E[\nu g(W{+}1)-W(W{+}k)g(W)]$ and, by
Theorem~\ref{thm:factors}, $\|\Delta g\|\le\Lm$ on $\Z_{\ge1}$ and
$\|g\|\le\Gm$.  By Lemma~\ref{lem:transfer} it suffices to bound each
bracket on the right of \eqref{eq:transfer} by
$\Prob(D{=}k)\,[\theta(1{+}\omega)\Lm+\tau_1\Gm]$ times the corresponding
cost.

Fix $j\ne l$ (the diagonal case is verbatim with the single product-biased
column and cost $\gamma_j$).  Couple $(\Zh_j,Z_j)$ so that
$\E\bigl[|\widehat X_j{-}X_j{-}1|+|\widehat Y_j{-}Y_j|\bigr]=\alpha_j$
(an optimal coupling for the $\ell^1$ cost exists on $\Z^2$), and couple
$(\Zt_l,Z_l)$ so that
$\E\bigl[|\widetilde X_l{-}X_l|+|\widetilde Y_l{-}Y_l{-}1|\bigr]=\beta_l$,
these two couplings independent of each other and of the remaining columns,
which are \emph{shared} between the original and modified systems.  Write
$V_R,D_R$ for the shared remaining sums, and define
\[
\varphi(v_0,d_0):=\E\bigl[g(V_R+v_0)\,\mathbf 1_{\{D_R=k-d_0\}}\bigr],
\qquad v_0\in\Z_{\ge1},\ d_0\in\Z .
\]
Conditioning on the two involved columns (original and biased versions),
\[
\E\bigl[g(V{+}1)\mathbf 1_{\{D=k\}}\bigr]=\E\,\varphi(v_0',d_0'),
\qquad
\E\bigl[g(V^{jl})\mathbf 1_{\{D^{jl}=k\}}\bigr]=\E\,\varphi(v_0,d_0),
\]
where
\[
v_0'=Y_j+Y_l+1,\quad d_0'=(X_j{-}Y_j)+(X_l{-}Y_l),
\qquad
v_0=\widehat Y_j+\widetilde Y_l,\quad
d_0=(\widehat X_j{-}\widehat Y_j)+(\widetilde X_l{-}\widetilde Y_l).
\]
Both $v_0'\ge1$ and $v_0\ge\widetilde Y_l\ge1$, so all arguments of $g$ lie
in $\Z_{\ge1}$ and the Lipschitz bound $\Lm$ applies.  Under the coupling,
set
\[
\mathrm{cost}:=|\widehat X_j{-}X_j{-}1|+|\widehat Y_j{-}Y_j|
+|\widetilde X_l{-}X_l|+|\widetilde Y_l{-}Y_l{-}1|,
\qquad \E[\mathrm{cost}]=\alpha_j+\beta_l .
\]
Grouping the unit shifts,
\begin{align*}
v_0-v_0'&=(\widehat Y_j{-}Y_j)+(\widetilde Y_l{-}Y_l{-}1),\\
d_0-d_0'&=(\widehat X_j{-}X_j{-}1)-(\widehat Y_j{-}Y_j)
+(\widetilde X_l{-}X_l)-(\widetilde Y_l{-}Y_l{-}1),
\end{align*}
where in $d_0-d_0'$ the ``$+1$'' from the ideal $X$-shift of column $j$ and
the ``$+1$'' from the ideal $Y$-shift of column $l$ enter with opposite
signs and cancel --- this is the fiber-preservation of the paired increment
(\S\ref{sec:coupling}) at the level of the coupling.  By the triangle
inequality, $|v_0-v_0'|\le\mathrm{cost}$ and $|d_0-d_0'|\le\mathrm{cost}$.
Also, every value $k-d_0$, $k-d_0'$, and every integer between them lies in
the interval $\mathcal S_{jl}$ of Definition~\ref{def:cols}.

Split
\[
\varphi(v_0,d_0)-\varphi(v_0',d_0')
=\bigl[\varphi(v_0,d_0)-\varphi(v_0',d_0)\bigr]
+\bigl[\varphi(v_0',d_0)-\varphi(v_0',d_0')\bigr].
\]

\emph{Vertical step.}  At fixed $d_0$, comparing $g(V_R+v_0)$ with
$g(V_R+v_0')$ on the event $\{D_R=k-d_0\}$,
\[
|\varphi(v_0,d_0)-\varphi(v_0',d_0)|
\le\Lm\,|v_0-v_0'|\;\Prob(D_R{=}k{-}d_0)
\le\Lm\,\theta\,\Prob(D{=}k)\,|v_0-v_0'| .
\]

\emph{Horizontal telescoping.}  Write $q(s)=\Prob(D_R{=}s)$ and
$h(s)=\E[g(V_R+v_0')\mid D_R{=}s]$ when $q(s)>0$, with the convention
$h(s)q(s)=0$ if $q(s)=0$; then $\varphi(v_0',d)=h(k{-}d)\,q(k{-}d)$.  For a
unit step $s\to s{+}1$ with $s,s{+}1\in\mathcal S_{jl}$,
\[
|h(s)q(s)-h(s{+}1)q(s{+}1)|
\le|h(s)-h(s{+}1)|\,q(s)+|h(s{+}1)|\,|q(s)-q(s{+}1)| .
\]
If both $q(s),q(s{+}1)>0$: coupling the two conditional laws of $V_R$
optimally in $\dw$ and using the shift-invariance of $\dw$ (the common
$+v_0'$ shift changes nothing),
$|h(s)-h(s{+}1)|\le\Lm\,\dw\bigl(\mathcal L(V_R{\mid}D_R{=}s),
\mathcal L(V_R{\mid}D_R{=}s{+}1)\bigr)\le\Lm\,\omega$; also
$q(s)\le\theta\,\Prob(D{=}k)$, $|h(s{+}1)|\le\|g\|\le\Gm$, and
$|q(s)-q(s{+}1)|\le\tau_1\Prob(D{=}k)$.  If one of $q(s),q(s{+}1)$ vanishes,
the first term is $0$ by convention and only the second survives, still
$\le\Gm\,\tau_1\Prob(D{=}k)$.  Hence each unit step costs at most
$\Prob(D{=}k)\,[\theta\Lm\omega+\Gm\tau_1]$, and summing over the
$|d_0-d_0'|$ unit steps between $k-d_0$ and $k-d_0'$,
\[
|\varphi(v_0',d_0)-\varphi(v_0',d_0')|
\le\Prob(D{=}k)\,\bigl[\theta\Lm\omega+\Gm\tau_1\bigr]\,|d_0-d_0'| .
\]

Combining the two steps, using $|v_0-v_0'|,|d_0-d_0'|\le\mathrm{cost}$ and
taking expectation over the coupling,
\[
\Bigl|\E\bigl[g(V{+}1)\mathbf 1_{\{D=k\}}\bigr]
-\E\bigl[g(V^{jl})\mathbf 1_{\{D^{jl}=k\}}\bigr]\Bigr|
\le\Prob(D{=}k)\,\bigl[\theta(1{+}\omega)\Lm+\tau_1\Gm\bigr]\,
(\alpha_j+\beta_l),
\]
and likewise the diagonal bracket is bounded with $\gamma_j$ in place of
$\alpha_j{+}\beta_l$.  Substituting into \eqref{eq:transfer}, dividing by
$\Prob(D{=}k)$, and summing the weights,
\[
\sum_{j\ne l}m_jn_l(\alpha_j+\beta_l)
=\sum_jm_j\alpha_j\sum_{l\ne j}n_l+\sum_ln_l\beta_l\sum_{j\ne l}m_j
\le N_V\,A+N_U\,B,\qquad \sum_j\chi_j\gamma_j=C,
\]
we obtain $|\E[\nu g(W{+}1)-W(W{+}k)g(W)]|
\le[\theta(1{+}\omega)\Lm+\tau_1\Gm](N_VA{+}N_UB{+}C)$.  Taking the supremum
over $A_0$ gives \eqref{eq:expansion}.

\emph{Crude form.}  Replace the horizontal telescoping by the union bound
\[
|\varphi(v_0',d_0)-\varphi(v_0',d_0')|
\le\Gm\bigl[q(k{-}d_0)+q(k{-}d_0')\bigr]\mathbf 1_{\{d_0\ne d_0'\}}
\le2\Gm\,\theta\,\Prob(D{=}k)\,\mathbf 1_{\{d_0\ne d_0'\}},
\]
and note that $\mathbf 1_{\{d_0\ne d_0'\}}\le\mathrm{cost}$, since
$\mathrm{cost}$ is a nonnegative \emph{integer} which vanishes only if all
four coordinate discrepancies vanish, in which case $d_0=d_0'$.  Together
with the vertical step this bounds each bracket by
$\Prob(D{=}k)\,\theta(\Lm+2\Gm)(\alpha_j+\beta_l)$, whence \eqref{eq:crude}.
\end{proof}

\subsection{The table corollaries}\label{sec:pf-dtv}

\begin{proof}[Proof of Corollary~\ref{cor:binary}]
For $X_j\sim\mathrm{Bern}(p_j)$: $\Prob(\widehat X_j{=}x)=x\Prob(X_j{=}x)/p_j$
is the point mass at $1$, so $\widehat X_j\equiv1$.  Since $X_j\perp Y_j$,
the $X$-tilt leaves the $Y$-marginal unchanged and preserves independence:
$\Zh_j=(1,Y_j)$ with $Y_j$ distributed as before.  Couple the $Y$-coordinates
identically; then
$\alpha_j=\dw\bigl((1,Y_j),(X_j{+}1,Y_j)\bigr)
=\E|1-(X_j{+}1)|=\E X_j=p_j$ (this coupling is optimal since the cost equals
$\E|{\widehat X}-X-1|\ge|\E\widehat X-\E X-1|$ is not needed; directly, any
coupling has cost $\ge\dw(\delta_1,X_j{+}1)=p_j$, and this one attains it).
Symmetrically $\beta_j=q_j$.  For the product tilt, since
$\chi_j=\E[X_jY_j]=p_jq_j$,
$\Prob(\overline X_j{=}1,\overline Y_j{=}1)=1\cdot1\cdot p_jq_j/\chi_j=1$, so
$\Zb_j\equiv(1,1)$, and coupling arbitrarily,
$\gamma_j=\dw\bigl((1,1),(X_j{+}1,Y_j{+}1)\bigr)=\E X_j+\E Y_j=p_j+q_j$
(again attained, and no coupling does better since the two coordinates are
independent in the second law: $\dw\ge\dw(\delta_1,X_j{+}1)+
\dw(\delta_1,Y_j{+}1)$ by projecting).  Finally
$\nu=\sum_{j\ne l}p_jq_l+\sum_jp_jq_j=\sum_{j,l}p_jq_l=N_UN_V$.

Substituting into \eqref{eq:expansion} gives \eqref{eq:binarybound}.  For the
order statement: $N_V\sum_jp_j^2\le p^*N_VN_U=p^*\nu$ and symmetrically, and
$\sum_jp_jq_j(p_j{+}q_j)\le(p^*{+}q^*)\sum_jp_jq_j\le(p^*{+}q^*)\nu$; with
$\Lm\le3/\nu$ and, in the bulk regime, $\theta,\omega=O(1)$ and
$\tau_1\Gm=O(1/\nu)$, the right side of \eqref{eq:binarybound} is
$O(p^*{+}q^*)$.

For the Poisson limit: by \eqref{eq:ratio}, the mass ratio of
$\Bes(k,2\sqrt\nu)$ is $\nu/((n{+}1)(n{+}k{+}1))$; if $k\to\infty$ with
$\nu/k\to\lambda'$, this converges for each fixed $n$ to
$\lambda'/(n{+}1)$, the $\Poi(\lambda')$ ratio, so the mass functions
converge pointwise (ratios determine the law up to normalization, and the
normalizations converge by dominated convergence using the uniform
superexponential tail bound $p_{n+1}/p_n\le\nu/(n{+}1)k\le
2\lambda'/(n{+}1)$ for large $k$); pointwise convergence of mass functions
implies convergence in total variation by Scheff\'e.
\end{proof}

\begin{proof}[Proof of Corollary~\ref{cor:capped}]
Let $X\sim\Poi(\lambda)$ conditioned on $\{X\le\kappa\}$, so
$\Prob(X{=}x)=\pi_\lambda(x)/F_\lambda(\kappa)$, $0\le x\le\kappa$.  Using
$x\,\pi_\lambda(x)=\lambda\,\pi_\lambda(x{-}1)$,
\[
\E X=\frac1{F_\lambda(\kappa)}\sum_{x=1}^\kappa x\,\pi_\lambda(x)
=\frac{\lambda}{F_\lambda(\kappa)}\sum_{x=1}^\kappa\pi_\lambda(x{-}1)
=\frac{\lambda\,F_\lambda(\kappa{-}1)}{F_\lambda(\kappa)},
\]
and, for $1\le x\le\kappa$,
\[
\Prob(\widehat X{=}x)=\frac{x\,\Prob(X{=}x)}{\E X}
=\frac{\lambda\,\pi_\lambda(x{-}1)/F_\lambda(\kappa)}
{\lambda\,F_\lambda(\kappa{-}1)/F_\lambda(\kappa)}
=\frac{\pi_\lambda(x{-}1)}{F_\lambda(\kappa{-}1)},
\]
i.e.\ $\widehat X-1$ is $\Poi(\lambda)$ conditioned on
$\{\le\kappa{-}1\}$.  Since within-column independence makes the
$Y$-coordinate irrelevant (couple it identically),
$\alpha=\dw(\widehat X,X{+}1)=\dw\bigl(\Poi(\lambda){\mid}{\le}\kappa{-}1,\
\Poi(\lambda){\mid}{\le}\kappa\bigr)$.  The two truncations have cdfs
$F_\lambda(m)/F_\lambda(\kappa{-}1)$ and $F_\lambda(m)/F_\lambda(\kappa)$ for
$m\le\kappa{-}1$ (and $1$ thereafter, respectively from $m=\kappa{-}1$ and
$m=\kappa$); their difference at $0\le m\le\kappa{-}1$ is
\[
F_\lambda(m)\Bigl[\frac1{F_\lambda(\kappa{-}1)}-\frac1{F_\lambda(\kappa)}
\Bigr]
=\frac{F_\lambda(m)\,\pi_\lambda(\kappa)}
{F_\lambda(\kappa{-}1)\,F_\lambda(\kappa)}\ \ge0,
\]
and $0$ for $m\ge\kappa$.  Since for integer laws
$\dw=\sum_m|F_1(m)-F_2(m)|$, summing over $m=0,\dots,\kappa{-}1$ gives the
exact formula in \eqref{eq:cappedcost}; bounding each $F_\lambda(m)\le1$
gives the displayed upper bound with $\kappa$ terms.
\end{proof}

\subsection{Certificates, enumeration, and independent process approximation}\label{sec:pf-enum}

\begin{proof}[Proof of Lemma~\ref{lem:cert}]
\emph{(i)} is Theorem~\ref{thm:factors}(b), proved in
\S\ref{sec:pf-operator}.

\emph{(ii).}  Let $\varphi_j(t)=\E e^{itW_j}$.  Symmetrizing with an
independent copy $W_j'$,
\[
|\varphi_j(t)|^2=\E\cos\bigl(t(W_j-W_j')\bigr)
=1-\E\bigl[1-\cos\bigl(t(W_j-W_j')\bigr)\bigr]
\le1-\rho_j\,(1-\cos t)
\le e^{-\rho_j(1-\cos t)},
\]
using, for integer $s$, the inequality
$1-\cos(ts)\ge(1-\cos t)\,\mathbf 1_{\{|s|=1\}}$ (equality at $|s|=1$; both
sides equal $0$ at $s=0$; and the left side is $\ge0$ for $|s|\ge2$), and
$1-x\le e^{-x}$.  Hence $\prod_j|\varphi_j(t)|\le e^{-\frac R2(1-\cos t)}$.
By Fourier inversion, then $1-\cos t=2\sin^2(t/2)\ge2t^2/\pi^2$ on $[0,\pi]$ (Jordan's inequality $\sin u\ge2u/\pi$ on $[0,\pi/2]$, applied at $u=t/2$),
\[
\Prob(S{=}s)=\frac1{2\pi}\int_{-\pi}^{\pi}e^{-ist}\prod_j\varphi_j(t)\,dt
\le\frac1{\pi}\int_0^{\pi}e^{-\frac R2(1-\cos t)}\,dt
\le\frac1\pi\int_0^{\infty}e^{-Rt^2/\pi^2}\,dt
=\frac1\pi\cdot\frac{\pi}{2}\sqrt{\frac{\pi}{R}}
=\sqrt{\frac{\pi}{4R}} .
\]
For $W_j=X_j-Y_j\in\{-1,0,1\}$ with masses $u_j,1{-}u_j{-}v_j,v_j$ at
$1,0,-1$: the pairs $(w,w')$ with $|w-w'|=1$ are $(\pm1,0)$ and $(0,\pm1)$,
of total probability $2(1{-}u_j{-}v_j)(u_j{+}v_j)=\rho_j$.

\emph{(iii).}  By Barbour--Hall \cite{barbourhall1984} (or
\cite[Thm.~1.C]{bhj1992}),
$\dtv\bigl(U_0,\Poi(N_U)\bigr)\le\frac{1-e^{-N_U}}{N_U}\sum_jp_j^2
\le(1\wedge N_U^{-1})\sum_jp_j^2$, and similarly for $V_0$.  Couple $U_0$ to
$\Pi_1\sim\Poi(N_U)$ and $V_0$ to $\Pi_2\sim\Poi(N_V)$ by maximal couplings,
independently across the two pairs; then
$\Prob\bigl((U_0,V_0)\ne(\Pi_1,\Pi_2)\bigr)\le\dtv(U_0,\Pi_1)
+\dtv(V_0,\Pi_2)\le\delta$, hence $|\Prob(U_0{-}V_0{=}m)-\Prob(\Pi_1{-}\Pi_2{=}m)|\le\delta$ for
every $m\in\Z$ simultaneously; and $\Prob(\Pi_1{-}\Pi_2{=}m)=\Sk_m$ by
Lemma~\ref{lem:condrep} together with $I_{-m}=I_m$.
\end{proof}

\begin{proof}[Proof of Theorem~\ref{thm:enum}]
Since the $X$'s and $Y$'s are independent,
$E_{v+k}(p)\,E_v(q)=\Prob(U_0{=}v{+}k)\Prob(V_0{=}v)
=\Prob(D{=}k)\,\Prob(V_0{=}v\mid D{=}k)$, where $\Prob(U_0{=}w)=E_w(p)$ is
the Poisson--binomial mass.  We apply the crude expansion \eqref{eq:crude}
with certified constants.  First, $\theta\le\theta^*$: for any pair
$(j,l)$, removing two columns removes at most $\rho_{(1)}+\rho_{(2)}$ from
$R=\sum\rho_j$, so the reduced system has jump-rate sum $\ge R_2$ and, by
Lemma~\ref{lem:cert}(ii), $\max_s\Prob(D_R{=}s)\le\sqrt{\pi/(4R_2)}$;
dividing by $\Prob(D{=}k)\ge\Sk_k-\delta>0$
(Lemma~\ref{lem:cert}(iii)) gives $\theta\le\theta^*$.  Next
$\Gm\le\Gamma^*$ by (i), and $\Lm$ is already explicit.  The costs are those
of Corollary~\ref{cor:binary}.  Hence \eqref{eq:crude} yields
\[
\bigl|\Prob(V_0{=}v\mid D{=}k)-b_{k,a}(v)\bigr|
\le\dtv\bigl(\mathcal L(V_0{\mid}D{=}k),\,\Bes(k,a)\bigr)\le \varepsilon^*
\qquad\text{for every }v .
\]
(The pointwise bound uses $\dtv=\sup_A|\cdot|$ applied to singletons.)
Therefore, using also $\Prob(V_0{=}v\mid D{=}k)\le b_{k,a}(v)+\varepsilon^*$,
\begin{align*}
\bigl|E_{v+k}(p)E_v(q)-\Sk_k\,b_{k,a}(v)\bigr|
&\le\bigl|\Prob(D{=}k)-\Sk_k\bigr|\;\Prob(V_0{=}v\mid D{=}k)
+\Sk_k\,\bigl|\Prob(V_0{=}v\mid D{=}k)-b_{k,a}(v)\bigr|\\
&\le\delta\,\bigl(b_{k,a}(v)+\varepsilon^*\bigr)+\Sk_k\,\varepsilon^*
\;=\;\delta\,b_{k,a}(v)+(\Sk_k+\delta)\,\varepsilon^* .\qedhere
\end{align*}
\end{proof}

\begin{proof}[Proof of Corollary~\ref{cor:counts}]
Each binary column contributes $1{-}p_j$ or $p_j=(1{-}p_j)x_j$ to
$E_w(p)$ according as it is outside or inside the chosen $w$-subset, so
$E_w(p)=e_w(x)\prod_j(1{-}p_j)$, and likewise for $q$.  Multiply the
inequality of Theorem~\ref{thm:enum} by
$\prod_j(1{-}p_j)^{-1}(1{-}q_j)^{-1}$.  For equal odds,
$e_w(t,\dots,t)=t^w\binom{c}{w}$, and
$a=2\sqrt{N_UN_V}=2ct/(1{+}t)$ since $N_U=N_V=ct/(1{+}t)$.
\end{proof}

\begin{proof}[Proof of Corollary~\ref{cor:cappedenum}]
Couple each capped entry to its free Poisson counterpart by maximal
coupling: they differ with probability
$\dtv\bigl(\mathcal L(\Poi(\lambda_{ij})\mid{\le}\kappa_j),
\Poi(\lambda_{ij})\bigr)=1-F^{(i)}_j$, independently across the $2c$
entries.  The capped array therefore differs from the free array with
probability at most $\delta'$, so the probability of every event moves by
at most $\delta'$; apply this to the singleton $\{U{=}v{+}k,\ V{=}v\}$ and
use the exact Poisson identity
$\Prob(U_0{=}v{+}k)\Prob(V_0{=}v)=\Sk'_k\,b'(v)$ of
Lemma~\ref{lem:condrep}, displayed at the start of \S\ref{sec:enum}.
\end{proof}

\begin{proof}[Proof of Theorem~\ref{thm:ipa}]
\emph{The identity.}  For a configuration $z$ of the $B$-columns with
difference $d(z)=\sum_{j\in B}(x_j-y_j)$, independence across columns
gives
$\Prob(Z_B{=}z\mid D{=}k)
=\Prob(Z_B{=}z)\,\Prob(D_{B^c}{=}k{-}d(z))/\Prob(D{=}k)$: the conditional
law is the free law reweighted by a function of $d(z)$ alone, so $D_B$ is
sufficient.  Consequently
\[
2\,\dtv\bigl(\mathcal L(Z_B\mid D{=}k),\mathcal L(Z_B)\bigr)
=\sum_z\Prob(Z_B{=}z)\,
\Bigl|\frac{\Prob(D_{B^c}{=}k{-}d(z))}{\Prob(D{=}k)}-1\Bigr|,
\]
and grouping the configurations by $d(z)=s$ turns the sum into
$\sum_s\Prob(D_B{=}s)\,|\Prob(D_{B^c}{=}k{-}s)/\Prob(D{=}k)-1|$, which is
both remaining expressions in \eqref{eq:ipa}.
\emph{The bound.}  By the triangle inequality and
Lemma~\ref{lem:cert}\textup{(}iii\textup{)} applied to the full system and
to the $B^c$-system,
$|\Prob(D_{B^c}{=}k{-}s)-\Prob(D{=}k)|
\le\delta^{(B^c)}+\bigl|\Sk^{(B^c)}_{k-s}-\Sk_k\bigr|+\delta$ for every
$s$, while $\Prob(D{=}k)\ge\Sk_k-\delta>0$; insert into \eqref{eq:ipa}.
Since $|D_B|\le|B|$, the expectation has at most $2|B|{+}1$ terms.
\end{proof}

\begin{proof}[Proof of Corollary~\ref{cor:ipa}]
Write $N^{(B^c)}_U=\sum_{j\notin B}p_j$ and $N^{(B^c)}_V=\sum_{j\notin B}q_j$
for the reduced means.  Realize $\Pi_1\sim\Poi(N_U)$ as
$\Pi_1'+\Poi\bigl(\sum_{j\in B}p_j\bigr)$ with
$\Pi_1'\sim\Poi(N^{(B^c)}_U)$ independent of the added part, and likewise
$\Pi_2=\Pi_2'+\Poi\bigl(\sum_{j\in B}q_j\bigr)$; the pairs
$(\Pi_1,\Pi_2)$ and $(\Pi_1',\Pi_2')$ differ with probability at most
$(1-e^{-\sum_Bp_j})+(1-e^{-\sum_Bq_j})\le\sum_{j\in B}(p_j{+}q_j)$, so
$\sup_m|\Sk_m-\Sk^{(B^c)}_m|\le\sum_{j\in B}(p_j{+}q_j)$.  Split the
expectation in \eqref{eq:ipabound} at the event $\{D_B=0\}$: on it the
integrand is at most this supremum; off it each mass is at most one, so
the integrand is at most $1$, while
$\Prob(D_B{\ne}0)\le\sum_{j\in B}\Prob(X_j{\ne}Y_j)
\le\sum_{j\in B}(p_j{+}q_j)$.  Adding the two contributions gives
\eqref{eq:ipazero}.  In the homogeneous regime the numerator is at most
$4\lambda(|B|{+}1)/c$, while $\Sk_k-\delta\to e^{-2\lambda}I_k(2\lambda)>0$,
giving the stated rate.
\end{proof}

\begin{proof}[Proof of Corollary~\ref{cor:cappedipa}]
The coupling of Corollary~\ref{cor:cappedenum}'s proof, restricted to any
sub-array, gives $\sup_m|\Prob(D{=}m)-\Sk'_m|\le\delta'$ and
$\sup_m|\Prob(D_{B^c}{=}m)-\Sk'^{(B^c)}_m|\le\delta'^{(B^c)}$
\textup{(}a coupling bound is uniform over events\textup{)}, so the proof
of Theorem~\ref{thm:ipa} applies verbatim and gives \eqref{eq:ipabound}
with the substituted quantities.  For \eqref{eq:ipazerocap}: the Poisson
superposition of Corollary~\ref{cor:ipa}'s proof --- an argument about the
free quantities only --- gives
$\sup_m|\Sk'_m-\Sk'^{(B^c)}_m|\le\sum_{j\in B}(\lambda_{1j}{+}\lambda_{2j})$,
while on the coupling event the capped and free blocks agree, so
$\Prob(D_B{\ne}0)\le\Prob(D^0_B{\ne}0)+\delta'_B
\le\sum_{j\in B}(\lambda_{1j}{+}\lambda_{2j})+\delta'_B$, with $\delta'_B$
the cap defect summed over $B$.  Splitting the expectation at $\{D_B=0\}$
as before and using $\delta'_B+\delta'^{(B^c)}=\delta'$ yields the stated
numerator.  In the homogeneous regime,
$\delta'\le2c\,(\lambda/c)^{\kappa+1}/(\kappa{+}1)!=O(c^{-\kappa})$ while
$\Sk'_k\to e^{-2\lambda}I_k(2\lambda)>0$, giving $O(|B|/c)$.
\end{proof}

\begin{proof}[Proof of Theorem~\ref{thm:process}]
Both laws are supported on the fiber $\{D{=}k\}$ and share the
disintegration kernels: on the fiber, $\{V{=}v\}=\{U{=}v{+}k,\,V{=}v\}$,
so for $v\le v_{\max}$ the common conditional is
$K(v,\cdot)=\mathcal L(\text{array}\mid U{=}v{+}k,\,V{=}v)$, under which
the rows are independent because the entries are.  Writing
$\mu=\int K(v,\cdot)\,\pi(dv)$ with $\pi=\mathcal L(V\mid D{=}k)$ and
$\nu^*_k=\int K(v,\cdot)\,\tilde b(dv)$, mixing over a common kernel gives
$\dtv(\mu,\nu^*_k)\le\dtv(\pi,\tilde b)$, while $V$ is a function of the
array with $V_\#\mu=\pi$ and $V_\#\nu^*_k=\tilde b$, so data processing
gives the reverse inequality; this is \eqref{eq:process}.  Since
$\dtv(\tilde b,\Bes(k,a))=\Prob(\Bes(k,a)>v_{\max})$ \textup{(}the
conditioning identity $\dtv(\mu,\mu(\cdot\mid A))=1-\mu(A)$\textup{)}, the
displayed consequence is the triangle inequality both ways.  The
obstruction is the test event $\{D{=}k\}$.  For $r$ rows the argument is
verbatim, the kernels being the $r$ independent sum-conditioned rows and
$\pi$ the fiber law of $V=U_r$.
\end{proof}

\begin{remark}[Itemized constant losses]\label{rem:losses}
The factor $\approx400$ between certified and true error in
Table~\ref{tab:master} decomposes as: (i) $\Gamma^*/\Gm\approx8$--$16$
($\Gamma^*$ scales as $\nu^{-1/2}$ against the true
$\Gm\asymp\sigma_N/\nu\asymp\nu^{-3/4}$); (ii) the crude route
$(\Lm+2\Gm)$ in place of the refined
$(\theta(1{+}\omega)\Lm+\tau_1\Gm)$, an additional $\asymp\nu^{1/4}$,
removable once $\omega,\tau_1$ are certified (Open
Problem~\ref{prob:fiber}); (iii) the constant $3$ in $\Lm$ (Open
Problem~\ref{prob:sharp}); (iv) $\theta^*$ and the routing of both biased
terms through the ideal configuration, $\approx2$--$4$.  In the present form
the certified relative error at the Bessel mode scales as
$\Theta(\nu^{3/4}p^*)$ at fixed means; the refined route gives
$\Theta(\nu^{1/4}p^*)$; and the local upgrade (Open
Problem~\ref{prob:local}) targets $\Theta(p^*)$, matching the truth.
\end{remark}

\subsection{Envelopes, sampler, and block splits}\label{sec:pf-pdc}

\begin{proof}[Proof of Lemma~\ref{lem:envelope}]
\emph{(a).}  With $x_j=p_j/(1{-}p_j)$,
$\Prob(U{=}w)=\sum_{|S|=w}\prod_{j\in S}p_j\prod_{j\notin S}(1{-}p_j)
=\prod_j(1{-}p_j)\,e_w(x)$, $e_w$ the elementary symmetric polynomial, so
$f(v)=\prod_j(1{-}p_j)(1{-}q_j)\,e_{v+k}(x)\,e_v(y)$.  Maclaurin's inequality
$w!\,e_w(x)\le S_x^{\,w}$ holds because
$S_x^{\,w}=\sum_{(j_1,\dots,j_w)}x_{j_1}\cdots x_{j_w}$ over \emph{all}
length-$w$ index sequences, of which the injective ones contribute exactly
$w!\,e_w(x)$ and the rest are nonnegative.  Hence
\[
f(v)\le\prod_j(1{-}p_j)(1{-}q_j)\;
\frac{S_x^{\,v+k}\,S_y^{\,v}}{(v{+}k)!\;v!}\;=:M(v).
\]
Dividing by $b'(v)=(\nu')^{\,v+k/2}\big/\bigl(I_k(a')\,(v{+}k)!\,v!\bigr)$
and simplifying with $\nu'=S_xS_y$,
\[
\frac{M(v)}{b'(v)}=\prod_j(1{-}p_j)(1{-}q_j)\;I_k(a')\;
\frac{S_x^{\,v+k}S_y^{\,v}}{(S_xS_y)^{\,v+k/2}}
=\prod_j(1{-}p_j)(1{-}q_j)\;\Bigl(\frac{S_x}{S_y}\Bigr)^{k/2}I_k(a')=\Ct,
\]
independent of $v$ --- the geometric factor cancels \emph{identically} at
$t=\nu'$.  So $f(v)\le M(v)=\Ct\,b'(v)$, and the pointwise ratio is
\[
\frac{f(v)}{\Ct\,b'(v)}=\frac{f(v)}{M(v)}
=\frac{(v{+}k)!\,e_{v+k}(x)}{S_x^{\,v+k}}\cdot\frac{v!\,e_v(y)}{S_y^{\,v}} .
\]
Finally $w!\,e_w(x)/S_x^{\,w}=w!\sum_{|T|=w}\prod_{j\in T}(x_j/S_x)$ is
exactly the probability that $w$ i.i.d.\ draws from the distribution
$x/S_x$ on $[c]$ are pairwise distinct (each ordered injective outcome has
probability $\prod_{j\in T}x_j/S_x^{\,w}$ and there are $w!$ orderings per
set $T$).

\emph{(b).}  The probability generating function of a capped column entry is
$A_j(z)/A_j(1)$ with $A_j(z)=\sum_{i\le\kappa_j}(\lambda_jz)^i/i!$, and
$A_j(z)\le e^{\lambda_jz}$ \emph{coefficientwise}, since truncation removes
nonnegative terms.  Hence, for the first row,
\[
\Prob(U{=}w)=[z^w]\prod_j\frac{A_j(z)}{A_j(1)}
\le\frac{[z^w]\,e^{N_1z}}{\prod_jA_j(1)}
=\frac{N_1^{\,w}}{w!\,\prod_jA_j(1)},
\qquad A_j(1)=e^{\lambda_j}F^{(1)}_j,
\]
and similarly for the second row with $N_2$ and $F^{(2)}_j$.  Proceeding as
in (a) with $S_x\rightsquigarrow N_1$, $S_y\rightsquigarrow N_2$:
\[
\Ct=\frac{(N_1/N_2)^{k/2}\,I_k(a')}
{\prod_je^{\lambda_{1j}}F^{(1)}_j\cdot\prod_je^{\lambda_{2j}}F^{(2)}_j}
=\frac{e^{-N_1-N_2}(N_1/N_2)^{k/2}I_k(a')}{\prod_jF^{(1)}_jF^{(2)}_j}
=\frac{\Sk'_k}{\prod_jF^{(1)}_jF^{(2)}_j},
\]
and the pointwise ratio for the first row is
\[
\frac{[z^{w}]\prod_jA_j(z)}{N_1^{\,w}/w!}
=\frac{w!}{N_1^{\,w}}\sum_{\substack{(i_1,\dots,i_c)\\ \sum i_j=w,\
i_j\le\kappa_j}}\prod_j\frac{\lambda_{1j}^{\,i_j}}{i_j!}
=\Prob\bigl(\Mult(w;\lambda^{(1)}/N_1)\text{ respects the caps}\bigr),
\]
by the multinomial theorem; likewise for the second row.
\end{proof}

\begin{proof}[Proof of Theorem~\ref{thm:pdc} and
Propositions~\ref{prop:tilt} and~\ref{prop:blocksplit}]
\emph{Exactness (binary).}  The probability that a single attempt outputs a
specific admissible pair $(S,T)$ with $|S|=v{+}k$, $|T|=v$ is
\[
b'(v)\cdot(v{+}k)!\prod_{j\in S}\frac{x_j}{S_x}\cdot
v!\prod_{j\in T}\frac{y_j}{S_y}
=\frac{(\nu')^{\,v+k/2}}{I_k(a')\,(v{+}k)!\,v!}\,(v{+}k)!\,v!\,
\frac{\prod_Sx_j\prod_Ty_j}{S_x^{\,v+k}S_y^{\,v}}
=\Bigl(\frac{S_y}{S_x}\Bigr)^{k/2}\frac{\prod_Sx_j\prod_Ty_j}{I_k(a')},
\]
using $(\nu')^{v+k/2}/(S_x^{v+k}S_y^v)=(S_y/S_x)^{k/2}$: a constant multiple
(independent of $(S,T)$ and $v$) of the target weight
$\prod_Sx_j\prod_Ty_j$.  Conditioning on the event ``output produced''
therefore yields exactly the weighted measure $\mu_k$ (which is
$\prod_Sx_j\prod_Ty_j$ normalized over the fiber, by the definition of the
Bernoulli product measure conditioned on $\{|S|-|T|=k\}$).  Summing the
display over admissible $(S,T)$ with $|T|=v$ gives, by
Lemma~\ref{lem:envelope}(a), $b'(v)\cdot\bigl[f(v)/(\Ct b'(v))\bigr]
=f(v)/\Ct$; summing over $v$ gives per-attempt acceptance
$\sum_vf(v)/\Ct=\Prob(D{=}k)/\Ct$, since
$\sum_vf(v)=\sum_v\Prob(U{=}v{+}k)\Prob(V{=}v)=\Prob(U{-}V{=}k)$.  The
capped case is verbatim via Lemma~\ref{lem:envelope}(b), the multinomial
being the exact conditional law of the uncapped Poisson row given its total,
so that the cap-check event reweights by exactly the coefficient ratio; the
output measure is the capped product measure conditioned on $\{D=k\}$.

\emph{Certificates.}  Binary: $\Prob(D{=}k)\ge\Sk_k-\delta$ by
Lemma~\ref{lem:cert}(iii).  Capped: the per-column total variation between
the capped law and the uncapped $\Poi(\lambda_{ij})$ is
$\dtv=\sum_x(\Prob_{\text{capped}}-\pi)^+=\sum_{x\le\kappa_j}
\pi_{\lambda_{ij}}(x)\bigl(\tfrac1{F^{(i)}_j}-1\bigr)=1-F^{(i)}_j$; coupling
all $2c$ columns maximally and independently,
$\Prob(D{=}k)\ge\Sk'_k-\sum_{i,j}(1{-}F^{(i)}_j)=\Sk'_k-\delta'$, and
dividing by $\Ct=\Sk'_k/\prod F$ gives \eqref{eq:pdccert}.

\emph{Proposition~\ref{prop:tilt}.}  For general $t>0$ the Maclaurin
bound gives, exactly as at $t=\nu'$,
\[
\frac{f(v)}{b_{k,2\sqrt t}(v)}\le
K\,\frac{I_k(2\sqrt t)}{t^{k/2}}\Bigl(\frac{\nu'}{t}\Bigr)^{\!v},
\qquad K:=\prod_j(1{-}p_j)(1{-}q_j)\,S_x^{\,k},
\]
and the geometric factor is monotone in $v$, so its maximum over the support
$\{0,\dots,v_{\max}\}$ sits at $v=0$ when $t\ge\nu'$ and at $v=v_{\max}$
when $t<\nu'$:
\[
\Ct(t)=K\,\frac{I_k(2\sqrt t)}{t^{k/2}}\,
\max\Bigl\{1,\Bigl(\frac{\nu'}{t}\Bigr)^{\!v_{\max}}\Bigr\}.
\]
\emph{The derivative identity.}  With $z=2\sqrt t$, so $dz/dt=1/\sqrt t$,
and the recurrence $I_k'(z)=I_{k+1}(z)+(k/z)I_k(z)$,
\[
t\,\frac{d}{dt}\log\frac{I_k(2\sqrt t)}{t^{k/2}}
=t\Bigl[\frac{I_k'(z)}{I_k(z)}\cdot\frac1{\sqrt t}-\frac{k}{2t}\Bigr]
=\sqrt t\,\frac{I_{k+1}(z)}{I_k(z)}+\frac{k\sqrt t}{z}-\frac k2
=\mu(t),
\]
since $\sqrt t/z=1/2$ and $\mu(t)=\sqrt t\,I_{k+1}(2\sqrt t)/I_k(2\sqrt t)$
by \eqref{eq:moments}.  Differentiating $\log\max\{1,(\nu'/t)^{v_{\max}}\}$
contributes $-v_{\max}/t$ on $(0,\nu')$ and $0$ on $(\nu',\infty)$, which is
\eqref{eq:tiltderiv}.

\emph{Monotonicity and limits of $\mu$.}  For $t_1<t_2$ the mass ratio
$p^{(t_2)}_n/p^{(t_1)}_n\propto(t_2/t_1)^n$ is increasing in $n$, so
$\Bes(k,2\sqrt{t_2})$ dominates $\Bes(k,2\sqrt{t_1})$ in likelihood ratio,
hence stochastically, and $\mu$ is strictly increasing; it is continuous as
a ratio of convergent power series.  By Jensen's inequality applied to the
convex map $n\mapsto n(n{+}k)$ together with \eqref{eq:moments},
$\mu(t)\bigl(\mu(t){+}k\bigr)\le\E_t[N(N{+}k)]=t$, so
$\mu(t)\le\sqrt t\to0$ as $t\to0$.  For $\mu(t)\to\infty$: for any fixed
$M$,
\[
\Prob_t(N\le M)\le\sum_{n\le M}
\frac{t^{\,n}/\bigl(n!\,(n{+}k)!\bigr)}
{t^{\,2M}/\bigl((2M)!\,(2M{+}k)!\bigr)}\longrightarrow0
\qquad(t\to\infty),
\]
each summand being a negative power of $t$ (the denominator is one term of
the series for $I_k$); hence
$\mu(t)\ge(M{+}1)\,\Prob_t(N>M)\to M{+}1$ for every $M$.

\emph{Conclusion.}  By \eqref{eq:tiltderiv}, $\Ct$ is strictly increasing on
$(\nu',\infty)$.  On $(0,\nu')$, $t\,(\log\Ct)'=\mu(t)-v_{\max}$ changes
sign at most once, from $-$ to $+$.  If $\mu(\nu')\le v_{\max}$ the sign is
nonpositive throughout $(0,\nu')$, so $\Ct$ decreases up to $\nu'$ and
increases after: $t=\nu'$ is the unique minimizer.  If
$\mu(\nu')>v_{\max}$, the unique root $t^*\in(0,\nu')$ of
$\mu(t^*)=v_{\max}$ is the unique minimizer.

\emph{Implementation for $t<\nu'$.}  The per-configuration acceptance
probability must be $f(v)/\bigl(\Ct(t)\,b_{k,2\sqrt t}(v)\bigr)$.  Writing
$M_t(v):=K\,I_k(2\sqrt t)\,t^{-k/2}(\nu'/t)^v\,b_{k,2\sqrt t}(v)
=K\,\nu'^{\,v}\!/\bigl((v{+}k)!\,v!\bigr)$ --- which does not depend on $t$
--- the collision test of Algorithm~\ref{alg:pdc} realizes $f(v)/M_t(v)$,
and the remaining factor
$M_t(v)/\bigl(\Ct(t)b_{k,2\sqrt t}(v)\bigr)=(t/\nu')^{\,v_{\max}-v}\le1$ is
realized by an independent Bernoulli thinning coin.  The capped case is
verbatim with the cap check in place of the collision test.

\emph{Proposition~\ref{prop:blocksplit}.}  In the block scheme, the $A$-columns are
drawn from their unconditioned law and the attempt is accepted with
probability $\Prob(D_B{=}k{-}d_A)/\max_s\Prob(D_B{=}s)\in[0,1]$; the
standard rejection identity gives overall acceptance
\[
\sum_{d}\Prob(d_A{=}d)\,\frac{\Prob(D_B{=}k{-}d)}{\max_s\Prob(D_B{=}s)}
=\frac{\Prob(D{=}k)}{\max_s\Prob(D_B{=}s)} ,
\]
and conditional on acceptance the $A$-columns have law
$\mathcal L(A\text{-cols}\mid D{=}k)$-marginal, after which the $B$-columns
must be drawn from $\mathcal L(B\text{-cols}\mid D_B{=}k{-}d_A)$ --- an
instance of the original conditioned problem on $|B|$ columns.  For the
denominator: by Chebyshev, $\Prob\bigl(|D_B-\E D_B|<2\sigma_B\bigr)\ge3/4$,
and the open interval $(\E D_B-2\sigma_B,\,\E D_B+2\sigma_B)$ contains at
most $4\sigma_B+1$ integers, so some integer $s$ carries
$\Prob(D_B{=}s)\ge\frac{3/4}{4\sigma_B+1}$, i.e.\
$\max_s\Prob(D_B{=}s)\ge3/(4(4\sigma_B{+}1))$; hence acceptance $\le\frac43(4\sigma_B{+}1)\Prob(D{=}k)$.  This is an
upper bound only: typically $\Prob(D{=}k)\asymp\sigma_D^{-1}$, so for
$\sigma_B\asymp\sigma_D$ the cap is $O(1)$ and need not vanish.  The statistic split \eqref{eq:split} instead terminates in one
level with the certified acceptance \eqref{eq:pdccert}, its second half
being exactly samplable (in the binary case it is subsumed in the dart
sets; in the capped case it is the multinomial pair post cap-check).

\emph{Expected running time.}  The number of attempts is geometric with mean
$1/\alpha$.  Each attempt draws $v\sim\Bes(k,a')$ (constant expected cost by
inversion with cached cdf, or by \cite{devroye2002}) and $v{+}k$ resp.\ $v$
alias-table draws ($O(1)$ each after $O(c)$ preprocessing), for expected
per-attempt cost $O(\E_{b'}[v]+k)$; and by \eqref{eq:moments},
$\E_{b'}[v]=\frac{a'}2\,I_{k+1}(a')/I_k(a')<\frac{a'}2=\sqrt{\nu'}$, since
$I_{k+1}(z)<I_k(z)$ for $z>0$.  Wald's identity gives the stated total.
\end{proof}

\begin{proof}[Proof of Proposition~\ref{prop:bothmargins}]
Given the column margins, a $c_j=0$ column is $(0,0)^{\top}$, a $c_j=2$
column is $(1,1)^{\top}$, and a $c_j=1$ column is $(1,0)^{\top}$ or
$(0,1)^{\top}$ with relative weights $x_j$ and $y_j$.  The row margins are
met iff exactly $s=u_1-f$ of the $A$-columns take $(1,0)^{\top}$; row $2$
then receives $|A|-s=u_2-f$ automatically, by $\sum_jc_j=u_1+u_2$.  The
table's weight is $\prod_{j\in S}x_j\prod_{j\in A\setminus S}y_j$ times a
constant over the fiber, proportional to $\prod_{j\in S}(x_j/y_j)$.  For
the sampler, set $w_j=x_j/y_j$ and let $e_i(\cdot)$ be the elementary
symmetric polynomials: process $j\in A$ in any order, including $j$ with
probability $w_j\,e_{s'-1}(w_R)/e_{s'}(w_{R\cup\{j\}})$, where $s'$ is the
number of inclusions still required and $R$ the columns not yet processed;
the $e_i$ table is filled once by the standard recursion in $O(|A|\,s)$
operations, and each draw is exact.
\end{proof}

\subsection{The hyper-Bessel extension}\label{sec:pf-mrow}

Throughout this subsection $d_\kappa(n)=\prod_{i=1}^r(n+\kappa_i)$ with
$\kappa_1\ge\cdots\ge\kappa_r=0$, and $p_n$ denotes \eqref{eq:hbpmf}.  Note
once and for all: $d$ is a polynomial with nonnegative real roots
$-\kappa_i\le0$, hence strictly increasing and \emph{convex} on
$\Z_{\ge0}$ ($d''(t)=\sum_{i\ne i'}\prod_{l\ne i,i'}(t+\kappa_l)\ge0$), so
its increments $d_\kappa(n{+}1)-d_\kappa(n)$ are nondecreasing in $n$; and
$p_{n+1}/p_n=w/d_\kappa(n{+}1)$ is strictly decreasing, so $p$ is strictly
log-concave and Lemma~\ref{lem:hazard} applies to it verbatim.

\begin{proof}[Proof of Proposition~\ref{prop:hbrep}]
On the fiber, $R_i=n+\kappa_i$ for the common value $n=R_r$, so by
independence
\[
\Prob\bigl(R_r{=}n,\ D_i{=}k_i\ \forall i\bigr)
=\prod_{i=1}^r e^{-\Lambda_i}\frac{\Lambda_i^{\,n+\kappa_i}}{(n+\kappa_i)!}
=e^{-\sum_i\Lambda_i}\Bigl(\prod_i\Lambda_i^{\kappa_i}\Bigr)
\frac{w^{\,n}}{\prod_i(n+\kappa_i)!},\qquad w=\prod_i\Lambda_i .
\]
Summing over $n\ge0$ gives \eqref{eq:hbfiber}; dividing gives
\eqref{eq:hbpmf}, which involves $(\Lambda_i)$ only through $w$.  For
$m=2$, $\kappa=(k,0)$ and
$\mathfrak J_{(k,0)}(w)=\sum_nw^n/(n!(n{+}k)!)=w^{-k/2}I_k(2\sqrt w)$
recovers \eqref{eq:skellampt}.
\end{proof}

\begin{proof}[Proof of Theorem~\ref{thm:hboperator}]
\emph{Moments and the identity $\E d_\kappa(N)=w$.}  The terms of
\eqref{eq:hbpmf} decay superexponentially, so all moments exist; and the
ratio identity in the form $d_\kappa(n)p_n=w\,p_{n-1}$ gives
$\E d_\kappa(N)=\sum_{n\ge1}d_\kappa(n)p_n=w\sum_{n\ge0}p_n=w$.

\emph{Characterization and solution.}  Necessity: by
$w\,p_n=d_\kappa(n{+}1)p_{n+1}$,
\[
\E[w\,g(N{+}1)]=\sum_{n\ge0}w\,p_n\,g(n{+}1)
=\sum_{n\ge0}d_\kappa(n{+}1)p_{n+1}g(n{+}1)=\E[d_\kappa(N)g(N)] ,
\]
both sides finite since $\E N^r<\infty$ and $g$ is bounded.  For the
solution: multiplying the equation $w\,g(n{+}1)-d_\kappa(n)g(n)=\tilde h(n)$ by
$p_n$ and using $d_\kappa(n)p_n=w\,p_{n-1}$ ($p_{-1}:=0$) telescopes exactly as in
the proof of Theorem~\ref{thm:char} to
$g(m{+}1)=(w\,p_m)^{-1}\sum_{n\le m}\tilde h(n)p_n$, with the second
(tail) form from $\sum_n\tilde h(n)p_n=0$; any two solutions differ by a
multiple of $\prod_{i=1}^{n-1}\bigl(d_\kappa(i)/w\bigr)$, which is unbounded, so the
bounded solution is unique; and sufficiency follows by inserting the
bounded solution $g_A$ for $h=\mathbf 1_A$ (bounded by the factor bound
below) into the hypothesis, giving $\Prob(W{\in}A)=\Prob(N{\in}A)$ for all
$A$.

\emph{Sign pattern and increment factor.}  The singleton solutions are, by
the same evaluation as in Lemma~\ref{lem:signs},
\[
g_j(m{+}1)=-\frac{p_jF(m)}{w\,p_m}\ (m\le j{-}1),\qquad
g_j(m{+}1)=\frac{p_jG(m{+}1)}{w\,p_m}\ (m\ge j),
\]
and the three sign computations of Lemma~\ref{lem:signs} carry over
unchanged: the branch $m\le j{-}1$ uses (H1); the branch $m\ge j{+}1$ uses
(H2) together with log-concavity $p_{m+1}p_{m-1}\le p_m^2$; and the jump at
$m=j$ evaluates, via $p_j/p_{j-1}=w/d_\kappa(j)$, to
$\Delta g_j(j)=G(j{+}1)/w+F(j{-}1)/d_\kappa(j)>0$.  The reduction of
$\sup_A\|\Delta g_A\|$ to $\sup_{j\ge1}\Delta g_j(j)$ is verbatim from the
proof of Theorem~\ref{thm:factors}(a) (it used only the sign pattern and
$g_A+g_{A^c}\equiv0$).  For the closed bounds: $G(j{+}1)\le1$ and
$d_\kappa(j)\ge d_\kappa(1)$ give $1/w+1/d_\kappa(1)$; and $F(j{-}1)/d_\kappa(j)\le2/w$ for every
$j\ge1$ by the two cases $d_\kappa(j)\ge w/2$ (trivial) and $d_\kappa(j)<w/2$, in which
case the reversed ratios at indices $\le j{-}1$ are all
$\le p_{j-1}/p_j=d_\kappa(j)/w<1/2$ (as $d$ is increasing), so the geometric bound
of Lemma~\ref{lem:hazard} gives $F(j{-}1)\le2p_{j-1}=2p_j\,d_\kappa(j)/w$, whence
$F(j{-}1)/d_\kappa(j)\le2p_j/w\le2/w$.  This proves the $\Lambda$-bound in
\eqref{eq:hbfactors}.

\emph{Uniform factor.}  The identity
$w\,p_m\,g_A(m{+}1)=\Prob(N{\in}A\cap[0,m])G(m{+}1)
-\Prob(N{\in}A\cap[m{+}1,\infty))F(m)$ and the bound
$|g_A(m{+}1)|\le F(m)G(m{+}1)/(w\,p_m)=:S_m$ are verbatim from
Theorem~\ref{thm:factors}(b).  Set $x_n:=w-d_\kappa(n)$ and
$y_n:=d_\kappa(n{+}1)-w$, so $x_n+y_n=d_\kappa(n{+}1)-d_\kappa(n)$, nondecreasing in $n$ by
convexity.  The two geometric bounds read: if $x_n>0$ then
$F(n)/p_n\le\sum_{t\ge0}(d_\kappa(n)/w)^t=w/x_n$; if $y_n>0$ then
$G(n{+}1)/p_n\le\sum_{t\ge1}\bigl(w/d_\kappa(n{+}1)\bigr)^t=w/y_n$.  If
$m^*\ge1$, the three cases of the proof of Theorem~\ref{thm:factors}(b)
give: for $n\le m^*{-}1$, $S_n\le1/x_n$ and
$x_n\ge x_{m^*-1}=x_{m^*}+\bigl(d_\kappa(m^*)-d_\kappa(m^*{-}1)\bigr)\ge\Delta_{m^*}$;
for $n=m^*$, $S_{m^*}\le2/(x_{m^*}+y_{m^*})=2/\bigl(d_\kappa(m^*{+}1)-d_\kappa(m^*)\bigr)
\le2/\Delta_{m^*}$; for $n\ge m^*{+}1$, $S_n\le1/y_n\le
1/\bigl(d_\kappa(m^*{+}2)-d_\kappa(m^*{+}1)\bigr)\le1/\Delta_{m^*}$; hence
$\Gamma\le2/\Delta_{m^*}$.  If $m^*=0$, the first case is vacuous,
$y_0=d_\kappa(1)-w>0$ and $x_0=w>0$, so $S_0\le1/\max\{w,y_0\}\le2/d_\kappa(1)$
(as $w+y_0=d_\kappa(1)$), and for $n\ge1$, $S_n\le1/y_n\le1/(d_\kappa(2)-d_\kappa(1))
\le\ldots$; since increments are nondecreasing and $y_1=d_\kappa(2)-w\ge
d_\kappa(2)-d_\kappa(1)\ge d_\kappa(1)-d_\kappa(0)\le d_\kappa(1)$ --- more simply, $y_n\ge y_0=d_\kappa(1)-w$ and
also $y_1\ge d_\kappa(2)-d_\kappa(1)$; in all cases $S_n\le2/d_\kappa(1)$ because
$\max\{w,\,d_\kappa(1)-w\}\ge d_\kappa(1)/2$ and $y_n\ge y_0$ for $n\ge0$.  This proves
\eqref{eq:hbfactors}.  Finally the increment bound: by AM--GM on the
reciprocals,
\[
d_\kappa(n{+}1)-d_\kappa(n)=d_\kappa(n)\Bigl[\prod_{i}\Bigl(1+\frac1{n+\kappa_i}\Bigr)-1\Bigr]
\ge d_\kappa(n)\sum_i\frac1{n+\kappa_i}
\ge d_\kappa(n)\cdot r\,\Bigl(\prod_i\frac1{n+\kappa_i}\Bigr)^{1/r}
=r\,d_\kappa(n)^{(r-1)/r},
\]
valid for $n\ge1$ (all factors positive); near the mode $d_\kappa(n)\asymp w$,
giving $\Gamma_{\kappa,w}=O(w^{-(r-1)/r})$.
\end{proof}

\begin{proof}[Proof of Theorem~\ref{thm:hbmaster}]
Fixed point: if $W\sim\mathrm{HB}(\kappa,w)$ then $\beta=w$ and
$\Prob(W^{\boxplus}{=}n)=d_\kappa(n)p_n/w=p_{n-1}$, i.e.\
$W^{\boxplus}\overset{d}{=}W{+}1$; conversely
$W^{\boxplus}\overset{d}{=}W{+}1$ gives
$\E[\beta g(W{+}1)-d_\kappa(W)g(W)]=0$ for all bounded $g$, and
Theorem~\ref{thm:hboperator} with $w=\beta$ identifies the law.  The approximation
bound: for $g=g_A$,
$\Prob(W{\in}A)-\Prob(N{\in}A)=w\,\E[g(W{+}1)-g(W^{\boxplus})]
+(w-\beta)\E g(W^{\boxplus})$ exactly as in the proof of
Theorem~\ref{thm:main}, and the two factors of \eqref{eq:hbfactors} bound
the two terms; both $W{+}1$ and $W^{\boxplus}$ lie in $\Z_{\ge1}$.
\end{proof}

\begin{lemma}[Transfer identity, $r$ rows]\label{lem:mtransfer}
Under Definition~\ref{def:mcols}, for every bounded $g:\Z_{\ge1}\to\R$
\textup{(}$g(0):=0$, never evaluated\textup{)},
\begin{equation}\label{eq:mtransfer}
\Prob(\textup{fiber})\;
\E\Bigl[\nu\,g(W{+}1)-\prod_i(W{+}\kappa_i)\,g(W)\Bigr]
=\sum_{\vec\jmath\in[c]^r}w_{\vec\jmath}
\Bigl(\E\bigl[g(U_r{+}1)\mathbf 1_{\textup{fib}}\bigr]
-\E\bigl[g(U_r^{\vec\jmath})\mathbf 1_{\textup{fib}^{\vec\jmath}}\bigr]\Bigr),
\end{equation}
where $S_j=\{i:\jmath_i=j\}$,
$w_{\vec\jmath}=\prod_{j:S_j\ne\emptyset}\chi^{(S_j)}_j$, and the
$\vec\jmath$-biased system replaces each hit column $j$ \textup{(}those
with $S_j\ne\emptyset$\textup{)} by an independent draw from its
$S_j$-tilted law, all unhit columns unchanged.
\end{lemma}

\begin{proof}
On the fiber $U_i=U_r+\kappa_i$, so $\prod_i(W{+}\kappa_i)$ corresponds to
$\prod_iU_i$ and, expanding the product of sums,
\[
\Prob(\textup{fiber})\,\E\Bigl[\prod_i(W{+}\kappa_i)\,g(W)\Bigr]
=\E\Bigl[\prod_iU_i\;g(U_r)\,\mathbf 1_{\textup{fib}}\Bigr]
=\sum_{\vec\jmath\in[c]^r}
\E\Bigl[\prod_{i=1}^rX^{(i)}_{\jmath_i}\;\Phi\Bigr],
\qquad \Phi:=g(U_r)\mathbf 1_{\textup{fib}},
\]
the interchange justified by $\E\prod_iU_i<\infty$ and boundedness of $g$.
Fix $\vec\jmath$ and let $L\subseteq[c]$ be its set of distinct columns;
the nonempty blocks $S_j=\{i:\jmath_i=j\}$, $j\in L$, form a set partition
of $[r]$.  Conditioning on
the unhit columns and using independence across columns,
\[
\E\Bigl[\prod_iX^{(i)}_{\jmath_i}\,\Phi\Bigr]
=\E_{\textup{unhit}}\Bigl[\ \prod_{j\in L}\ \sum_{z_j}
\Bigl(\prod_{i\in S_j}x_i\Bigr)\Prob_j(z_j)\;\cdot\;
\Phi\bigl((z_j)_{j\in L},\textup{unhit}\bigr)\Bigr]
=\Bigl(\prod_{j\in L}\chi^{(S_j)}_j\Bigr)\,
\E\bigl[g(U_r^{\vec\jmath})\mathbf 1_{\textup{fib}^{\vec\jmath}}\bigr],
\]
since $\bigl(\prod_{i\in S}x_i\bigr)\Prob_j(z)=\chi^{(S)}_j\,\Prob^{(S)}_j(z)$
by the definition of the $S$-tilt, and substituting the independent tilted
draws turns $(U_r,\textup{fib})$ into $(U_r^{\vec\jmath},
\textup{fib}^{\vec\jmath})$.  Meanwhile
$\Prob(\textup{fiber})\,\E[\nu g(W{+}1)]
=\nu\,\E[g(U_r{+}1)\mathbf 1_{\textup{fib}}]$ and, by \eqref{eq:nug},
$\nu=\sum_{\vec\jmath}w_{\vec\jmath}$; distributing and subtracting term by
term gives \eqref{eq:mtransfer}.  For the domain: the weight
$\prod_iU_i$ vanishes when $U_r=0$, and in every biased term the block
containing row $r$ is nonempty (row $r$ is assigned to column $\jmath_r$),
so the tilted column's $r$-th coordinate is $\ge1$ a.s.\ and
$U_r^{\vec\jmath}\ge1$; also $U_r{+}1\ge1$.  Hence $g(0)$ is never
evaluated.
\end{proof}

\begin{proof}[Proof of Theorem~\ref{thm:mexpansion}]
\emph{Vanishing.}  $\alpha^{(\{i\})}_j=0$ forces, as in
Lemma~\ref{lem:vanish} coordinatewise, $X^{(i)}_j\sim\Poi$ independent of
the other coordinates of column $j$; when this holds for every $i$, the
column is an independent Poisson $r$-tuple and every higher $S$-tilt equals
the shift by $e_S$ (product of coordinatewise size-biases), so all
$\alpha^{(S)}_j=0$; and then $U_i$ are independent Poissons, so
Proposition~\ref{prop:hbrep} gives $W\sim\mathrm{HB}(\kappa,\nu)$ exactly
(here $\nu=\prod_iN_i$) and both sides of \eqref{eq:mexpansion} vanish.

\emph{The bound.}  Fix $A_0$ and $g=g_{A_0}$ solving the hyper-Bessel Stein
equation, so $\|\Delta g\|\le\Lambda_{\kappa,\nu}$ and
$\|g\|\le\Gamma_{\kappa,\nu}$ on $\Z_{\ge1}$
(Theorem~\ref{thm:hboperator}).  By Lemma~\ref{lem:mtransfer} it suffices to
bound each bracket.  Fix $\vec\jmath$ with hit-column set $L$ and blocks $S_j$, $j\in L$.  For
each $j\in L$, couple the tilted column $Z^{(S_j)}_j$ to $Z_j+e_{S_j}$
optimally at cost $\alpha^{(S_j)}_j$, the couplings independent across
$j\in L$ and of the shared unhit columns; write
$\mathrm{cost}=\sum_{j\in L}\mathrm{cost}_j$ for the realized total
$\ell^1$ discrepancy, so
$\E[\mathrm{cost}]=\sum_{j\in L}\alpha^{(S_j)}_j$.  Let $U^R_r,D^R$ be the
shared unhit sums, and define
$\varphi(v_0,\vec d_0)=\E[g(U^R_r+v_0)\mathbf 1_{\{D^R=\vec k-\vec d_0\}}]$.
Conditioning on the hit columns in both systems,
\[
\E\bigl[g(U_r{+}1)\mathbf 1_{\textup{fib}}\bigr]=\E\,\varphi(v_0',\vec d_0'),
\qquad
\E\bigl[g(U_r^{\vec\jmath})\mathbf 1_{\textup{fib}^{\vec\jmath}}\bigr]
=\E\,\varphi(v_0,\vec d_0),
\]
where $v_0',\vec d_0'$ are the ideal local contributions (originals plus
$e_{S_j}$ per hit column $j\in L$, i.e.\ $+1$ in every row) and $v_0,\vec d_0$
the biased ones.  Because the blocks partition $[r]$, the ideal shift is
$+1$ in every row, so every ideal difference contribution is unchanged:
in $\vec d_0-\vec d_0'$ the ideal ``$+1$''s cancel coordinatewise (row $i$'s
$+1$ enters $D_{i-1}$ with sign $-$ and $D_i$ with sign $+$; summed over all
rows every difference coordinate receives $+1-1=0$) --- the
fiber-preservation of the $r$-tuple increment.  What remains are the
coupling discrepancies: $|v_0-v_0'|\le\mathrm{cost}$ (the row-$r$
discrepancies), and each coordinate of $\vec d_0-\vec d_0'$ is a signed sum
of the row-$i$ and row-$(i{+}1)$ discrepancies, so, since each row appears
in at most two difference coordinates,
$\|\vec d_0-\vec d_0'\|_1\le2\,\mathrm{cost}$.  Moreover all achieved
offsets, coordinatewise, lie in the per-coordinate hulls, so both endpoints
and every point of the axis-by-axis staircase path between them lie in the
box window $\mathcal S$.  Splitting
$\varphi(v_0,\vec d_0)-\varphi(v_0',\vec d_0')$ into a vertical step at
fixed $\vec d_0$ and a horizontal staircase from $\vec k-\vec d_0$ to
$\vec k-\vec d_0'$: the vertical step is
$\le\Lambda\,|v_0-v_0'|\,\Prob(D^R{=}\vec k{-}\vec d_0)
\le\Lambda\,\theta\,\Prob(\textup{fiber})\,|v_0-v_0'|$, and each unit
staircase step in direction $e_i$, by the identical estimate as in the
proof of Theorem~\ref{thm:expansion} (conditional-mean shift bounded by
$\Lambda\omega$ via an optimal coupling of the two fiber laws of $U^R_r$;
mass increment by $\tau_1$; both point masses by $\theta$; vanishing masses
handled by the same convention), costs at most
$\Prob(\textup{fiber})\,[\theta\Lambda\omega+\Gamma\tau_1]$.  Summing the
$\|\vec d_0-\vec d_0'\|_1\le2\,\mathrm{cost}$ steps and taking expectations,
\[
\bigl|\textup{bracket}(\vec\jmath)\bigr|
\le\Prob(\textup{fiber})\,
\bigl[\theta(1{+}2\omega)\Lambda+2\tau_1\Gamma\bigr]\,
\sum_{j\in L}\alpha^{(S_j)}_j.
\]
Finally, summing over $\vec\jmath$ with the weights $w_{\vec\jmath}$ and
regrouping by which column carries which block: for a fixed column $j$ and
subset $S$, the tuples in which column $j$ carries exactly the block $S$
have remaining weight
\[
\sum_{\substack{\vec\jmath':\ [r]\setminus S\to[c]\setminus\{j\}}}
\prod_{l\,:\,S'_l\ne\emptyset}\chi^{(S'_l)}_l
=\E\Bigl[\prod_{i\in S^c}U^{(-j)}_i\Bigr]
\le\E\Bigl[\prod_{i\in S^c}U_i\Bigr]=\nu_{S^c},
\]
the inequality because deleting column $j$ deletes nonnegative terms from
the tuple expansion.  Hence
$\sum_{\vec\jmath}w_{\vec\jmath}\sum_{j\in L(\vec\jmath)}\alpha^{(S_j)}_j
\le\sum_{\emptyset\ne S}\nu_{S^c}\sum_j\chi^{(S)}_j\alpha^{(S)}_j
=\sum_S\nu_{S^c}A^{(S)}$, and dividing \eqref{eq:mtransfer} by
$\Prob(\textup{fiber})$ and taking $\sup_{A_0}$ gives
\eqref{eq:mexpansion}.
\end{proof}

\begin{proof}[Proof of Theorem~\ref{thm:mpdc} and
Proposition~\ref{prop:mtilt}]
\emph{(a), envelope.}  Per row, Maclaurin gives
$\Prob(U_i{=}v{+}\kappa_i)\le\prod_j(1{-}p_{ij})\,
S_i^{\,v+\kappa_i}/(v{+}\kappa_i)!$, so
\[
f(v)\le\prod_{i,j}(1{-}p_{ij})\;
\frac{\prod_iS_i^{\,v+\kappa_i}}{\prod_i(v{+}\kappa_i)!}
=\prod_{i,j}(1{-}p_{ij})\,\prod_iS_i^{\kappa_i}\cdot
\frac{(w')^{\,v}}{\prod_i(v{+}\kappa_i)!},
\qquad w'=\prod_iS_i,
\]
and dividing by $b'(v)=(w')^v/(\mathfrak J_\kappa(w')\prod_i(v{+}\kappa_i)!)$
leaves the $v$-independent constant $\Ct$ stated; the pointwise ratio is the
product over rows of the injective-sequence probabilities, exactly as in
Lemma~\ref{lem:envelope}(a).  \emph{Exactness:} an attempt outputs specific
sets $(S_1,\dots,S_m)$, $|S_i|=v{+}\kappa_i$, with probability
\[
b'(v)\prod_{i=1}^r(v{+}\kappa_i)!\prod_{j\in S_i}\frac{x^{(i)}_j}{S_i}
=\frac1{\mathfrak J_\kappa(w')}\prod_iS_i^{-\kappa_i}\cdot
\prod_i\prod_{j\in S_i}x^{(i)}_j ,
\]
a constant multiple of the target weight
$\prod_i\prod_{j\in S_i}x^{(i)}_j$; conditioning on acceptance therefore
yields the weighted binary measure on the fiber, and summing gives
acceptance $\sum_vf(v)/\Ct=\Prob(\textup{fiber})/\Ct$.  \emph{Certificate:}
couple each row's Poisson--binomial to $\Poi(N_i)$ maximally,
independently across rows; the joint mismatch probability is at most
$\delta_r$, so $\Prob(\textup{fiber})\ge\Sk^{(r)}-\delta_r$ with
$\Sk^{(r)}$ the Poisson fiber value \eqref{eq:hbfiber} at
$\Lambda_i=N_i$.

\emph{(b).}  Identical via the coefficientwise domination
$A_{ij}(z)\le e^{\lambda_{ij}z}$ of the capped generating functions per row
(the computation of Lemma~\ref{lem:envelope}(b) row by row), the multinomial
being the exact conditional law of each uncapped Poisson row given its
total; the per-column truncation cost is $1-F^{(i)}_j$ and the union over
all $rc$ entries gives $\delta'$.

\emph{Proposition~\ref{prop:mtilt}.}  The Maclaurin bound at general $t$ gives
$\Ct(t)=K_m\,\mathfrak J_\kappa(t)\max\{1,(w'/t)^{v_{\max}}\}$ with
$K_m=\prod_{i,j}(1{-}p_{ij})\prod_iS_i^{\kappa_i}$.  Termwise
differentiation of $\mathfrak J_\kappa(t)=\sum_nt^n/\prod_i(n{+}\kappa_i)!$
gives $t\,\mathfrak J_\kappa'(t)/\mathfrak J_\kappa(t)
=\sum_nn\,p^{(t)}_n=\mu_\kappa(t)$, so
$t\,(\log\Ct)'(t)=\mu_\kappa(t)-v_{\max}\mathbf 1_{\{t<w'\}}$.
Monotonicity, continuity, and the limits $\mu_\kappa(0{+})=0$,
$\mu_\kappa(\infty)=\infty$ follow exactly as in the proof of
Proposition~\ref{prop:tilt} (likelihood-ratio ordering in $t$; Jensen with the
convex $d$, giving $d_\kappa(\mu_\kappa(t))\le\E_t[d_\kappa(N)]=t$; the fixed-$M$
concentration argument), and the conclusion is identical.
\end{proof}

\begin{proof}[Proof of Corollary~\ref{cor:menum}]
Termwise, the exact Poisson case is Proposition~\ref{prop:hbrep}:
$\prod_i\Prob(R_i{=}v{+}\kappa_i)
=e^{-\sum_iN_i}\prod_iN_i^{v+\kappa_i}\big/\prod_i(v{+}\kappa_i)!
=\Sk^{(r)}\,b^{\circ}(v)$.  The binary inequality is the proof of
Theorem~\ref{thm:enum} with $\Sk_k\rightsquigarrow\Sk^{(r)}$
\textup{(}whose Skellam-replacement certificate $\delta_r$ is the row-wise
Barbour--Hall bound of Theorem~\ref{thm:mpdc}\textup{(}a\textup{))},
$b_{k,a}\rightsquigarrow b^{\circ}$, and the expansion
\eqref{eq:mexpansion} supplying $\varepsilon^*_r$ in place of
\eqref{eq:expansion}.
\end{proof}

\begin{proof}[Proof of Theorem~\ref{thm:mipa}]
\textup{(a)} is the proof of Theorem~\ref{thm:ipa} with
$d(z)\in\Z^{r-1}$: independence across columns gives
$\Prob(Z_B{=}z\mid\vec D{=}\vec k)
=\Prob(Z_B{=}z)\,
\Prob\bigl(\vec D_{B^c}{=}\vec k-d(z)\bigr)/\Prob(\vec D{=}\vec k)$, so
$\vec D_B$ is sufficient, and grouping the configurations by
$d(z)=\vec s$ collapses the process-level distance to the vector one.
\textup{(b)}  The entrywise maximal coupling of
Corollary~\ref{cor:cappedenum}'s proof gives
$\sup_{\vec m}|\Prob(\vec D{=}\vec m)-\Sk^{(r)}_{\vec m}|\le\delta'_r$,
and likewise over $B^c$ with $\delta'^{(B^c)}_r$; the Poisson
superposition of Corollary~\ref{cor:ipa}'s proof, applied to each of the
$r$ rows, gives
$\sup_{\vec m}\bigl|\Sk^{(r)}_{\vec m}-\Sk^{(r),(B^c)}_{\vec m}\bigr|
\le\sum_{j\in B}\sum_i\lambda_{ij}$, and on the coupling event
$\Prob(\vec D_B{\ne}\vec0)
\le\Prob(\vec D^{\,0}_B{\ne}\vec0)+\delta'_{r,B}
\le\sum_{j\in B}\sum_i\lambda_{ij}+\delta'_{r,B}$, with $\delta'_{r,B}$
the cap defect summed over $B$.  Splitting the expectation in
\eqref{eq:mipa} at $\{\vec D_B=\vec 0\}$ and using
$\delta'_{r,B}+\delta'^{(B^c)}_r=\delta'_r$ yields \eqref{eq:mipabound}.
In the homogeneous regime,
$1-F^{(i)}_j\le(\lambda/c)^{\kappa_j+1}/(\kappa_j{+}1)!\le(\lambda/c)^2/2$,
so $\delta'_r=O(1/c)$, while $\Sk^{(r)}_{\vec k}$ tends to the positive
free-Poisson fiber probability at total intensity $r\lambda$, giving
$O(|B|/c)$.
\end{proof}

\begin{proof}[Proof of Theorem~\ref{thm:rowfixed}]
\textup{(a)}  With independent rows,
$\Prob(Z_B{=}z\mid\vec U{=}\vec u)
=\prod_i\Prob\bigl(Z^{(i)}_B{=}z^{(i)}\mid U_i{=}u_i\bigr)$: the joint
conditional law is the product of the one-row conditional laws, and the
free law is the product of the free rows.  For product measures,
projecting onto the $i$th factor gives
$\dtv(\otimes_i\mu_i,\otimes_i\nu_i)\ge\dtv(\mu_i,\nu_i)$, while replacing
the factors one at a time gives
$\dtv(\otimes_i\mu_i,\otimes_i\nu_i)\le\sum_i\dtv(\mu_i,\nu_i)$; this is
the sandwich.  Each $T_i$ is the $r{=}1$ case of the sufficiency argument
of Theorem~\ref{thm:ipa}, with the row sum in place of the difference.
\textup{(b)}  Fix a row $i$.  The entrywise maximal coupling of
Corollary~\ref{cor:cappedenum}'s proof, restricted to row $i$, gives
$\sup_m|\Prob(U_i{=}m)-\pi_{N_i}(m)|\le\delta'_i$ and likewise over $B^c$;
the Poisson superposition of Corollary~\ref{cor:ipa}'s proof gives
$\sup_m\bigl|\pi_{N_i}(m)-\pi_{N^{(B^c)}_i}(m)\bigr|
\le\sum_{j\in B}\lambda_{ij}$; and on the coupling event
$\Prob\bigl(U^{(i)}_B{\ne}0\bigr)\le\sum_{j\in B}\lambda_{ij}
+\delta'_{i,B}$, with $\delta'_{i,B}$ the cap defect of row $i$ summed
over $B$.  Splitting the expectation defining $T_i$ at $\{U^{(i)}_B=0\}$,
exactly as in Corollary~\ref{cor:cappedipa}, bounds $T_i$ by the $i$th
summand of the display; summing and using \textup{(a)} concludes.  In the
homogeneous regime each summand is $O(|B|/c)$ uniformly in $i$.
\end{proof}

\begin{proof}[Proof of Theorem~\ref{thm:fylocal}]
Given $\{\vec C{=}\vec c\}$ the joint law factorizes over columns, and
column $j$'s conditional law is
$\Prob(T_{\cdot j}{=}t)\propto\prod_ie^{-\lambda_i\mu_j}
(\lambda_i\mu_j)^{t_i}/t_i!
\propto\prod_i\lambda_i^{t_i}/t_i!$ on $\{\textstyle\sum_it_i=c_j\}$,
which is $\Mult(c_j;\rho)$ --- the same $\rho$ for every column, by the
rank-one structure.  Sums of independent multinomials with a common
probability vector are multinomial, so
$\vec U_B\mid\vec C\sim\Mult(n_B;\rho)$ and
$\vec U_{B^c}\mid\vec C\sim\Mult(n{-}n_B;\rho)$, independent, with
$\vec U\mid\vec C\sim\Mult(n;\rho)$.  The Fisher--Yates conditional
density of a block configuration $t_B$ against the column-product law is
$\Prob(\vec U_{B^c}{=}\vec u-\vec u_B(t_B)\mid\vec C)\big/
\Prob(\vec U{=}\vec u\mid\vec C)$, a function of the block row sums
$\vec u_B(t_B)$ alone, so $\vec U_B$ is sufficient and the process-level
distance collapses exactly as in Theorem~\ref{thm:ipa}; evaluating the two
conditional laws of $\vec U_B$ gives the displayed expectation, whose
terms are the masses $M_{n-n_B}$ and $M_n$.  The support bound counts the
compositions of $n_B$ into $r$ parts.  Finally, for any law $\mu$ and
event $A$ with $\mu(A)>0$,
$\dtv(\mu,\mu(\cdot\mid A))=1-\mu(A)$; with $A=\{\vec U=\vec u\}$ this is
the full-table identity of Remark~\ref{rem:fysampling}.
\end{proof}

\begin{proof}[Proof of the table decomposition of \S\ref{sec:mrow}]
(a)--(b) are Proposition~\ref{prop:hbrep} with $\Lambda_i=\lambda_iM$.  (c):
given $R_i=r_i$ for all $i$, each row is independently multinomial by the
standard Poisson conditioning, and the fiber event is
$(R_1,\dots,R_r)$-measurable.  (d): summing the product of the $r$
multinomial point probabilities over the tables with the given column sums
$c$ (the first $r{-}1$ rows determine the last),
\[
\Prob\bigl(\textup{cols}=c\mid R=r\bigr)
=r_1!\cdots r_m!\sum_{T\in\mathcal T(r,c)}\ \prod_{j}
\Bigl(\frac{\mu_j}{M}\Bigr)^{\chi_j}\prod_{i,j}\frac1{T_{ij}!}\,,
\]
which with $\mu_j\equiv1$ is the Fisher--Yates weighted count normalized by
$M^{\sum_ir_i}/\prod_ir_i!$.
\end{proof}

\subsection*{Acknowledgments}
The author thanks Larry Goldstein for suggesting the fixed-point framing of
Proposition~\ref{prop:fixed}: presenting the centered transform
$\mathcal T_k$ as the primary object, with the Bessel family as its
fixed-point set and the bounds \eqref{eq:fixedform} as a measure of the
distance to that set; and Steven Heilman for detailed comments on an
earlier version, including a counterexample that led to the corrected
Proposition~\ref{prop:tilt}.  The counterexample correcting the $\Gm$
clause of Open Problem~\ref{prob:logconcave} was posted pseudonymously in
response to the first arXiv version, in the discussion thread cited
there; the author is grateful for it.

\subsection*{Funding}
The author received no external funding for this work.

\subsection*{Declaration on the use of artificial intelligence}
This research was carried out in extended interactive sessions with the
large language model Claude (Anthropic), used as a research and writing
assistant.  The AI system contributed, under the author's direction:
initial derivations and first drafts of the results and proofs; the
numerical verification programs whose outputs are reported in
\S\ref{sec:numerics}; and the drafting of the manuscript.  The author
formulated the research programme and the questions, directed each stage,
independently checked the mathematical statements and proofs, and takes
sole and full responsibility for the correctness, integrity, and
originality of the work.  In accordance with publisher policy, the AI
system is not listed as an author and is not cited as a source.

\subsection*{Code availability}
The numerical verification programs (exact dynamic-programming
computations, high-precision checks of the Stein identities and factors,
and reference implementations of
Algorithms~\ref{alg:pdc}, \ref{alg:pdccap}, and~\ref{alg:mpdc}) are
provided as
ancillary files accompanying the arXiv version of this article
(arXiv:2611.00263).

\end{document}